\documentclass[12pt]{amsart}
\usepackage{amssymb}
\title[Arakelov inequalities and uniformization]
{Arakelov inequalities and the uniformization of certain rigid Shimura varieties}
\author[Eckart Viehweg]{Eckart Viehweg}
\address{Universit\"at Duisburg-Essen, Mathematik, 45117 Essen, Germany}
\email{viehweg@uni-essen.de}
\thanks{This work has been supported by the ``DFG-Schwerpunktprogramm
Globale Methoden in der Komplexen Geometrie'', and by the DFG-Leibniz program.}
\author[Kang Zuo]{Kang Zuo}
\address{Universit\"at Mainz,
Fachbereich 17, Mathematik,
55099 Mainz, Germany}
\email{kzuo@mathematik.uni-mainz.de}
\begin{document}
\theoremstyle{plain}
\newtheorem{thm}{Theorem}[section]
\newtheorem{theorem}[thm]{Theorem}
\newtheorem{lemma}[thm]{Lemma}
\newtheorem{corollary}[thm]{Corollary}
\newtheorem{proposition}[thm]{Proposition}
\newtheorem{addendum}[thm]{Addendum}
\newtheorem{variant}[thm]{Variant}
\newtheorem{lemdef}[thm]{Lemma and Definition}
\theoremstyle{definition}
\newtheorem{construction}[thm]{Construction}
\newtheorem{notations}[thm]{Notations}
\newtheorem{question}[thm]{Question}
\newtheorem{problem}[thm]{Problem}
\newtheorem{remark}[thm]{Remark}
\newtheorem{remarks}[thm]{Remarks}
\newtheorem{definition}[thm]{Definition}
\newtheorem{claim}[thm]{Claim}
\newtheorem{assumption}[thm]{Assumption}
\newtheorem{assumptions}[thm]{Assumptions}
\newtheorem{properties}[thm]{Properties}
\newtheorem{example}[thm]{Example}
\newtheorem{setup}[thm]{Set-up}
\numberwithin{equation}{section}
\catcode`\@=11
\def\opn#1#2{\def#1{\mathop{\kern0pt\fam0#2}\nolimits}}
\def\bold#1{{\bf #1}}%
\def\underrightarrow{\mathpalette\underrightarrow@}
\def\underrightarrow@#1#2{\vtop{\ialign{$##$\cr
 \hfil#1#2\hfil\cr\noalign{\nointerlineskip}%
 #1{-}\mkern-6mu\cleaders\hbox{$#1\mkern-2mu{-}\mkern-2mu$}\hfill
 \mkern-6mu{\to}\cr}}}
\let\underarrow\underrightarrow
\def\underleftarrow{\mathpalette\underleftarrow@}
\def\underleftarrow@#1#2{\vtop{\ialign{$##$\cr
 \hfil#1#2\hfil\cr\noalign{\nointerlineskip}#1{\leftarrow}\mkern-6mu
 \cleaders\hbox{$#1\mkern-2mu{-}\mkern-2mu$}\hfill
 \mkern-6mu{-}\cr}}}
\let\amp@rs@nd@\relax
\newdimen\ex@
\ex@.2326ex
\newdimen\bigaw@
\newdimen\minaw@
\minaw@16.08739\ex@
\newdimen\minCDaw@
\minCDaw@2.5pc
\newif\ifCD@
\def\minCDarrowwidth#1{\minCDaw@#1}
\newenvironment{CD}{\@CD}{\@endCD}
\def\@CD{\def\A##1A##2A{\llap{$\vcenter{\hbox
 {$\scriptstyle##1$}}$}\Big\uparrow\rlap{$\vcenter{\hbox{%
$\scriptstyle##2$}}$}&&}%
\def\V##1V##2V{\llap{$\vcenter{\hbox
 {$\scriptstyle##1$}}$}\Big\downarrow\rlap{$\vcenter{\hbox{%
$\scriptstyle##2$}}$}&&}%
\def\={&\hskip.5em\mathrel
 {\vbox{\hrule width\minCDaw@\vskip3\ex@\hrule width
 \minCDaw@}}\hskip.5em&}%
\def\verteq{\Big\Vert&&}%
\def\noarr{&&}%
\def\vspace##1{\noalign{\vskip##1\relax}}\relax\let\amp@rs@nd@&\iffalse}\fi
 \CD@true\vcenter\bgroup\relax\let\\=\cr\iffalse}\fi\tabskip\z@skip\baselineskip20\ex@
 \lineskip3\ex@\lineskiplimit3\ex@\halign\bgroup
 &\hfill$\m@th##$\hfill\cr}
\def\@endCD{\cr\egroup\egroup}
\def\>#1>#2>{\amp@rs@nd@\setbox\z@\hbox{$\scriptstyle
 \;{#1}\;\;$}\setbox\@ne\hbox{$\scriptstyle\;{#2}\;\;$}\setbox\tw@
 \hbox{$#2$}\ifCD@
 \global\bigaw@\minCDaw@\else\global\bigaw@\minaw@\fi
 \ifdim\wd\z@>\bigaw@\global\bigaw@\wd\z@\fi
 \ifdim\wd\@ne>\bigaw@\global\bigaw@\wd\@ne\fi
 \ifCD@\hskip.5em\fi
 \ifdim\wd\tw@>\z@
 \mathrel{\mathop{\hbox to\bigaw@{\rightarrowfill}}\limits^{#1}_{#2}}\else
 \mathrel{\mathop{\hbox to\bigaw@{\rightarrowfill}}\limits^{#1}}\fi
 \ifCD@\hskip.5em\fi\amp@rs@nd@}
\def\<#1<#2<{\amp@rs@nd@\setbox\z@\hbox{$\scriptstyle
 \;\;{#1}\;$}\setbox\@ne\hbox{$\scriptstyle\;\;{#2}\;$}\setbox\tw@
 \hbox{$#2$}\ifCD@
 \global\bigaw@\minCDaw@\else\global\bigaw@\minaw@\fi
 \ifdim\wd\z@>\bigaw@\global\bigaw@\wd\z@\fi
 \ifdim\wd\@ne>\bigaw@\global\bigaw@\wd\@ne\fi
 \ifCD@\hskip.5em\fi
 \ifdim\wd\tw@>\z@
 \mathrel{\mathop{\hbox to\bigaw@{\leftarrowfill}}\limits^{#1}_{#2}}\else
 \mathrel{\mathop{\hbox to\bigaw@{\leftarrowfill}}\limits^{#1}}\fi
 \ifCD@\hskip.5em\fi\amp@rs@nd@}
\newenvironment{CDS}{\@CDS}{\@endCDS}
\def\@CDS{\def\A##1A##2A{\llap{$\vcenter{\hbox
 {$\scriptstyle##1$}}$}\Big\uparrow\rlap{$\vcenter{\hbox{%
$\scriptstyle##2$}}$}&}%
\def\V##1V##2V{\llap{$\vcenter{\hbox
 {$\scriptstyle##1$}}$}\Big\downarrow\rlap{$\vcenter{\hbox{%
$\scriptstyle##2$}}$}&}%
\def\={&\hskip.5em\mathrel
 {\vbox{\hrule width\minCDaw@\vskip3\ex@\hrule width
 \minCDaw@}}\hskip.5em&}
\def\verteq{\Big\Vert&}
\def\novarr{&}
\def\noharr{&&}
\def\SE##1E##2E{\slantedarrow(0,18)(4,-3){##1}{##2}&}
\def\SW##1W##2W{\slantedarrow(24,18)(-4,-3){##1}{##2}&}
\def\NE##1E##2E{\slantedarrow(0,0)(4,3){##1}{##2}&}
\def\NW##1W##2W{\slantedarrow(24,0)(-4,3){##1}{##2}&}
\def\slantedarrow(##1)(##2)##3##4{%
\thinlines\unitlength1pt\lower 6.5pt\hbox{\begin{picture}(24,18)%
\put(##1){\vector(##2){24}}%
\put(0,8){$\scriptstyle##3$}%
\put(20,8){$\scriptstyle##4$}%
\end{picture}}}
\def\vspace##1{\noalign{\vskip##1\relax}}\relax\let\amp@rs@nd@&\iffalse}\fi
 \CD@true\vcenter\bgroup\relax\let\\=\cr\iffalse}\fi\tabskip\z@skip\baselineskip20\ex@
 \lineskip3\ex@\lineskiplimit3\ex@\halign\bgroup
 &\hfill$\m@th##$\hfill\cr}
\def\@endCDS{\cr\egroup\egroup}
\newdimen\TriCDarrw@
\newif\ifTriV@
\newenvironment{TriCDV}{\@TriCDV}{\@endTriCD}
\newenvironment{TriCDA}{\@TriCDA}{\@endTriCD}
\def\@TriCDV{\TriV@true\def\TriCDpos@{6}\@TriCD}
\def\@TriCDA{\TriV@false\def\TriCDpos@{10}\@TriCD}
\def\@TriCD#1#2#3#4#5#6{%
\setbox0\hbox{$\ifTriV@#6\else#1\fi$}
\TriCDarrw@=\wd0 \advance\TriCDarrw@ 24pt
\advance\TriCDarrw@ -1em
\def\SE##1E##2E{\slantedarrow(0,18)(2,-3){##1}{##2}&}
\def\SW##1W##2W{\slantedarrow(12,18)(-2,-3){##1}{##2}&}
\def\NE##1E##2E{\slantedarrow(0,0)(2,3){##1}{##2}&}
\def\NW##1W##2W{\slantedarrow(12,0)(-2,3){##1}{##2}&}
\def\slantedarrow(##1)(##2)##3##4{\thinlines\unitlength1pt
\lower 6.5pt\hbox{\begin{picture}(12,18)%
\put(##1){\vector(##2){12}}%
\put(-4,\TriCDpos@){$\scriptstyle##3$}%
\put(12,\TriCDpos@){$\scriptstyle##4$}%
\end{picture}}}
\def\={\mathrel {\vbox{\hrule
   width\TriCDarrw@\vskip3\ex@\hrule width
   \TriCDarrw@}}}
\def\>##1>>{\setbox\z@\hbox{$\scriptstyle
 \;{##1}\;\;$}\global\bigaw@\TriCDarrw@
 \ifdim\wd\z@>\bigaw@\global\bigaw@\wd\z@\fi
 \hskip.5em
 \mathrel{\mathop{\hbox to \TriCDarrw@
{\rightarrowfill}}\limits^{##1}}
 \hskip.5em}
\def\<##1<<{\setbox\z@\hbox{$\scriptstyle
 \;{##1}\;\;$}\global\bigaw@\TriCDarrw@
 \ifdim\wd\z@>\bigaw@\global\bigaw@\wd\z@\fi
 \mathrel{\mathop{\hbox to\bigaw@{\leftarrowfill}}\limits^{##1}}
 }
 \CD@true\vcenter\bgroup\relax\let\\=\cr\iffalse}\fi
 \tabskip\z@skip\baselineskip20\ex@
 \lineskip3\ex@\lineskiplimit3\ex@
 \ifTriV@
 \halign\bgroup
 &\hfill$\m@th##$\hfill\cr
#1&\multispan3\hfill$#2$\hfill&#3\\
&#4&#5\\
&&#6\cr\egroup%
\else
 \halign\bgroup
 &\hfill$\m@th##$\hfill\cr
&&#1\\%
&#2&#3\\
#4&\multispan3\hfill$#5$\hfill&#6\cr\egroup
\fi}
\def\@endTriCD{\egroup}
\newcommand{\sA}{{\mathcal A}}
\newcommand{\sB}{{\mathcal B}}
\newcommand{\sC}{{\mathcal C}}
\newcommand{\sD}{{\mathcal D}}
\newcommand{\sE}{{\mathcal E}}
\newcommand{\sF}{{\mathcal F}}
\newcommand{\sG}{{\mathcal G}}
\newcommand{\sH}{{\mathcal H}}
\newcommand{\sI}{{\mathcal I}}
\newcommand{\sJ}{{\mathcal J}}
\newcommand{\sK}{{\mathcal K}}
\newcommand{\sL}{{\mathcal L}}
\newcommand{\sM}{{\mathcal M}}
\newcommand{\sN}{{\mathcal N}}
\newcommand{\sO}{{\mathcal O}}
\newcommand{\sP}{{\mathcal P}}
\newcommand{\sQ}{{\mathcal Q}}
\newcommand{\sR}{{\mathcal R}}
\newcommand{\sS}{{\mathcal S}}
\newcommand{\sT}{{\mathcal T}}
\newcommand{\sU}{{\mathcal U}}
\newcommand{\sV}{{\mathcal V}}
\newcommand{\sW}{{\mathcal W}}
\newcommand{\sX}{{\mathcal X}}
\newcommand{\sY}{{\mathcal Y}}
\newcommand{\sZ}{{\mathcal Z}}
\newcommand{\A}{{\mathbb A}}
\newcommand{\B}{{\mathbb B}}
\newcommand{\C}{{\mathbb C}}
\newcommand{\D}{{\mathbb D}}
\newcommand{\E}{{\mathbb E}}
\newcommand{\F}{{\mathbb F}}
\newcommand{\G}{{\mathbb G}}
\newcommand{\HH}{{\mathbb H}}
\newcommand{\I}{{\mathbb I}}
\newcommand{\J}{{\mathbb J}}
\renewcommand{\L}{{\mathbb L}}
\newcommand{\M}{{\mathbb M}}
\newcommand{\N}{{\mathbb N}}
\newcommand{\BP}{{\mathbb P}}
\newcommand{\Q}{{\mathbb Q}}
\newcommand{\R}{{\mathbb R}}
\newcommand{\BS}{{\mathbb S}}
\newcommand{\T}{{\mathbb T}}
\newcommand{\U}{{\mathbb U}}
\newcommand{\V}{{\mathbb V}}
\newcommand{\W}{{\mathbb W}}
\newcommand{\X}{{\mathbb X}}
\newcommand{\Y}{{\mathbb Y}}
\newcommand{\Z}{{\mathbb Z}}
\newcommand{\id}{{\rm id}}
\newcommand{\rk}{{\rm rk}}
\newcommand{\END}{{\mathbb E}{\rm nd}}
\newcommand{\End}{{\rm End}}
\newcommand{\Hg}{{\rm Hg}}
\newcommand{\tr}{{\rm tr}}
\newcommand{\Sl}{{\rm Sl}}
\newcommand{\Gl}{{\rm Gl}}
\newcommand{\Sp}{{\rm Sp}}
\newcommand{\MT}{{\rm MT}}
\newcommand{\Cor}{{\rm Cor}}
\newcommand{\Hom}{{\sH}{\rm om}}
\newcommand{\Mon}{{\rm Mon}}
\newcommand{\s}{{\rm sl}}
\newcommand{\ch}{{\rm c}}
\newcommand{\gr}{{\mathfrak g \mathfrak r}}
\newcommand{\doubledot}{{\mbox{\Tiny\textbullet,\textbullet}}}
\begin{abstract}
Let $Y$ be a non-singular projective manifold with an ample canonical sheaf, and let $\V$ be a 
$\Q$-variation of Hodge structures of weight one on $Y$ with Higgs bundle $E^{1,0} \oplus E^{0,1}$, 
coming from a family of Abelian varieties. If $Y$ is a curve the Arakelov inequality says that 
the slopes satisfy $\mu(E^{1,0})-\mu(E^{0,1})\leq \mu(\Omega^1_Y)$.
 
We prove a similar inequality in the higher dimensional case. If the latter is an equality, 
and if the discriminant of $E^{1,0}$ or the one of $E^{0,1}$ is zero, one hopes that $Y$ is a Shimura
variety, and $\V$ a uniformizing variation of Hodge structures. 
This is verified, in case the universal covering of $Y$ does not contain factors of rank $>1$.
Part of the results extend to variations of Hodge structures over quasi-projective manifolds $U$.
\end{abstract}
\maketitle
\tableofcontents
\renewcommand{\thethm}{\arabic{thm}} \renewcommand{\theequation}{\arabic{equation}}     
Let $Y$ be a complex $n$-dimensional projective manifold, $S\subset Y$ a 
reduced normal crossing divisor,
$U=Y\setminus S$, and let $f:V \to U$ be a smooth family of $g$-dimensional Abelian varieties.
Assume that the local system $R^1f_*\C_V$ has uni-potent monodromy along the
components of $S$. Let $\V$ be a $C$-sub-variation of Hodge structures in
$R^1f_*\C_V$.

The Deligne extension of $\V\otimes \sO_U$ to $Y$
carries a Hodge filtration. Taking the graded sheaf one obtains
the (logarithmic) Higgs bundle
$$
(E, \theta)=(E^{1,0}\oplus E^{0,1}, \theta),
$$
where $\theta: E \to E\otimes \Omega_Y^1(\log S)$ is zero on $E^{0,1}$ and
factors through
$$
\theta:E^{1,0}\>>> E^{0,1}\otimes \Omega_Y^1(\log S)
$$
on $E^{1,0}$. Define for a torsion free coherent sheaf $\sF$ on $Y$
\begin{gather*}
\Upsilon(\sF)=\frac{\ch_1(\sF)}{\rk(\sF)} \in H^2(Y,\Q) \mbox{ \ \ and}\\
\Delta(\sF)= 2\cdot \rk(\sF) \cdot \ch_2(\sF) - (\rk(\sF)-1)\cdot \ch_1(\sF)^2 \in H^4(Y,\Q).
\end{gather*}
Over the $n$-dimensional variety $Y$ choose an invertible sheaf $\sN$, or more generally
an $\R$-divisor $\sN$, and define the slope and the discriminant of $\sF$ as
$$
\mu_\sN(\sF)=\Upsilon(\sF).\ch_1(\sN)^{n-1} \mbox{ \ \ and \ \ }
\delta_\sN(\sF)=\Delta(\sF).\ch_1(\sN)^{n-2},
$$
respectively. For the Higgs bundle $(E,\theta)$ of the variation of Hodge structures
$\V$ we define
\begin{gather*}
\mu_\sN(\V)=\mu_\sN(E^{\doubledot})=\mu_\sN(E^{1,0})-\mu_\sN(E^{0,1})\mbox{ \ \ and}\\
\notag
\delta_\sN(\V)=\delta_\sN(E^{\doubledot})={\rm Min}\{\delta_\sN(E^{1,0}), \ \delta_\sN(E^{0,1})\}. 
\end{gather*}
We will choose $\sN=\omega_Y(S)$ in the introduction and we will write $\mu$ and $\delta$ instead of 
$\mu_{\omega_Y(S)}$ and $\delta_{\omega_Y(S)}$.

For $Y$ a curve and $\V=R^1f_*\C_V$ the Arakelov inequality due to Faltings \cite{Fa} says that
$$
\deg(E^{1,0})-\deg(E^{0,1})=2\cdot \deg(E^{1,0}) \leq \rk(E^{1,0}) \cdot \deg(\Omega^1_Y(\log S)),
$$
or equivalently that
\begin{equation}\label{eqin.1}
\mu(R^1f_*\C_V) \leq \mu(\Omega^1_Y(\log S)).
\end{equation}
If (\ref{eqin.1}) is an equality the Higgs field $\theta$ is an isomorphism. As shown in
\cite{VZ1}, this forces $U$ to be a Shimura curve and $R^1f_*\Z_V$ to be a uniformizing
variation of Hodge structures.
An intermediate result says, that $R^1f_*\C_V=\L\otimes \U$, for a unitary bundle $\U$
and for a uniformizing variation of Hodge structures $\L$ of weight one and rank two, i.e.
for some $\L$ with Higgs bundle 
$(\sL \oplus \sL^{-1}, \tau:\sL \to \sL^{-1}\otimes \omega_Y(S))$,
where $\sL^2\cong \omega_Y(S)$. In particular $E^{1,0}=\sL\otimes_\C\U$
and $E^{0,1}=\sL^{-1}\otimes_\C\U$ are both poly-stable.

We want to obtain inequalities similar to (\ref{eqin.1}), hoping that equality forces $Y\setminus S$ to be a locally Hermitian symmetric domain, and the variation of Hodge structures to be again the standard one, up to the tensor product with a unitary local systems. 

Before stating the results, let us consider the example of a two dimensional compact ball
quotient $Y$. Replacing $Y$ by an \'etale cover, the uniformizing $\R$-variation of
Hodge structures splits over $\C$ as a direct sum $\V\oplus \bar\V$, and interchanging
$\V$ and $\bar\V$, if necessary, the Higgs bundles of $\V$ and $\bar\V$ are given by
$(E^{1,0}\oplus E^{0,1},\theta)$ and $({E'}^{1,0}\oplus {E'}^{0,1}\theta')$ for
$$
E^{1,0}=\omega_Y ^{\frac{1}{3}} , \ \ E^{0,1}=T^1_Y\otimes \omega_Y^{\frac{1}{3}}
\mbox{ \ and \ } {E'}^{1,0}=\Omega_Y^1\otimes \omega_Y^{-\frac{1}{3}} , \ \ {E'}^{0,1}=
\omega_Y^{-\frac{1}{3}},
$$
respectively (see \cite[4.1]{Lo}, \cite[9.1]{Sim} and Section \ref{co}).
One finds
$$
\mu(\V)= \mu(\bar{\V}) = \mu(\Omega_Y^1),
$$
whereas for the whole variation of Hodge structures one has a strict inequality
$\mu(\V\oplus \bar\V) < \mu_{\omega_Y}(\Omega_Y^1)$. So one should expect optimal Arakelov type inequalities only for the Hodge bundles of irreducible local sub-systems.

\begin{theorem}\label{in.1} Assume that $\omega_Y(S)$ is nef and
ample with respect to $U$, and let $f:V\to U$ be a smooth family of polarized $g$-dimensional Abelian varieties
such that the local monodromy of $R^1f_*\C_V$ around the components of $S$ is uni-potent.
Let $\V$ be a sub-variation of Hodge structures in $R^1f_*\C_V$ without
a unitary direct factor and with Higgs bundle 
$$(E=E^{1,0}\oplus E^{0,1},\theta).$$ 
Then 
\begin{equation}\label{eqin.2}
\mu(\V) \leq \mu(\Omega^1_Y(\log S)) 
\end{equation}
The equality
\begin{equation}\label{eqin.3}
\mu(\V) = \mu(\Omega^1_Y(\log S))
\end{equation}
implies that $E^{1,0}$ and $E^{0,1}$ are both semi-stable, and that
\begin{equation}\label{eqin.4}
\delta(\V)\geq 0.
\end{equation}
\end{theorem} 
Here ``semi-stable'' refers to semi-stability for the slope $\mu$. The definition, as well as
the condition ``ample with respect to $U$'' will be recalled in \ref{ya.1}. 

The inequality (\ref{eqin.4}) follows immediately from the semi-stability
of $E^{1,0}$ and $E^{0,1}$ and the Bogomolov inequality, saying that the discriminant of a semi-stable locally free sheaves is non negative (see for example \cite[7.3.1]{HL}).

\begin{remark}\label{in.2} As we will see in the proof of Theorem \ref{in.1} at the end of Section \ref{sl}
one can allow $\V$ to have unitary direct factors $\U$ which are invariant under complex conjugation. In particular the inequality (\ref{eqin.2}) holds for all $\R$-sub-variations $\V$ of Hodge structures, and for $\V=R^1f_*\C_V$ itself. 

If the equality (\ref{eqin.3}) holds for some $\V$, it holds for all irreducible $\C$-sub-variations of
Hodge structures in $\V$. So it will never hold if $\V$ contains non-trivial unitary direct factors. 
\end{remark}
The assumption ``$\omega_Y(S)$ nef and ample with respect to $U$'' allows to apply
Yau's uniformization theorem \cite{Y}, recalled in \ref{ya.4}. It implies that the sheaf
$\Omega_Y^1(\log S)$ is $\mu$-poly-stable. Hence one has a direct sum decomposition
$$
\Omega_Y^1(\log S)=\Omega_1\oplus \cdots \oplus \Omega_{s''}\oplus \cdots \oplus \Omega_{s'}
\oplus \cdots \oplus \Omega_s
$$
in stable sheaves. We choose the indices such that for $i=1, \ldots , s''$
the sheaf $\Omega_i$ is invertible. For $i=s''+1, \ldots , s'$, we assume that the sheaves 
$S^m(\Omega_i)$ are stable for all $m>1$, and not invertible. Finally
for $i=s'+1, \ldots , s$ we have the remaining stable direct factors,
i.e. those with $S^{m_i}(\Omega_i)$ non-stable, for some $m_i>1$.

We will need stronger positivity conditions. First of all we have to require
the sheaves $\Omega_i$ to be nef, or equivalently $\Omega_Y^1(\log S)$ to be nef.
If $S=\emptyset$ this conditions holds true for projective submanifolds $Y$ of the
moduli stack $\sA_g$ of polarized Abelian varieties (see Lemma \ref{av.1}).
In general, the assumption ``$\Omega^1_Y(S)$ nef'' depends on the choice of a ``good'' compactification $Y$ of $U$. For the moduli space itself, such a compactification is described in \cite{Fa}.

Secondly we will frequently need the following properties of the slope and the discriminant.
\begin{enumerate}
\item[($*$)] The the tensor product of $\sF\otimes\sG$ of $\mu$-poly-stable sheaves
$\sF$ and $\sG$ is is again $\mu$-poly-stable.
\item [($**$)] A locally free $\mu$-poly-stable sheaf $\sF$ is unitary if 
$\mu(\sF)=\delta(\sF)=0$.
\end{enumerate}
Note that both, ($*$) and ($**$), hold true, if $\omega_Y(S)$ is ample (see Lemma \ref{ya.3}).
Again, as we will see in  Lemma \ref{av.1}, this condition holds for projective submanifolds $Y$ of the moduli stack $\sA_g$. However for $S\neq \emptyset$ ampleness is too much to expect. S.T. Yau informed us, that ($*$) and ($**$) hold true under the assumption that $\omega_Y(S)$ is nef and ample with respect to some open dense subscheme, and that the proof will be given in a forthcoming article by Sun and Yau. 

So we will work with the following Set-up, noting that at present, i.e. without using
the unpublished result of Sun and Yau, the conditions are only reasonable for $S=\emptyset$.
\begin{setup}\label{in.3} 
The sheaf $\Omega^1_Y(\log S)$ is nef, $\omega_Y(S)$ is ample with respect to
$U$ and the conditions ($*$) and ($**$) hold true for $\mu=\mu_{\omega_Y(S)}$ and
$\delta=\delta_{\omega_Y(S)}$. 

Let $f:V\to U$ be a smooth family of polarized $g$-dimensional Abelian varieties
such that the local monodromy of $R^1f_*\C_V$ around the components of $S$
is uni-potent. 
\end{setup} 
For surfaces $Y$ with $U\varsubsetneq Y$ we will consider in \ref{su.2} a slightly different Set-up. 

\begin{proposition}\label{in.4}
In the Set-up \ref{in.3} assume that $s=s'$, that for all irreducible $\C$-sub-variations $\V$ of 
Hodge structures in $R^1f_*\C_V$ with logarithmic Higgs bundle $(E^{1,0} \oplus E^{0,1},\theta)$
one has $\mu(\V) = \mu(\Omega^1_Y(\log S))$, and that one of the following conditions holds true
\begin{enumerate}
\item[i.] $\delta(\V)=0,$
\item[ii.] $s''=s$.
\end{enumerate}
Then there exists some $i\in \{1,\ldots, s'\}$ such that: 
\begin{enumerate}
\item[a.] The Higgs field $\theta$ factors like
$$
E^{1,0} \> \theta_i  >> E^{0,1}\otimes \Omega_i \> \subset >> E^{0,1}\otimes \Omega_Y^1(\log S).
$$
\item[b.] $E^{1,0}$ and $E^{0,1}$ are stable.
\item[c.] Either $\rk(E^{1,0})=\rk(E^{0,1})+\rk(\Omega_i)$ or
$\rk(E^{0,1})=\rk(E^{1,0})+\rk(\Omega_i)$.
\end{enumerate}
\end{proposition}
If $\omega_Y(S)$ is ample, one can replace the condition c) in \ref{in.4} by:
\begin{enumerate}
\item[c.] Either $\theta_i$ is an isomorphism or $E^{0,1}\cong E^{1,0}\otimes \Omega_i^\vee$
and $\theta_i$ is the natural map
$$
E^{1,0} \>>> E^{1,0}\otimes \Omega_i^\vee \otimes \Omega_i\cong E^{0,1}\otimes \Omega_i. 
$$
\end{enumerate}
We hope that the conclusion a) of Proposition \ref{in.4} remains true if one allows
stable factors $\Omega_i$ of the third type, i.e. those with $S^m(\Omega_i)$
non-stable for some $m>0$. Moreover the conditions i) and ii) should not be necessary at this place.
So we will not use the condition $s=s'$ in the next Theorem, and refer instead to the conclusion of Theorem \ref{in.4}.

\begin{theorem}\label{in.5} In the Set-up \ref{in.3} let $\V$ be an irreducible $\C$-sub-variations of Hodge structures in $R^1f_*\C_V$ with Higgs bundle $(E^{1,0} \oplus E^{0,1},\theta)$.
Assume that $\mu(\V)=\mu(\Omega^1_Y(\log S))$ and that for some $i\in \{1,\ldots, s'\}$ the conditions a)--c) in Proposition \ref{in.4}
hold true. If either $i\leq s''$, or if $s''<i\leq s'$ and $\delta(\V)=0$, there exists an \'etale covering $\phi:Y'\to Y$, and an invertible sheaf $\sL_i$ with:
\begin{enumerate}
\item[a.] $\displaystyle \sL_i^{n_i+1}=\phi^*\det(\Omega_i),$ for $n_i=\rk(\Omega_i)$.
\item[b.] $\V'=\phi^*(\V)$ or its dual is the tensor product of a unitary
local system $\U_i$, regarded as a variation of Hodge structures of bidegree $(0,0)$, with a
variation of Hodge structures $\L_i$ with Higgs bundle
$$
\big( \sL_i \ \oplus \ \sL_i \otimes \phi^*\Omega^\vee_i, \ \tau\big)
$$
where $\tau$ is given by the morphism
\begin{gather*}
\sL_i \hookrightarrow \sL_i
\otimes \sE nd(\phi^*\Omega_i)= \sL_i\otimes
\phi^*(\Omega^\vee_i\otimes \Omega_i)\hspace{2cm} \\ \hspace{6cm}
\>\subset >>  \sL_i\otimes \phi^*(\Omega^\vee_i) \otimes
\Omega^1_{Y'}(\log S'),
\end{gather*}
induced by the homotheties $\sO_Y \hookrightarrow \sE nd(\Omega_i)$.
\end{enumerate}
\end{theorem}
The explicite form of the variation of Hodge structures given in Theorem \ref{in.5} will allow
in Section \ref{sh} to calculate the derived Mumford-Tate group of $\W=R^1f_*\Q_V$. To this aim, we have to study
in Section \ref{rg} the decomposition of certain wedge products of $\W$, and to determine the possible
Hodge cycles. This will finally allow to prove:
\begin{theorem}\label{in.6}
Under the assumptions made in Proposition \ref{in.4}, assume that the morphism 
$\varphi:U\to \sA_g$ to the moduli stack of polarized $g$-dimensional Abelian varieties, induced by $f:V\to U$, is generically finite. Then $U$ is a rigid Shimura subvariety of $\sA_g$. The universal covering $\tilde{U}$ of $U$ decomposes as the product of $s=s'$ complex balls of dimensions $n_i=\rk(\Omega_i)$.
\end{theorem}
The assumption that $s=s'$, i.e. that there are no direct factors $\Omega_i$ of the third type,
automatically holds if $Y$ is a surface. We will discuss this case in the second half of Section \ref{su}, and we will 
formulate and prove variants of Theorems \ref{in.5} and \ref{in.6} enforcing the conditions
($*$) and ($**$) in the Set-up \ref{in.3} by considering certain small twists of the slope $\mu$.

In Section \ref{ya} we will recall Yau's Uniformization Theorem, and some of its consequences.
The proof of Theorem \ref{in.1} will be given in Section \ref{sl}. 
In Section \ref{st} we prove Proposition \ref{in.4}, as well as the stability of the
Hodge bundles $E^{1,0}$ and $E^{0,1}$.

The next Section \ref{av} recalls some well known properties of the moduli stack of Abelian varieties,
and a first application of Proposition \ref{in.4}. At this stage we will also discuss the relation between our approach and the one of Moonen in \cite{Mo}. Moreover we will outline a possible approach towards a generalization of Theorem \ref{in.6} allowing factors $\Omega_i$ of the third type, i.e. with $S^m(\Omega_i)$ non-stable for some $i$.
 
As indicated in Proposition \ref{in.4} things are nicer if one assumes that all the direct factors
$\Omega_i$ of $\Omega_Y^1(\log S)$ are invertible. In this case, one obtains a numerical
characterization of generalized Hilbert modular varieties, stated and discussed in Section \ref{su}.
We will show in this section as well that for $i\leq s''$ the condition b) in Theorem \ref{in.5} 
is a consequence of a). 

As a first step towards Theorem \ref{in.5} for $s'' < i \leq s'$ we will show in
Section \ref{pr} that the condition $\delta(\V)=0$ implies that the factor of the universal
covering $\tilde{U}$ of $U$, corresponding to $\Omega_i$ is a complex ball.

At this stage we do not know the existence of an invertible sheaf of the form
$\sL_i=\det(\Omega_i)^{\frac{1}{n_i+1}}$, asked for in \ref{in.5}, a).
If $S=\emptyset$, knowing that $Y$ is a quotient of products of balls, will allow in Lemma \ref{co.2}
to apply the Simpson correspondence, and to construct the sheaf $\sL_i$.
At the end of Section \ref{co} we finish the proof of \ref{in.5}.

Beforehand, in Section \ref{nc} we consider the case $S\neq \emptyset$, or more precisely
the one where $S$ meets the leaves of the foliation defined by the direct factor $\Omega_i$
of $\Omega^1_Y(\log S)$. As it will turn out, in this case the unitary bundle in \ref{in.5}, b),
is trivial and one obtains the existence of $\sL_i$ ``for free''.

The last step is to show that the quotient of products of complex balls in Theorem
\ref{in.5} is a rigid Shimura variety. The rigidity is shown in Section \ref{rg}. There
we study the decomposition of $R^1f_*\R_V$ or $R^1f_*\C_V$ in $\R$ and $\C$ irreducible direct factors
in more detail. In particular the first one can be realized over a totally real number field.
We calculate the possible bidegrees of global sections of the wedge products of
$R^1f_*\C_V$. This will imply in Section \ref{sh} finally that $U$ is a Shimura variety, and allow to end
the proof of Theorem \ref{in.6}.\vspace{.2cm}

\noindent
This article relies on C. Simpson's correspondence between Higgs bundles and local systems
(\cite{Sim}, \cite{Sim1}, \cite{Sim2} and \cite{Sim3}). The second main ingredient is S.T. Yau's
uniformization theorem, recalled in Section \ref{ya}. It is based on the existence of K\"ahler-Einstein
metrics, due to Yau in the projective case, and extended to the quasi-projective case jointly with
G. Tian. We thank both of them for explaining how to use their results to study the uniformization of manifolds.

Thanks also to Y-H. Yang, who clarified and verified several arguments from differential geometry,
and to F. Bogomolov who explained the first named author his view of ball
quotients, and who told us Lemma \ref{co.1}, as well as its proof.

We are grateful to the referee for pointing out several mistakes and ambiguities in the first version
of this article, and for suggestions how to improve the presentation of our results.

The first steps towards Arakelov inequalities over a higher dimensional bases were done when the
first named author visited the Courant Institute, New York, the final version of this article was
written during his stay at the IAS, Princeton. He would like to thank the members of both
Institutes for their hospitality.\vspace{.2cm}

{\bf Notations.} As in the Introduction we will consider up to Section \ref{sh}
$\C$-sub-variations of Hodge structures $\V$ in $R^1f_*\C_V$. We will say that $\V$
is defined over some subfield $K\subset \C$, if there exists a $K$-sub-variation of Hodge structures
$\V_K\subset R^1f_*K_V$ with $\V_K\otimes \C=\V$.

If $(E,\theta)$ is a Higgs bundle, and $\sG \subset E$ a subsheaf, we call $\sG$ a Higgs subsheaf,
if $\theta(\sG)\subset \sG\otimes \Omega^1_Y(\log S)$. We will call $\sG$ a saturated
Higgs subsheaf if in addition $E/\sG$ is torsion free. If $\sG$ is a saturated Higgs subsheaf
$E/\sG$ will be called a quotient Higgs sheaf.

The dual of $\Omega^1_Y(\log S)$ will be denoted by $T^1_Y(-\log S)$.

As in the introduction, $\mu$ and $\delta$ usually stand for the slope and discriminant with respect to the invertible sheaf $\omega_Y(S)$. However in parts of the article we will allow 
$\mu=\mu_\sN$ and $\delta=\delta_\sN$, for ample invertible sheaves (or $\R$-divisors) $\sN$, provided 
the assumptions madein \ref{sl.5} hold true.

Stability, semi-stability and poly-stability will always be for the slope $\mu$. Just in case we want to underline that we allow different polarization, we write $\mu_\sN$, $\delta_\sN$ and we will talk about $\mu_\sN$-stability.

\section{Stability and locally Hermitian symmetric spaces}\label{ya}
\renewcommand{\thethm}{\arabic{section}.\arabic{thm}} \renewcommand{\theequation}{\arabic{section}.\arabic{equation}}
Let us recall some properties of locally free sheaves $\sF$ on a manifold $Y$ of dimension $n$.
\begin{definition}\label{ya.1} Let $\sH$ be an ample invertible sheaf on $Y$.
\begin{enumerate}
\item[i.] $\sF$ is numerically effective (nef), if for all curves $\tau:C\to Y$ and for all
invertible quotient sheaves $\sL$ of $\tau^*\sF$ one has $\deg(\sL)\geq 0$.
\item[ii.] $\sF$ is ample with respect to an open subscheme $U'$ of $Y$, if for some $\nu \gg 0$
there exists a morphism $\oplus \sH \to S^\nu(\sF)$, which is surjective on $U'$.
If $\sF$ is invertible, this is equivalent to:\\
For some $\eta > 0$ the sheaf $\sF^\eta$ is generated by $H^0(Y,\sF^\eta)$ in all points $u\in U'$,
and the induced morphism $U'\to \BP(H^0(Y,\sF^\eta)$ is an embedding.
\item[iii.] $\sF$ is big, if it is ample with respect to some open dense subscheme.
\end{enumerate}
\end{definition}
Let us also recall the notion of stability for a torsion free coherent sheaf $\sF$ with respect to an
invertible sheaf or an $\R$-divisor $\sN$. 
\begin{enumerate}
\item[iv.] $\sF$ is $\mu_\sN$-stable, if $\mu_\sN(\sG) < \mu_\sN(\sF)$
for all subsheaves $\sG$ of $\sF$ with $\rk(\sG)<\rk(\sF)$.
\item[v.] $\sF$ is $\mu_\sN$-semi-stable, if $\mu_\sN(\sG) \leq \mu_\sN(\sF)$ for all subsheaves $\sG$ of $\sF$.
\item[vi.] A $\mu_\sN$-semi-stable sheaf $\sF$ is $\mu_\sN$-poly-stable, if it is the direct sum of
$\mu_\sN$-stable sheaves.
\end{enumerate}
In iv), v) and vi) we will assume that $\sN$ is nef and big. If $\sN$ is not ample, and if
$\sF$ is a locally free $\mu_\sN$-stable sheaf it might happen that $\sF$ contains a stable
locally free subsheaf $\sG \varsubsetneq \sF$ with $\mu_\sN(\sG)=\mu_\sN(\sF)$. By definition this is only possible if $\rk(\sG)=\rk(\sF)$ and if the cokernel of $\sG \hookrightarrow \sF$ is supported on divisors $D$ with $D.\ch_1(\sN)^{n-1}=0$.

Recall that the slope is $\mu_\sN(\sF)=\Upsilon(\sF).\ch_1(\sN)^{n-1}$.
So the standard properties of the slope follow from:
\begin{lemma}\label{ya.2} Let $\sF$ and $\sG$ be torsion free coherent sheaves on $Y$.
\begin{gather*}
\Upsilon(\sF\otimes \sG)=\frac{\rk(\sG)\cdot\ch_1(\sF)+
\rk(\sF)\cdot\ch_1(\sG)}{\rk(\sF)\cdot\rk(\sG)}=\Upsilon(\sF)+\Upsilon( \sG),\\
\Upsilon(S^m(\sF))=m\cdot \Upsilon(\sF), \hspace{1cm}
\Upsilon(\bigwedge^m(\sF))=m\cdot \Upsilon(\sF)
\end{gather*}
where for the third equality one assumes $m\leq \rk(\sF)$. If
$$
0\>>> \sF \>>> \sK \>>> \sG \>>> 0
$$
is an exact sequence, $\Upsilon(\sK)$ is equal to
$$
\frac{\ch_1(\sF)+\ch_1(\sG)}{\rk(\sF) + \rk(\sG)}=
\frac{\rk(\sF)}{\rk(\sF) + \rk(\sG)}\Upsilon(\sF)+\frac{\rk(\sG)}{\rk(\sF) + \rk(\sG)}\Upsilon( \sG).
$$
\end{lemma}
\begin{lemma}\label{ya.3} \
\begin{enumerate}
\item[a.] There exists the Harder-Narasimhan filtration for all torsion free
coherent sheaves $\sF$, i.e. a filtration
$$
0=\sF_0 \subset \sF_1 \subset \cdots \subset \sF_m
$$
with $\sF_\ell/\sF_{\ell-1}$ torsion free, $\mu_\sN$-semi-stable and with
\begin{gather*}
\mu_{{\rm max},\sN}(\sF)=\mu_\sN(\sF_1) >
\mu_\sN(\sF_2/\sF_{1}) >\hspace*{2cm} \\ \hspace*{2cm}
\cdots > \mu_\sN(\sF_m/\sF_{m-1})=\mu_{{\rm min},\sN}(\sF).
\end{gather*}
\item[b.] The tensor product of $\mu_\sN$-semi-stable sheaves is again $\mu_\sN$-semi-stable.
\item[c.] If $\sN$ is ample, the tensor product of $\mu_\sN$-poly-stable sheaves is again
$\mu_\sN$-poly-stable, as well as the pullback under finite morphisms.
\item[d.] If $\sN$ is ample and if $\sF$ is a locally free sheaf, $\mu_\sN$-poly-stable of slope
$\mu_\sN(\sF)=)$ and with $\delta_\sN(\sF)=0$, then $\sF$ is unitary, i.e. there exists a unitary 
local system $\U$ on $Y$ with $\sF=\U\otimes \sO_Y$.
\end{enumerate}
\end{lemma}
\begin{proof}
If $\sN$ is ample, the existence of the Harder-Narasimhan filtration is well known
(see \cite[1.6.7]{HL}), and the proof carries over to the case $\sN$ nef and big.
c) is shown in \cite[3.2.3 and 3.2.11]{HL}, for example, and d) in \cite{UY} and \cite{D}.

For $\sN$ ample, b) is well known (see \cite[3.1.4]{HL}).
In general one can use the following argument. Given an ample invertible sheaf $\sH$, consider
the Harder-Narasimhan filtration
$$
0=\sF_0 \subset \sF_1 \subset \cdots \subset \sF_m
$$
with respect to $\ch_1(\sN\otimes \sH^\epsilon)$.  By \cite[p. 263]{La} for $\epsilon>0$ and sufficiently
small, this filtration is independent of $\epsilon$. One has
$$
\mu_{\sN}(\sF_1) = \lim_{\epsilon \to 0} \mu_{\sN\otimes \sH^\epsilon}(\sF_1)
\geq \lim_{\epsilon \to 0} \mu_{\sN\otimes \sH^\epsilon}(\sF)=
\mu_{\sN}(\sF).
$$
If $\sF$ is $\mu_\sN$-stable, one knows that $\mu_{\sN}(\sF_1) \leq \mu_{\sN}(\sF)$
with equality if and only if $\sF_1=\sF$. So $\sF$ is $\mu_{\sN\otimes \sH^\epsilon}$-semi-stable.
Applying this to the graded sheaf with respect to a Jordan-H\"older filtration, one
finds that a $\mu_\sN$ semistable sheaf $\sF$ is also $\mu_{\sN\otimes \sH^\epsilon}$-semi-stable. 

So if $\sG$ is a second $\mu_\sN$-semi-stable sheaf, for $\epsilon$ small enough, $\sF$ and $\sG$ are both
$\mu_{\sN\otimes \sH^\epsilon}$-semi-stable, hence $\sF\otimes \sG$ as well. Taking the limit
$\epsilon \to 0$ one finds $\sF\otimes \sG$ to be $\mu_\sN$-semi-stable.
\end{proof}

Besides of Simpson's correspondence the main technical tool used in this article is
Yau's uniformization theorem, explained in \cite{Y} and \cite{Y2}. So for the rest of this section,
if not stated otherwise, we will choose $\sN=\omega_Y(S)$ and $\mu=\mu_\sN$.
\begin{theorem}\label{ya.4}
Let $Y$ be a complex projective manifold and $S\subset Y$ a reduced normal
crossing divisor. Assume that that $\omega_Y(S)$ is nef and ample with respect to
$Y\setminus S$. Then:
\begin{enumerate}
\item[a.] For all $m \geq 0$ the sheaves $S^m(\Omega^1_Y(\log S))$ are poly-stable.
\item[b.] Let
$$
\Omega_Y^1(\log S)=\Omega_1\oplus \cdots \oplus \Omega_s
$$
be the decomposition of $\Omega_Y^1(\log S)$ in stable direct factors of the same slope
and $n_i=\rk(\Omega_i)$. Then for $i=1,\ldots,s$ the $(n_i,n_i)$ current
$$
2(n_i+1)\cdot \ch_2(\Omega_i).\ch_1(\Omega_i)^{n_i-2} -
n_i\cdot \ch_1(\Omega_i)^{n_i}.
$$
is semi-positive.
\item[c.] Choose in b) $0 \leq s'' \leq s' \leq s$ with:
\begin{enumerate}
\item[i.] For $1 \leq i \leq s''$  the sheaf $\Omega_i$ is of rank one.
\item[ii.] For $s'' < i \leq s'$ the rank $n_i$ of $\Omega_i$ is larger than one
and for all $m>0$ the sheaf $S^m(\Omega_i)$ is stable.
\item[iii.] For $s' < i \leq s$ and for some $m_i>0$ the sheaf $S^{m_i}(\Omega_i)$
is not stable.
\end{enumerate}
Then $M_i$ is a one dimensional ball for $i=1,\ldots,s''$ and a 
a bounded symmetric domain of rank $>1$ for $i=s'+1,\ldots,s$.

For $i=s''+1,\ldots,s'$ one has the equality
\begin{gather}\label{eqya.1}
2(n_i+1)\cdot \ch_2(\Omega_i).\ch_1(\Omega_i)^{n_i-2}.\ch_1(\omega_Y(S))^{n-n_i}
\hspace*{2cm}  \\ \notag \hspace*{4cm} =
n_i\cdot \ch_1(\Omega_i)^{n_i}.\ch_1(\omega_Y(S))^{n-n_i}
\end{gather}
if and only if $M_i$ is an $n_i$-dimensional complex ball.
\item[d.] In particular, if in c) the equation (\ref{eqya.1}) holds for 
$i=s''+1,\ldots,s'$, then $U=Y\setminus S$ is a quotient of a bounded symmetric domain by a discrete group with finite volume. 
\end{enumerate}
\end{theorem}
In \cite{Y} the assumptions i) and iii) in c) are expressed differently. There
it is required, that $S^{m'_i}(\Omega_i)$ contains a direct factor of rank one
of the same slope as $S^{m'_i}(\Omega_i)$. Obviously this condition holds true if $n_i=1$.
For $n_i>1$ it is equivalent to iii).

The proof of Theorem \ref{ya.4} gives in fact additional informations.
\begin{addendum}\label{ya.5}
Let $\Omega_i$ be one of the stable factors in \ref{ya.4}, b), and let $\U$ be an
irreducible unitary local system on $Y$. Assume that $\U\otimes_\C S^m(\Omega_i)$
is not $\mu$-stable for some $m>0$. Then for some $m_i>0$ the sheaf
$S^{m_i}(\Omega_i)$ is not $\mu$-stable, hence $s' < i \leq s$.
\end{addendum}

\begin{proof}[Some notations from the proof of \ref{ya.4}]
Let $\tilde U = M_1\times \cdots \times M_s$ be the decomposition of the universal covering of
$U$ which corresponds to the decomposition in \ref{ya.4}, b). In particular, $\dim M_i=n_i$
for all $i$ (see \cite[p. 272]{Y2}). The holonomy group $H_i$ of $T_i=\Omega_i^\vee$
with respect to the K\"ahler-Einstein metric (more precisely, the projection of the
K\"ahler-Einstein metric on $T^1_Y(-\log S)$ to $T_i$) is contained in ${\rm U}(n_i)$.

As explained in \cite[p. 479]{Y}, the condition in \ref{ya.4}, c), iii) on the non-stability
of $S^m(\Omega_i)$ is equivalent to the condition that $H_i\neq {\rm U}(n_i)$.
The latter holds if and only if $M_i$ is a Hermitian symmetric space of rank $>1$.

If $H_i={\rm U}(n_i)$, by \cite[p. 272]{Y2} (with some misprint in the sign),
$$
2\cdot (n_i +1) \cdot \ch_2(T_i)- n_i\cdot \ch_1(T_i)^2
$$
is a semi-positive $(2,2)$ form, hence 
\begin{equation}\label{eqya.2}
2\cdot (n_i +1) \cdot \ch_2(\Omega_i)- n_i\cdot \ch_1(\Omega_i)^2,
\end{equation}
as well. Then, for $n_i\geq 2$,
\begin{equation}\label{eqya.3}
2\cdot (n_i +1) \cdot \ch_2(\Omega_i).\ch_1(\Omega_i)^{n_i-2}
- n_i\cdot \ch_1(\Omega_i)^{n_i}
\end{equation}
is a semi-positive $(n_i,n_i)$ current. It is zero, if and only if $M_i$ has constant negative 
holomorphic sectional curvature, hence if it is isomorphic to the complex ball.
\end{proof}

\begin{proof}[Proof of the Addendum \ref{ya.5}]
Assume $S^{m'}(\Omega_i)$ is stable, for all $m'>0$. If $n_i = 1$, there is
nothing to show. For $n_i >1$ we have just seen, that the holonomy group of $T_i$ is ${\rm U}(n_i)$,
hence the one of $S^m(T_i)$ is $S^m({\rm U}(n_i))$.

The holonomy group of $\U$ with respect to a locally constant metric is trivial, hence
$\{{\rm Id}_\ell \}$, the identity of $\Gl_\ell$. By functoriality the holonomy group
of $\U\otimes S^m(T_i)$  with respect to the product of the locally
constant metric and the K\"ahler Einstein metric is
$$
\{{\rm Id}_\ell \} \otimes S^m({\rm U}(n_i))=S^m({\rm U}(n_i))\times \cdots \times S^m({\rm U}(n_i)).
$$
Suppose $\U\otimes S^m(T_i)$ is not stable. Then a splitting gives rise to a splitting of the
product metric, hence to a splitting of the holonomy group $\{{\rm Id}_\ell\} \otimes S^m({\rm U}(n_i))$.

Since $\U$ is irreducible, such a splitting of $\U\otimes S^m(T_i)$ can not
arise from the natural splitting  $S^m({\rm U}(n_i))\times...\times S^m({\rm U}(n_i))$. Thus, it forces
one of the factors $S^m({\rm U}(n_i))$ to split. As in \cite[p. 272]{Y2} this contradicts
the irreducibility of $S^{m'}(T_i)$ for all $m'>0$.
\end{proof}
For the next Lemma we need more than the positivity of $\omega_Y(S)$.
\begin{lemma}\label{ya.6} Assume that $\omega_Y(S)$ is ample with respect to $U$, and that the sheaf
$\Omega^1_Y(\log S)$ is nef.
\begin{enumerate}
\item[i.] Then all the stable direct factors $\Omega_i$ in \ref{ya.4}, b), and their determinants
$\det(\Omega_i)$ are nef, and $\ch_1(\Omega_j)^{n_j+1}$ is numerically trivial.
\item[ii.] For $\nu_1, \ldots,\nu_s$ with $ \nu_1+\cdots+\nu_s=n$ the product 
$$
\ch_1(\Omega_1)^{\nu_1}. \cdots . \ch_1(\Omega_s)^{\nu_s}
$$ 
is a positive multiple of $\ch_1(\omega_Y(S))^{n}$,
if $ \nu_\iota=n_\iota$ for $\iota=1,\ldots,s$, and zero otherwise.
\item[iii.] $\ch_1(\Omega_1)^{n_1}. \cdots . \ch_1(\Omega_s)^{n_s}>0$.
\item[iv.] If $D$ is an effective $\Q$ divisor with $D.\ch_1(\omega_Y(S))^{n-1}=0$, then
$$D.\ch_1(\Omega_1)^{\nu_1}. \cdots . \ch_1(\Omega_s)^{\nu_s}=0$$ for all 
$\nu_1, \ldots,\nu_s$ with $ \nu_1+\cdots+\nu_s=n-1$. 
\item[v.] Let ${\rm NS}_0$ denote the subspace of the Neron-Severi group ${\rm NS}(Y)_\Q$ of $Y$ which is generated by all effective divisors $D$ satisfying the condition in iv). 
If for some $\alpha\in \Q$ one has $\ch_1(\Omega_i)-\alpha\cdot\ch_1(\Omega_j) \in {\rm NS}_0$ then $i = j$.
\item[vi.] If for some $m>0$ there is an injection $\rho: S^m(\Omega_i)\hookrightarrow S^m(\Omega_j)$, then $i = j$. 
\item[vii.] The equality (\ref{eqya.1}) in Theorem \ref{ya.4}, c), holds if and only if
$$
2(n_i+1)\cdot \ch_2(\Omega_i).\ch_1(\omega_Y(S))^{n-2}
= n_i\cdot \ch_1(\Omega_i)^{2}.\ch_1(\omega_Y(S))^{n-2}.
$$
\end{enumerate}
\end{lemma}
\begin{proof}
In Theorem \ref{ya.4} the sheaf $\Omega_Y^1(\log S)$ is poly-stable, and the assumption
obviously implies that all the $\Omega_i$ are nef, hence $\det(\Omega_j)$ as well. The Bogomolov-Sommese vanishing theorem (see \cite[6.9]{EV}) implies that $\kappa(\det(\Omega_j)) \leq n_j$,
hence $\ch_1(\Omega_j)^{n_j+1}$ is numerically trivial and i) holds true.

The cycle $\ch_1(\omega_Y(S))^\eta$ is a linear combination of expressions of the form
$$
\ch_1(\Omega_1)^{\nu_1}. \cdots . \ch_1(\Omega_s)^{\nu_s}
\mbox{ \ \ with \ \ } \nu_1+\cdots+\nu_s=\eta,
$$
with non-negative coefficients. For $\eta=n$ all those intersection cycles are zero, except the one with $\nu_j=n_j$, for all $j$. Since $\omega_Y(S)$ is big, one finds for some positive rational number $\alpha$
$$
\ch_1(\Omega_1)^{n_1}. \cdots . \ch_1(\Omega_s)^{n_s}
=\alpha \cdot \ch_1(\omega_Y(S))^{n} > 0,
$$
hence ii) and iii) hold true.

By i) for $\underline{\nu}=(\nu_1,\ldots,\nu_s)$ with $\nu_1+\cdots+\nu_s=\eta=n-1$ one finds that
$$
D.\ch_1(\Omega_1)^{\nu_1}. \cdots . \ch_1(\Omega_s)^{\nu_s}\geq 0.
$$
So $D.\ch_1(\omega_Y(S))^{n-1}=0$ implies that all those numbers must be zero, as claimed in iv).

The equality $\ch_1(\Omega_i)=\alpha\cdot \ch_1(\Omega_j)+D$ for  $i < j$ and $D\in {\rm NS}_0$ implies that
$$ 
\ch_1(\Omega_1)^{n_1}. \cdots .\ch_1(\Omega_i)^{n_i}. \cdots . \ch_1(\Omega_j)^{n_j}.
\cdots . \ch_1(\Omega_s)^{n_s} 
$$
is a multiple of an intersection number containing $\ch_1(\Omega_j)^{n_j+1}=0$. Hence by part i) it is zero, contradicting iii).

Assume there is an injection $\rho$ in vi). The poly-stability of the two sheaves implies that $S^m(\Omega_j)\cong S^m(\Omega_i) \oplus R$. So for some $\alpha>0$ one finds 
$\ch_1(\Omega_j)=\ch_1(\Omega_i) + \alpha \cdot \ch_1(R)$. 
The sheaf $R$ as a quotient of a nef sheaf has to be nef, and
$$
\ch_1(\Omega_j)^{n_j+1}.\ch_1(\omega_Y(S))^{n-1-n_j}
\geq \ch_1(\Omega_j)^{n_j}.\ch_1(\Omega_i).\ch_1(\omega_Y(S))^{n-1-n_j}
$$
is positive, in contradiction to part i).

For vii) write $\gamma$ for the semi-positive $(2,2)$ form in (\ref{eqya.2}), i.e.
$$
\gamma=2(n_i+1)\cdot \ch_2(\Omega_i)-n_i\cdot \ch_1(\Omega_i)^{2}.
$$
Since $\gamma.\ch_1(\Omega_1)^{\nu_1}. \cdots . \ch_1(\Omega_s)^{\nu_s}\geq 0$, for all
tuples $\underline{\nu}$ with $\nu_1+\cdots+\nu_s=\eta=n-2$, the equality in vii)
implies that all those intersection numbers are zero, in particular
the equation (\ref{eqya.1}) holds. On the other hand, the equation (\ref{eqya.1})
implies that $M_i$ is a complex ball, hence the $(n_i,n_i)$ current (\ref{eqya.3}) is zero.

Then $\gamma.\ch_1(\Omega_1)^{\nu_1}. \cdots . \ch_1(\Omega_s)^{\nu_s}=0$ whenever
$\nu_i \geq n_i-2$. For the other tuples $\underline{\nu}$ there is some $\nu_j>n_j$,
hence again the intersection number is zero, and the equation in vii) holds true.
\end{proof}
We will consider in this article stable and semi-stable sheaves with respect
to slopes defined by non-ample invertible sheaves $\sN$.  
\begin{definition}\label{ya.7}
Assume that $\sN$ is nef and big. 
\begin{enumerate}
\item[a.] Let $\sF$ be a coherent torsion free sheaf, and $\sG$ a subsheaf.
We call $\sF$ and $\sG$ equivalent (or $\mu_\sN$-equivalent) to $\sF$, if $\sF/\sG$ is a torsion sheaf and if $\ch_1(\sF)-\ch_1(\sG)$ is the class of an effective divisor $D$ with $D.\ch_1(\sN)^{n-1}=0$.
\item[b.] If $\theta:\sF\to \sE$ is a morphism between coherent torsion free sheaves,
the saturated image ${\rm Im}'(\theta)$ is the kernel of the map 
$$\sE \to (\sE/{\rm Im}(\theta))/_{\rm torsion}.$$ 
\item[c.] If in b) $\theta$ is injective, we call ${\rm Im}'(\theta)$ the saturated hull of $\sF$ in $\sE$.
\item[d.] We call a coherent torsion free sheaf $\sF$ weakly poly-stable, if it is equivalent
to a poly-stable subsheaf.
\end{enumerate}
\end{definition}
\begin{lemma}\label{ya.8} Assume that $\sN$ is nef and big.  
\begin{enumerate}
\item If $\sF$ is $\mu_\sN$-stable and if $\sG\subset \sF$ is a subsheaf
with $\mu_\sN(\sG)=\mu_\sN(\sF)$ then $\sF$ and $\sG$ are $\mu_\sN$-equivalent.
\item A $\mu_\sN$-stable sheaf $\sF$ is $\mu_\sN$-semi-stable.
\item If $\theta:\sF\to \sE$ is a morphism between torsion free $\mu_\sN$-semi-stable sheaves of the same slope $\mu_0$, then the saturated image ${\rm Im}'(\theta)$ is a $\mu_\sN$-semi-stable
subsheaf of slope $\mu_0$, and ${\rm Im}(\theta)$ and ${\rm Im}'(\theta)$ are $\mu_\sN$-equivalent. 
\item If in 3) the sheaf $\sF$ is weakly poly-stable, then the saturated image is weakly
poly-stable.
\end{enumerate}
\end{lemma}
\begin{proof}
Let $\sF$ be a coherent torsion free sheaf, and let $\sG$ be a subsheaf with $\rk(\sG)=\rk(\sF)$. Then
the cokernel of $\sG\hookrightarrow\sF$ is supported on an effective divisor $D$
and $\ch_1(\sF)-\ch_1(\sG)$ is an effective $\Q$-divisor $D$. Since
$\sN$ is nef, one finds $D.\ch_1(\sN)^{n-1}\geq 0$ and $\mu_\sN(\sG) \leq \mu_\sN(\sF)$.
In particular 2) holds true.

If $\sG$ is a subsheaf of a stable sheaf $\sF$ with $\mu_\sN(\sG)=\mu_\sN(\sF)$, then by Definition
$\rk(\sF)=\rk(\sG)$, hence $D.\ch_1(\sN)^{n-1}= 0$ as claimed in 1).

In 3) write $\sG={\rm Im}(\theta)$ and $\sG'={\rm Im}'(\theta)$. 
Then 
$$
\mu_0 \leq \mu_\sN(\sG) \leq \mu_\sN(\sG') \leq \mu_\sN(\sE) =\mu_0,
$$ 
and $\sG'$ is a $\mu_\sN$-semi-stable subsheaf of slope $\mu_0$. 

In d) we may replace $\sF$ by the equivalent poly-stable subsheaf. Then $\sG$ is poly-stable, hence
$\sG'$ by definition weakly poly-stable.
\end{proof}
Later $\sN$ will either be ample, hence in \ref{ya.7}, a), the divisor $D$ will be zero,
or $\sN=\omega_Y(S)$. In the second case we can make a more precise statement.
\begin{lemma}\label{ya.9} Let $\mu=\mu_{\omega_Y(S)}$, and let $\sG \hookrightarrow \sF$ be an inclusion
of semi-stable sheaves of the same slope and rank. Then $\ch_1(\sF)-\ch_1(\sG)$ lies in the subspace ${\rm NS}_0$
defined in \ref{ya.6}, v).
\end{lemma}
\begin{proof}
This follows from part 3) of \ref{ya.8} and from \ref{ya.6}, iv).
\end{proof}
We will write $C\equiv C'$ for two classes $C, C' \in {\rm NS}(Y)$ with $C-C'\in {\rm NS}_0$.  
\section{The slope of a Higgs bundle}\label{sl}
Let $\V$ be a polarized $\C$-variation of Hodge structures on $U=Y\setminus S$ with
uni-potent local monodromies and with logarithmic Higgs bundle
$$
(E=E^{1,0}\oplus E^{0,1}, \ \theta:E^{1,0}\to E^{0,1}\otimes \Omega_Y^1(\log S)).
$$
We may choose in the first part of this section $\mu=\mu_\sN$, where $\sN$ is a nef invertible sheaf (or an $\R$-divisor) on $Y$, ample with respect to $U$. Later we will use the assumptions stated in \ref{sl.5}.

Since $(E,\theta)$ is the Higgs bundle of a local system with unipotent local monodromy, 
$\ch_1(E)=\ch_1(E^{1,0})+\ch_1(E^{0,1})=0$, hence $\mu(\det(E))=0$. 
We will need two slightly different results on the behavior of slopes under filtrations of Higgs bundles.
\begin{lemma}\label{sl.1}
Let $0=E_0 \subset E_1 \subset E_2 \subset \cdots \subset E_\ell=E$
be a filtration of $E$ by saturated sub-Higgs sheaves.
\begin{enumerate}
\item[a.] If $\mu(\det(F_i))=0$ for all $i$ then
$$
\mu(E^{\doubledot})\leq {\rm Max}\{\mu(F_i^{\doubledot}); \ i=1,\ldots, \ell \},
$$
with equality if and only if
\begin{gather*}
\mu(F_1^{\doubledot})=\mu(F_2^{\doubledot})=\cdots = \mu(F_\ell^{\doubledot})
\mbox{ \ \ and \ \ } \\
\rk(F^{1,0}_i)\cdot\rk(F^{0,1}_1)=
\rk(F^{0,1}_i)\cdot\rk(F^{1,0}_1).
\end{gather*}
\item[b.] Assume that $\mu(\det(E))=0$, that $\mu(\det(E_i))\leq 0$ and 
that $F_i^{0,1} \neq 0$, for $i=1,\ldots, \ell$. If
$$
0> \mu(E_1^{0,1}) > \mu(F_2^{0,1}) >  \ldots  > \mu(F^{0,1}_\ell)
$$
one has $\mu(E^{\doubledot}) \leq {\rm Max}\{\mu(F_i^{\doubledot}); i=1,\ldots, \ell \}$,
and the equality implies that $\mu(\det(E_i))= 0$, for $i=1,\ldots,\ell$.
\end{enumerate}
\end{lemma}
\begin{proof} Let us write $\s(E)$ and $\s(F_i)$ for $\mu(E^{\doubledot})$ and $\mu(F^{\doubledot}_i)$. 
Part a) can be easily shown by induction on $\ell$, whereas the assumptions made in b) prevent a similar argument. We first introduce some correction terms allowing to handle both parts at once.

\begin{itemize}
\item $c_i=\mu(\det(F_i))$ (hence $c_i=0$ in a) and $c_1+\cdots+c_\ell=0$ in b)).\vspace{.2cm}
\item $\mu_i^{p,q}=\mu(F_i^{p,q})$, for $(p,q)=(1,0)$ or $(0,1)$.\vspace{.2cm}
\item $r_i^{1,0}=\rk(F_i^{1,0})$.\vspace{.2cm}
\item $\displaystyle r_i^{0,1}=\rk(F_i^{0,1})-\frac{c_i}{\mu_i^{0,1}}$ in b) and $\displaystyle r_i^{0,1}=\rk(F_i^{0,1})$ in a).
\end{itemize}

\begin{claim}\label{sl.2} One has
\begin{enumerate}
\item[I.] $\mu_i^{0,1}\cdot r_i^{0,1}=-\mu_i^{1,0}\cdot r_i^{1,0}$, and in particular
$r_i^{0,1}>0$.\vspace{.2cm}
\item[II.]  \ \hspace*{\fill} $r_1^{0,1}+ \cdots + r_\ell^{0,1} \leq \rk(E^{0,1}).$
 \hspace*{\fill} \ \vspace{.2cm}
\item[III.] \ \hspace*{\fill} $\displaystyle
\mu(E^{0,1})\geq \frac{\sum_{i=1}^\ell \mu_i^{0,1}\cdot r_i^{0,1}}{\sum_{i=1}^\ell
r_i^{0,1}}.$ \hspace*{\fill} \ \vspace{.2cm}
\item[IV.] \ \hspace*{\fill} $\displaystyle
\s(E) \leq \frac{\sum_{i=1}^\ell \mu_i^{1,0}\cdot r_i^{1,0}}{\sum_{i=1}^\ell
r_i^{1,0}}-\frac{\sum_{i=1}^\ell \mu_i^{0,1}\cdot r_i^{0,1}}{\sum_{i=1}^\ell
r_i^{0,1}}$ \hspace*{\fill} \ \\[.2cm]
with equality if and only if  $c_1=c_2=\cdots = c_\ell=0$.
\end{enumerate}
\end{claim}
\begin{proof} All this is obvious under the assumption a). 
For example, I. is just saying that $c_i=\mu^{1,0}_i\cdot r^{1,0}_i+\mu^{0,1}_i\cdot r^{0,1}_i$.
and IV. is the definition of $\s(E)$.

Under the assumptions made in b), part I. follows from
$$
c_i=\mu_1^{1,0}\cdot r_i^{1,0}+\mu_i^{0,1}\cdot \rk(F_i^{0,1})
$$
and from the choice of $r_i^{0,1}$. By assumption, for all $r>0$
$$
\mu(\det(E_r))=\sum_{i=1}^r c_i \leq 0
$$
with equality for $r=\ell$. Since $\displaystyle\frac{-1}{\mu^{0,1}_{\ell}}>0$ and
$\displaystyle\frac{\mu^{0,1}_{i}-\mu^{0,1}_{i+1}}
{\mu^{0,1}_{i}\cdot \mu^{0,1}_{i+1}}>0$ one finds
\begin{multline}\label{eqsl.1} 
\Big( \sum_{i=1}^\ell r_i^{0,1}\Big) - \rk(E^{0,1})=
\Big( \sum_{i=1}^\ell r_i^{0,1} - \rk(F_i^{0,1})\big)=
\sum_{i=1}^\ell \frac{-c_i}{\mu_i^{0,1}}=\\
\frac{-1}{\mu^{0,1}_{\ell}}\cdot\Big(
\sum_{j=1}^\ell c_i \Big) +
\sum_{i=1}^{\ell-1}\frac{\mu^{0,1}_{i}-\mu^{0,1}_{i+1}}
{\mu^{0,1}_{i}\cdot \mu^{0,1}_{i+1}}\cdot\Big(\sum_{j=1}^i c_i \Big) \leq 0,
\end{multline}
as claimed in II. Finally,
$$
\mu(E^{0,1})= \frac{1}{\rk(E^{0,1})}\cdot\sum_{i=1}^\ell \mu_i^{0,1}\cdot\rk(F_i^{0,1})<0.
$$
Using II. one finds
\begin{multline}\label{eqsl.2}
\mu(E^{0,1})\geq \frac{\sum_{i=1}^\ell \mu_i^{0,1}\cdot\rk(F_i^{0,1})}{\sum_{i=1}^\ell r_i^{0,1}}
=\\
\frac{\sum_{i=1}^\ell \mu_i^{0,1}\cdot r_i^{0,1}}{\sum_{i=1}^\ell r_i^{0,1}}+\frac{\sum_{i=1}^\ell
c_i}{\sum_{i=1}^\ell r_i^{0,1}} = \frac{\sum_{i=1}^\ell \mu_i^{0,1}\cdot r_i^{0,1}}{\sum_{i=1}^\ell
r_i^{0,1}},
\end{multline}
hence III and IV. If IV or equivalently III are equalities 
the same holds true for (\ref{eqsl.2}) and hence for (\ref{eqsl.1}).
Obviously the latter implies that $\sum_{j=1}^i c_j=0$ for all $i$.
\end{proof}
Let us write
$$
s_r=\frac{\sum_{i=1}^r \mu_i^{1,0}\cdot r_i^{1,0}}{\sum_{i=1}^r
r_i^{1,0}}-\frac{\sum_{i=1}^r \mu_i^{0,1}\cdot r_i^{0,1}}{\sum_{i=1}^r
r_i^{0,1}}.
$$
\begin{claim}\label{sl.3}
One has 
$$s_m \leq {\rm Max}\{s_{m-1}, \mu_m^{1,0}-\mu_m^{0,1}\}={\rm Max}\{s_{m-1}, \s(F_m)\},
$$
with equality if and only if $s_{m-1}= \s(F_m)$ and
$$
r_m^{1,0}\cdot \sum_{i=1}^m r_i^{0,1}=r_m^{0,1}\cdot \sum_{i=1}^m r_i^{1,0}.
$$
\end{claim}
Obviously \ref{sl.3} implies Lemma \ref{sl.1} a) and b. In fact,
\begin{multline*}
\s(E) \leq s_\ell \leq {\rm Max}\{s_{\ell-1}, \s(F_\ell)\}\leq {\rm Max}\{s_{\ell-2},
\s(F_{\ell-1}),
\s(F_\ell)\} \leq \\
\cdots \leq {\rm Max}\{\s(F_i); i=1,\ldots, \ell \}.
\end{multline*}
The equality implies  
$\s(F_\ell)=s_{\ell-1}=\s(F_{\ell-1})=s_{\ell-2}= \cdots =\s(F_{1})$, and
\begin{gather*}
r_\ell^{1,0}\cdot r_{\ell-1}^{0,1} \cdot \sum_{i=1}^\ell \frac{r_i^{0,1}}{r_{\ell-1}^{0,1}}
=r_\ell^{1,0}\cdot \sum_{i=1}^\ell r_i^{0,1}=\hspace{2cm}\\
\hspace{2cm} r_\ell^{0,1}\cdot \sum_{i=1}^\ell r_i^{1,0}= r_\ell^{1,0}\cdot \sum_{i=1}^\ell r_i^{0,1}=
r_\ell^{0,1}\cdot r_{\ell-1}^{1,0} \sum_{i=1}^\ell \frac{r_i^{1,0}}{r_{\ell-1}^{1,0}}\\
\cdots \cdots \cdots\\
\cdots \cdots \cdots\\
r_3^{1,0}\cdot r_2^{0,1}(1+\frac{r_1^{0,1}}{r_2^{0,1}})=r_3^{1,0}\cdot (r_2^{0,1}+r_1^{0,1})=\hspace{2cm}\\
\hspace{2cm}r_3^{0,1}\cdot (r_2^{1,0}+r_1^{1,0})=
r_3^{0,1}\cdot r_2^{1,0}(1+\frac{r_1^{1,0}}{r_2^{1,0}})\\
\mbox{and \ \ }r_2^{1,0}\cdot r_1^{0,1}=r_2^{0,1}\cdot r_1^{1,0}.
\end{gather*}
One finds the last condition stated in a).
In case b) we have seen already in \ref{sl.2}, IV), that the equality $\s(E) = s_\ell $ implies 
$c_1=c_2=\cdots=c_\ell=0$, as claimed in \ref{sl.1}, b).\\[.3cm]
\noindent
{\it Proof of Claim \ref{sl.3}.} \
Let us first handle the case $m=2$. We write $\s_1=\s(F_1), \ \s_2=\s(F_2)$ and
$s={\rm Max}\{s_1, \s(F_2)\}={\rm Max}\{\s_1 , \s_2 \}$. Let us choose
$$
h=(r_1^{1,0}+r_2^{1,0})\cdot (r_1^{1,0}+r_1^{0,1})\cdot
(r_2^{1,0}+r_2^{0,1})\cdot (r_1^{0,1}+r_2^{0,1}) >0
$$
and
$$
f=\big(r_1^{1,0}\cdot r_1^{0,1}\cdot(r_2^{1,0}+r_2^{0,1})+
r_2^{1,0}\cdot r_2^{0,1}\cdot(r_1^{1,0}+r_1^{0,1})\big)\cdot
\big(r_1^{1,0}+r_2^{1,0}+r_1^{0,1}+r_2^{0,1}\big).
$$
It is an easy exercise (by hand or using any computer algebra program) to see that $f=h-g^2$ for
$g=\big(r_1^{1,0}\cdot r_2^{0,1}-r_2^{1,0}\cdot r_1^{0,1}\big)$.

The Property I) implies that
$$
\mu_i^{1,0} = (\mu_i^{1,0}-\mu_i^{0,1})\cdot\frac{r_i^{0,1}}{r_i^{1,0}+r_i^{0,1}}=
\s(F_i)\cdot\frac{r_i^{0,1}}{r_i^{1,0}+r_i^{0,1}}
$$
and $- \mu_i^{0,1}\cdot r_i^{0,1} =\mu_i^{1,0}\cdot r_i^{1,0}$. Then 
\begin{multline*}
s_2=
\frac{\mu_1^{1,0}\cdot r_1^{1,0}+\mu_2^{1,0}\cdot r_2^{1,0} }{
r_1^{1,0}+r_2^{1,0}}+\frac{\mu_1^{1,0}\cdot r_1^{1,0}+\mu_2^{1,0}\cdot r_2^{1,0}}{
r_1^{0,1}+r_2^{0,1}}=\\
\frac{\s_1\cdot\frac{r_1^{1,0} \cdot r_1^{0,1}}{r_1^{1,0}+r_1^{0,1}}+\s_2\cdot\frac{r_2^{1,0}
\cdot r_2^{0,1}}{r_2^{1,0}+r_2^{0,1}}}{r_1^{1,0}+r_2^{1,0}}+
\frac{\s_1\cdot\frac{r_1^{1,0} \cdot r_1^{0,1}}{r_1^{1,0}+r_1^{0,1}}+
\s_2\cdot\frac{r_2^{1,0} \cdot r_2^{0,1}}{r_2^{1,0}+r_2^{0,1}}}
{r_1^{0,1}+r_2^{0,1}}=\\
\Big(\frac{\s_1\cdot r_1^{1,0}\cdot r_1^{0,1}}{(r_1^{1,0}+r_2^{1,0})\cdot
(r_1^{1,0}+r_1^{0,1})}+
\frac{\s_2\cdot r_2^{1,0}\cdot r_2^{0,1}}{(r_1^{1,0}+r_2^{1,0})\cdot
(r_2^{1,0}+r_2^{0,1})}\Big)\cdot\\
\frac{r_1^{1,0}+r_2^{1,0}+r_1^{0,1}+r_2^{0,1}}{r_1^{0,1}+r_2^{0,1}},
\end{multline*}
and $s_2 \leq s\cdot\frac{f}{h}=s \cdot (1-\frac{g^2}{h})\leq s.$
If $s_2 = s$ the polynomial $g$ must be zero, and $\s(F_1)=\s(F_2)$.

For $m>2$ we argue by induction. By definition one can write $s_r$ as
$$
\frac{\big(\sum_{i=1}^{r-1} \mu_i^{1,0}\cdot r_i^{1,0}\big) + \mu_r^{1,0} \cdot r_r^{1,0}}
{\sum_{i=1}^r r_i^{1,0}}-
\frac{\big( \sum_{i=1}^{r-1} \mu_i^{0,1}\cdot r_i^{0,1}\big) + \mu_r^{0,1} \cdot r_r^{0,1}}
{\sum_{i=1}^r r_i^{0,1}}.
$$
So writing ${r'_2}^{p,q}=r_m^{p,q}$, ${\mu'_2}^{p,q}=\mu_m^{p,q}$, ${\s'}_2=\s(F_m)$,
$$
{r'_1}^{p,q}=\sum_{i=1}^{m-1}r_i^{p,q}, \mbox{ \ \ and \ \ }{\mu'_1}^{p,q}=\frac{
\sum_{i=1}^{r-1} \mu_i^{p,q}\cdot r_i^{p,q}}{\sum_{i=1}^{r-1} r_i^{p,q}}
$$
one gets
$$
s'_2=s_m=
\frac{{\mu'}_1^{1,0}\cdot {r'}_1^{1,0}+{\mu'}_2^{1,0}\cdot {r'}_2^{1,0} }{
{r'}_1^{1,0}+{r'}_2^{1,0}}+\frac{{\mu'}_1^{1,0}\cdot {r'}_1^{1,0}+{\mu'}_2^{1,0}\cdot
{r'}_2^{1,0}}{{r'}_1^{0,1}+{r'}_2^{0,1}}.
$$
So repeating the argument for $m=2$ with $'$ added, we obtain \ref{sl.3} for all $m$.
\end{proof}

We will frequently use Simpson's correspondence for sub-Higgs bundles of a given Higgs bundle of a variation of Hodge structures. Since we do not require $\sN$ to be ample, we have to work with saturated subsheaves
$G \subset E$, i.e. with subsheaves such that $E/G$ is torsion free.

\begin{proposition}\label{sl.4} Let $\sN$ be nef and ample with respect to $U$. 
Let $E$ be a logarithmic Higgs bundle induced by a $\C$-variation of Hodge structures $\V$
on $U$ with uni-potent local monodromy. If $G\subset E$ is a sub-Higgs sheaf then
for all $n-1\geq m \geq 0$ and for all ample invertible sheaves $\sH$ on $Y$ one has
\begin{equation}\label{eqsl.3} 
\ch_1(G).\ch_1(\sN)^{n-m-1}.\ch_1(\sH)^m \leq 0.
\end{equation}
Moreover, if $G\subset E$ is saturated the following conditions are equivalent:
\begin{enumerate}
\item For some $m\geq 0$ and for all ample invertible sheaves $\sH$ the equality holds in (\ref{eqsl.3}).
\item For all $m$ and for all ample invertible sheaves $\sH$ the equality holds in (\ref{eqsl.3}).
\item $G$ is induced by a local sub-system of $\V$.
\end{enumerate}
In particular, if one of the conditions 1)--3) holds true, $\sG$ is a logarithmic Higgs bundle
and a direct factor of $E$. 
\end{proposition}
\begin{proof}
Consider for $r=\rk G$ the rank one sub-Higgs sheaf
$$
(\det G,0)=\bigwedge^r(G,\theta) \subset \bigwedge^r(E,\theta).
$$
The curvature of the Hodge metric $h$ on $\det G$ is negative semi-definite, and
the Chern form $\ch_1(G,h)$ represents the Chern class of $\det(G)$.
So one obtains (\ref{eqsl.3}) for $m=n-1$. Since $\sN$ is in the closure of the ample cone,
one obtains (\ref{eqsl.3}) for all $m$.

Assume now, that $G \subset E$ is saturated. Obviously 3) implies 2) and 2) implies 1).
If (\ref{eqsl.3}) is an equality for some $m\leq n-1$ and for all $\sH$, then it is an equality for $m=0$, since $\sN$ lies in the closure of the ample cone. The invertible sheaf $\det G$ is a saturated subsheaf of $\bigwedge^r(E)$. Replacing $Y$ by some blowing up with centers outside of $U$, the polarization $\sN$ by its pullback and $\sG$ by the saturated hull of the pullback, the equality (\ref{eqsl.3}) remains true. Hence we may assume that the inclusion $\det(G)\hookrightarrow \bigwedge^r(E)$ splits locally. Moreover, 
since $\ch_1(\sN)$ is ample with respect to $U$, one can choose the blowing up such that $\ch_1(\sN)$
is the sum $\alpha\cdot\ch_1(\sH)+\beta\cdot D$, where $\alpha$ and $\beta$ are positive real numbers, 
where $D$ is an effective $\Q$-divisor, and where $\sH$ is ample.

We will show by induction on $m$ that (\ref{eqsl.3}) holds for $\sH$ and for all
$m$. For $0\leq m_0 \leq n-2$ write $\ch_1(G).\ch_1(\sN)^{n-m_0-1}.\ch_1(\sH)^{m_0}$ as
$$
\alpha\cdot\ch_1(G).\ch_1(\sN)^{n-m_0-2}.\ch_1(\sH)^{m_0+1} + \\ \beta\cdot\ch_1(G).\ch_1(\sN)^{n-m_0-2}.\ch_1(\sH)^{m_0}.D.
$$
By \cite{Zu} none of the terms can be positive. If (\ref{eqsl.3}) is an equality for
$m=m_0$ both terms must be zero, and one obtains (\ref{eqsl.3}) for $m=m_0+1$.
By induction (\ref{eqsl.3}) holds for all $m$, in particular for $m=n-1$. \cite{Sim} implies that $G$ is induced by a local system on $U$.
\end{proof}
In the sequel we will need the compatibility of ``poly-stability''
with tensor products, as stated in condition ($*$) in the Set-up \ref{in.3}.
\begin{assumption}\label{sl.5}
Either $\mu=\mu_{\omega_Y(S)}$ and $\delta=\delta_{\omega_Y(S)}$ or
$\mu=\mu_\sN$ and $\delta=\delta_\sN$ for an ample invertible sheaf
(or an ample $\R$-divisor). In the first case we assume that the assumptions made in Set-up \ref{in.3} hold. In the second case we assume that $\omega_Y(S)$ is nef and ample with respect to $U$ and that $S^m(\Omega_Y^1(\log S))$ is $\mu_\sN$-poly-stable for all $m>0$, a condition which automatically holds in the first case by Theorem \ref{ya.4}.
\end{assumption}
Note that under those assumptions $\mu(\Omega^1_Y(\log S))>0$.
\begin{proposition}\label{sl.6} Assume that \ref{sl.5} holds. Let $(E=E^{1,0}\oplus E^{0,1}, \theta)$
be the logarithmic Higgs bundle of a $\C$-variation of Hodge structures $\V$ on $U=Y\setminus S$
of weight one, and assume that $\V$ has no unitary direct factor. Then:
\begin{enumerate}
\item[a.] $\mu(\V)\leq \mu(\Omega^1_Y(\log S))$.
\item[b.] The equality $\mu(\V) = \mu(\Omega^1_Y(\log S))$ implies that
$E^{1,0}$ and $E^{0,1}$ are both semi-stable.
\end{enumerate}
\end{proposition}
\begin{proof}
\ref{sl.1}, a), allows us to assume that $\V$ is irreducible.
Let
\begin{gather*}
0=E^{1,0}_0 \subset E^{1,0}_1 \subset \cdots \subset E^{1,0}_\ell=E^{1,0}\\
\mbox{and \ \ }
0=E^{0,1}_0 \subset E^{0,1}_1 \subset \cdots \subset E^{0,1}_{\ell'}=E^{0,1}
\end{gather*}
be the saturated Harder-Narasimhan filtrations of $E^{1,0}$ and $E^{0,1}$, respectively.
Replacing $\V$ by its dual, if necessary, we may assume that $\ell \leq \ell'$.
Since $(E^{0,1}_1,0)$ is a sub-Higgs sheaf,
and since $(E^{1,0}_\ell/E^{1,0}_{\ell-1},0)$ a quotient Higgs sheaf one finds
$\mu(E^{0,1}_1)<0<\mu(E^{1,0}_\ell/E^{1,0}_{\ell-1})$ and
$$
\mu(E_1^{1,0}) >\cdots >  \mu(E^{1,0}_\ell/E^{1,0}_{\ell-1}) >0 >
\mu(E^{0,1}_{1}) > \cdots > \mu(E^{0,1}_{\ell'}/E^{0,1}_{\ell'-1}).
$$
In particular by Proposition \ref{sl.4} $(E_1^{1,0},0)$ can not occur as a sub-Higgs sheaf 
and $\theta(E_1^{1,0})\neq 0$. The morphism
$$
E^{1,0}\otimes T^1_Y(-\log S) \>>> E^{0,1}_{\ell'}/E^{0,1}_{\ell'-1},
$$
induced by $\theta$ is non zero; otherwise $(E^{0,1}_{\ell'}/E^{0,1}_{\ell'-1},0)$
would be a quotient Higgs sheaf of negative degree, contradicting again the inequality (\ref{eqsl.3}) in Proposition \ref{sl.4}. 

Choose two sequences
$$
0=j_0 < j_1 < \cdots < j_r=\ell \mbox{ \ \ and \ \ } 0=j'_0 < j'_1 < \cdots < j'_r=\ell',
$$
in the following way:\\[.1cm]
Assume one has defined $j_{m-1}$ and $j'_{m-1}$. Then
$j'_m$ is the minimal number with $\theta(E^{1,0}_{j_{m-1}+1}) \subset E^{0,1}_{j'_m}\otimes
\Omega^1_Y(\log S)$, and $j_m$ is the maximum of all $j$ with
$\theta(E^{1,0}_{j}) \subset E^{0,1}_{j'_m}\otimes \Omega^1_Y(\log S)$.

In different terms, one has
\begin{gather*}
\theta(E^{1,0}_{j_{m}}) \subset E^{0,1}_{j'_m}\otimes \Omega^1_Y(\log S), \hspace{1cm}
\theta(E^{1,0}_{j_{m}}) \not\subset E^{0,1}_{j'_m-1}\otimes\Omega^1_Y(\log S)\\
\mbox{ and \ \ }
\theta(E^{1,0}_{j_{m-1}+1}) \not\subset E^{0,1}_{j'_m-1}\otimes\Omega^1_Y(\log S).
\end{gather*}
One has a non-trivial morphisms
\begin{gather*}
E^{1,0}_{j_{m-1}+1}/E^{1,0}_{j_{m-1}}\>>>
E^{0,1}_{j'_m}/E^{0,1}_{j'_m-1}\otimes \Omega^1_Y(\log S) \mbox{ \ \ and \ \ }\\
E^{1,0}_{j_{m-1}+1}/E^{1,0}_{j_{m-1}}\otimes T^1_Y(-\log S) \>>>
E^{0,1}_{j'_m}/E^{0,1}_{j'_m-1}.
\end{gather*}
Hence for a stable sheaf $\sC$ with $\mu(\sC)=\mu(E^{1,0}_{j_{m-1}+1}/E^{1,0}_{j_{m-1}})$
one obtains a non-zero morphism
$$
\sC \otimes T^1_Y(-\log S) \>>>
E^{0,1}_{j'_m}/E^{0,1}_{j'_m-1}.
$$
Since $\sC \otimes T^1_Y(-\log S)$ is semi-stable, this implies that
$$
\mu(\sC)=\mu(E^{1,0}_{j_{m-1}+1}/E^{1,0}_{j_{m-1}})-\mu(\Omega_Y^1(\log S)) \leq
\mu(E^{0,1}_{j'_m}/E^{0,1}_{j'_m-1}).
$$
Consider the filtration of $E$ by saturated sub-Higgs sheaves
$$
E_{1}=E_{j_1}^{1,0}\oplus E_{j'_1}^{0,1} \subset E_{2}=E_{j_2}^{1,0}\oplus E_{j'_2}^{0,1}\subset
\cdots \subset E_{r}=E_{j_r}^{1,0}\oplus E_{j'_r}^{0,1}=E,
$$
with successive quotients $F_m=F_m^{1,0}\oplus F_m^{0,1}$ for
$$
F_{m}^{1,0}=E_{j_m}^{1,0}/E_{j_{m-1}}^{1,0} \mbox{ \ \ and \ \ } F_m^{0,1}=
E_{j'_m}^{0,1}/E_{j'_{m-1}}^{0,1}.
$$
Of course, $\mu(E_i)\leq 0$, and
$$
\mu(E_{j_{m-1}+1}^{1,0}/E_{j_{m-1}}^{1,0} ) \geq \mu(F_m^{1,0}) > 0 > \mu(F_m^{0,1}) \geq
\mu(E_{j'_m}^{0,1}/E_{j'_{m}-1}^{0,1}).
$$
Hence
$$
0 > \mu(F_1^{0,1}) > \mu(F_2^{0,1}) > \cdots > \mu(F_r^{0,1}),
$$
and
\begin{equation}\label{eqsl.4}
\mu(F_m^{1,0}) \leq \mu(E_{j_{m-1}+1}^{1,0}/E_{j_{m-1}}^{1,0} )-
\mu(E_{j'_m}^{0,1}/E_{j'_{m}-1}^{0,1})\leq \mu(\Omega_Y^1(\log S)).
\end{equation}
Lemma \ref{sl.1}, b), implies that
\begin{equation}\label{eqsl.5}
\mu(E^{\doubledot})\leq \mu(\Omega_Y^1(\log S)).
\end{equation}
If (\ref{eqsl.5}) is an equality, then $\mu(\det(E_1))=0$, and by Proposition \ref{sl.4} $E_1$ must correspond to a local sub-system. Since we assumed $\V$ to be irreducible one finds $r=1$. Moreover,
(\ref{eqsl.4}) must be an equality, which implies that $\ell=\ell'=1$, and both, $E^{1,0}$ and
$E^{0,1}$ are semi-stable.
\end{proof}
\begin{proof}[Proof of Theorem \ref{in.1}]
Let $Y$ be a projective manifold, $S$ a reduced normal crossing divisor and let
$f:V\to U$ be a smooth family of Abelian varieties, satisfying the assumptions made in
Theorem \ref{in.1}. Let $\V$ be a $\C$-sub-variation of Hodge structures of 
$R^1f_*\C_V$ without a unitary local sub-system. Then the Higgs field
$\theta:E^{1,0}\to E^{0,1}\otimes \Omega_Y^1(\log S)$ is injective and the Arakelov inequality (\ref{eqin.2}) in Theorem \ref{in.1}, as well as the interpretation of the equality (\ref{eqin.3}), follow directly from
\ref{ya.4}, a) and \ref{sl.6}. 
\end{proof}
Let us verify what we stated in Remark \ref{in.2}. Consider
a real sub-variation $\V'$ of Hodge structures and the largest unitary sub-system $\U$ of $\V'$. Let us write
$\V'=\V\oplus\U$, and $(E'^{1,0}\oplus E'^{0,1},\theta')$ and $(E^{1,0}\oplus E^{0,1},\theta)$
for the Higgs bundles of $\V'$ and $\V$. Since $\U$ is invariant under complex conjugation, the same holds true for 
$\V$.Then $\rk(E'^{1,0})=\rk(E'^{0,1})$ and $\rk(E^{1,0})=\rk(E^{0,1})$. Since 
$\ch_1(E'^{p,1-p})=\ch_1(E^{p,1-p})$ one obtains by Theorem \ref{in.1} 
\begin{multline*}
\mu(\V')=
\frac{\ch_1(E'^{1,0})- \ch_1(E'^{0,1})}{\rk(E'^{1,0})}.\ch_1(\omega_Y(S))^{n-1}
\leq\\ \frac{\ch_1(E^{1,0})-\ch_1(E^{0,1})}{\rk(E^{1,0})}.\ch_1(\omega_Y(S))^{n-1}= \mu(\V)\leq
\mu(\Omega_Y^1(\log S)),
\end{multline*}
and the Arakelov inequality holds for $\V'$. It can only be an equality if 
$\V=\V'$, hence if $\U=0$.\vspace{.2cm}

Recall that the socle $\sS'_1(E)$ of a semi-stable sheaf $E$ is the unique largest poly-stable
subsheaf of slope $\mu(E)$ (see \cite[1.5.5]{HL}). If $\mu$ is the slope with respect to
an ample invertible sheaf, $\sS'_1(E)$ is saturated in $E$. In general one chooses
$\sS_1(E)$ as the saturated hull of $\sS'_1(E)$. Doing so, one perhaps looses the
poly-stability, but one still has the weak poly-stability, as defined in \ref{ya.7}.

So $\sS_1(E)$ is a maximal weakly poly-stable subsheaf of $E$. Applying the same construction to
$E/\sS_1(E)$ one finds the socle-filtration
$$
0=\sS_0(E) \varsubsetneqq \sS_1(E) \varsubsetneqq \cdots \varsubsetneqq \sS_{\rho(E)}(E)=E
$$
such that $\sS_i(E)/\sS_{i-1}(E)$ is the saturation of the socle of $E/\sS_{i-1}(E)$. In particular
the graded sheaf $\gr_\sS(E)$ with respect to this filtration is weakly poly-stable with slope $\mu(E)$. 

\begin{lemma}\label{sl.7} Keeping the assumptions made in \ref{sl.5}, let $\Omega$ be a poly-stable sheaf and let $E$ and $F$ be semi-stable locally free sheaves, of slopes
$\mu(\Omega)$, $\mu(E)$, and $\mu(F)$, respectively. Assume that $\mu(E)=\mu(F)+\mu(\Omega)$,
and consider a morphism $\theta:E \to F\otimes \Omega$. Then:
\begin{enumerate}
\item[a.] If $\theta$ is injective, it respects the socle filtration, i.e. for all $i$
$$
\theta(\sS_i(E)) \subset \sS_i(F) \otimes \Omega \mbox{ \ \ and \ \ }
\theta^{-1}(\sS_i(F) \otimes \Omega) = \sS_i(E)).
$$
\item[b.] In a) the induced morphism $\mathfrak{gr}_\sS E \to (\mathfrak{gr}_\sS F)\otimes \Omega$
is again injective.
\item[c.] If $E$ and $F$ are weakly poly-stable for $T=\Omega^\vee$ the saturated image of
$$
\theta':E\otimes T \>\theta\otimes{\rm id}>> F\otimes \Omega\otimes T \>{\rm id}\otimes {\rm tr}>> F.
$$ 
is a weakly poly-stable subsheaf of $F$ of slope $\mu(F)$.
\end{enumerate}
\end{lemma}
\begin{proof} For c) consider a poly-stable subsheaf $E'$ of $E$ of slope $\mu(E)$ of maximal rank. Then $E'\otimes T$ is poly-stable of slope $\mu(E)-\mu(\Omega)=\mu(F)$, hence its image
in $F$ as well. Then by definition the saturated image is weakly poly-stable.
 
To prove a) we proceed by induction on the length of the socle filtration $\rho(E)$. 
For $i=0$, in particular for $\rho(E)=0$, the first inclusion is obvious 
and the second equality is just the injectivity of $\theta$. 

The sheaf $\sS_1(E)$ is weakly poly-stable, and $\theta'$ is a morphism between semi-stable sheaves of the same slope. Then $\theta'(\sS_1(E)\otimes T)$ is weakly poly-stable
by part c), hence contained in $\sS_1(F)$ . This implies that $\theta(\sS_1(E)) \subset \sS_1(F)\otimes \Omega$ and that $\sS_1(E)\subset \theta^{-1}(\sS_1(F)\otimes \Omega)$.

The injectivity of $\theta$ implies that $\theta^{-1}(\sS_1(F) \otimes \Omega)$ is
weakly poly-stable, hence contained in $\sS_1(E)$. 

The first inclusion shows that $E \to (F/_{\sS_1(F)})\otimes \Omega$ factors through
$$
\tilde\theta: E/_{\sS_1(E)} \>>> (F/_{\sS_1(F)})\otimes \Omega
$$
and the second equality says that $\tilde\theta$ is again injective.
Since $\rho(E/_{\sS_1(E)}) = \rho(E)-1$, we obtain a). 
Part b) follows directly from a). 
\end{proof}
\begin{corollary}\label{sl.8}
Under the assumptions made in \ref{sl.6}, assume that\\
$\mu(E^\doubledot)=\mu(\Omega^1_Y(\log S))$.  Then there is a filtration
$$0=F_0^{p,1-p}\subset F_1^{p,1-p}\subset \cdots \subset F_\nu ^{p,1-p}=E^{p,1-p}$$ 
with
\begin{enumerate}
\item[i.] $\theta(F_\eta^{1,0})\subset F_\eta^{0,1}\otimes \Omega^1_Y(\log S)$ and
$\theta^{-1}(F_\eta^{0,1}\otimes \Omega^1_Y(\log S)) = F_\eta^{1,0}$.
\item[ii.] $E_{\eta+1}^{p,1-p}=F_{\eta+1}^{p,1-p}/F_\eta^{p,1-p}$ is weakly poly-stable of slope
$\mu(E^{p,1-p})$.
\end{enumerate}
\end{corollary}
\begin{proof}
Since we assumed that the local system has no unitary factor, $\theta$ must be injective,
and \ref{sl.7} applies.
\end{proof}

\section{Stability of Hodge bundles}\label{st}

Before stating the main result of this section in \ref{st.4}, let us recall some 
facts about Chern classes.
\begin{lemma}\label{st.1}
There exist non-negative rational numbers $a_1$ and $a_2$, depending on $m$, $r$ and on $n=\dim(Y)$,
such that for all locally free sheaves $\sF$ on $Y$ of rank $r$ one has
\begin{equation}\label{eqst.1}
\ch_2(S^m(\sF)) \equiv a_1\cdot \ch_2(\sF) + a_2\cdot (\ch_1(\sF)^2 -\ch_2(\sF)).
\end{equation}
Moreover, if $m>1$, $r>1$ and $n>1$ one has $a_2 >0$ and $a_1 - a_2 > 0$.
\end{lemma}
\begin{proof}
Recall that in degree $2$ there are two Schur polynomials,
$$
s_{(2,0)}=\ch_2 \mbox{ \ \ and \ \ } s_{(1,1)}=\ch_1^2 - \ch_2,
$$
and, as explained in \cite[II, Chapter 8]{L}, they generate the cone of 
degree $2$ positive polynomials. Since $\ch_2(S^m)$ is a positive polynomial,
i.e. since $\ch_2(S^m(\sF))$ is positive for $\sF$ an ample locally free sheaf, one finds 
$a_1, a_2 \geq 0$ for which (\ref{eqst.1}) holds.
 
On the other hand, the equation (\ref{eqst.1}) is a very special case of universal 
relations between Schur polynomials of tensor bundles, studied by Pragacz 
(see for example \cite[II. p. 121]{L}). In particular the constants occurring can be 
chosen to be independent of the 
bundle $\sF$. To verify $a_2>0$ and $a_1-a_2>0$ we can consider special bundles.

For example, if $m$, $n$ and $r$ are strictly larger than one, for 
$$
\sF=\sH\oplus \bigoplus^{r-1} \sO_Y
$$ 
one finds $\ch_2(\sF)=0$, $\ch_2(S^m(\sF))>0$, hence $a_2>0$.

Consider next the bundle 
$$
\sF= \sH \oplus \sH^{-1} \oplus \bigoplus^{r-2}\sO_Y.
$$
One has $\ch_1(\sF)\equiv 0$ and $\ch_2(\sF).\ch_1(\sH)^{n-2} < 0$, 
hence in order to show that $a_1 - a_2 > 0$, it suffices to show that
$$
(a_1-a_2)\ch_2(\sF).\ch_1(\sH)^{n-2}=\ch_2(S^m(\sF)).\ch_1(\sH)^{n-2} < 0.
$$
The sheaf $S^m(\sF)$ is equal to
$$
\bigoplus_{i=0}^m S^{m-i}(\sH \oplus \sH^{-1})  \otimes 
S^i(\bigoplus^{r-2}\sO_Y)
= \bigoplus^{\binom{r+i-2}{i}}\big(\bigoplus_{i=0}^m S^{m-i}(\sH \oplus \sH^{-1})\big).  
$$
Since for a direct sum of sheaves with zero first Chern class the second Chern class 
is additive, $\ch_2(S^m(\sF)).\ch_1(\sH)^{n-2} < 0$ follows from
\begin{multline*}
\ch_2(S^m(\sH \oplus \sH^{-1})).\ch_1(\sH)^{n-2}=
\sum_{j=0}^{[\frac{m}{2}]}\ch_2(\sH^{m-j} \oplus \sH^{-m+j})).\ch_1(\sH)^{n-2}=\\
\ch_1(\sH)^{n}\cdot \sum_{j=0}^{[\frac{m}{2}]}-(m-j)^2 <0.
\end{multline*}
\end{proof}
Recall that in Section \ref{sl} we considered for a semi-stable sheaf $\sF$ of slope $\mu_0$ the socle filtration
$\sS_\bullet(\sF)$, and that the direct factors of the corresponding
graded sheaf $\mathfrak{gr}_\sS\sF$ are all torsion free and weakly poly-stable of slope $\mu_0$.
Obviously one can refine the filtration $\sS_\bullet(\sF)$ to obtain a Jordan-H\"older
filtration ${\rm JH}_\bullet(\sF)$. By definition the direct factors of the graded sheaf
$\mathfrak{gr}_{{\rm JH}}\sF$ are all stable of slope $\mu_0$. One can be more precise:
\begin{lemma}\label{st.2} Let $\sB$ be a stable saturated subsheaf of
$\mathfrak{gr}_\sS\sF$ of slope $\mu_0$. 
There exists a Jordan-H\"older filtration ${\rm JH}_\bullet(\sF)$ of $\sF$, refining $\sS_\bullet(\sF)$, such that $\sB$ is a direct factor of $\mathfrak{gr}_{{\rm JH}}\sF$.
\end{lemma}
Note that $\mathfrak{gr}_\sS\sF$ contains a poly-stable subsheaf $\sP$ of slope $\mu_0$, such that the cokernel is a torsion sheaf. If $\sB'$ is one of the stable direct factors of $\sP$ we may choose in Lemma \ref{st.2} for $\sB$ the saturated hull of $\sB'$.
\begin{lemma}\label{st.3} Let $\mu=\mu_\sN$ with $\sN$ nef and big, and let $\sF$ and $\sF'$ be two locally free sheaves.
\begin{enumerate} 
\item[a.] Then 
$$
\frac{\Delta(\sF\otimes \sF')}{\rk(\sF)^2\cdot\rk(\sF')^2}=
\frac{\Delta(\sF)}{\rk(\sF)^2} + \frac{\Delta(\sF')}{\rk(\sF')^2}
$$
In particular, if $\sL$ is invertible,
$$
\Delta(\sF)=\Delta(\sF\otimes\sL) \mbox{ \ \ hence \ \ }
\delta(\sF)=\delta(\sF\otimes\sL).
$$
\item[b.] For $m>0$ one has $\Delta(S^m(\sF))=0$ if and only if 
$\Delta(\sF)=0$.
\item[c.] For $m>0$ one has $\delta(S^m(\sF))=0$ if and only if 
$\delta(\sF)=0$.
\item[d.] If $\sF$ is semi-stable of slope $\mu_0$, then the following conditions are equivalent:
\begin{enumerate}
\item[1)] $\delta(\sF)=0$.
\item[2)] $\delta(\mathfrak{gr}_{{\rm JH}}\sF)=0$, where ${\rm JH}$ is a Jordan-H\"older filtration.
\item[3)] For all stable direct factors $\sG$ of $\mathfrak{gr}_{{\rm JH}}\sF$ one has
$\delta(\sG)=0$.
\item[4)] $\delta(\mathfrak{gr}_\sS\sF)=0$, where $\sS$ is the socle filtration.
\item[5)] For all $i$ one has $\delta(\sS_{i-1}(\sF)/\sS_i(\sF))=0$.
\item[6)] For all stable subsheaves $\sG$ of $\sS_{i-1}(\sF)/\sS_i(\sF)$
of slope $\mu_0$ one has $\delta(\sG)=0$.
\end{enumerate}
\end{enumerate}
\end{lemma}
\begin{proof}
a) is well known and shown in \cite[p. 72]{HL}, for example.
In order to prove \ref{st.3}, b) and c), we may replace $Y$ by a finite covering, hence assume that
$\det(\sF)=\sL^{r}$ for some invertible sheaf $\sL$ and for $r=\rk(\sF)$. Obviously
$\ch_1(\sF\otimes \sL^{-1})=0$ and $\Delta(\sF)=0$ if and only if
$\ch_2(\sF\otimes \sL^{-1})=0$. By \ref{st.1} the latter is equivalent to 
$$
\Delta(S^m(\sF))=\ch_2(S^m(\sF\otimes \sL^{-1}))=0.
$$
For c) we use the same argument: $\delta(\sF)=0$ if and only if
$$
\ch_2(\sF\otimes \sL^{-1}).\ch_1(\sN)^{n-2}=0,
$$
and this is equivalent to $\delta(S^m(\sF))=\ch_2(S^m(\sF\otimes \sL^{-1})).\ch_1(\sN)^{n-2}=0$.

If $\sF$ in d) is semi-stable but not stable consider a stable subsheaf $\sG_1$
of slope $\mu(\sF)$ and the exact sequence
$$
0\>>> \sG_1 \>>> \sF \>>> \sG_2 \>>> 0.
$$
One may assume again that $\ch_1(\sG_2)=0$. Then $\ch_1(\sF)=\ch_1(\sG_1)$
and $\ch_2(\sF)=\ch_2(\sG_1)+\ch_2(\sG_2)$. Writing $r=\rk(\sF)$ and $r_i=\rk(\sG_i)$, one finds 
\begin{multline*}
\frac{1}{r}\cdot\Delta(\sF)=\frac{1}{r_1}\cdot\Delta(\sG_1)+\frac{1}{r_2}\cdot\Delta(\sG_2)
- \frac{r-1}{r}\cdot\ch_1(\sF)^2 + \frac{r_1-1}{r_1}\cdot\ch_1(\sG_1)^2=\\
\frac{1}{r_1}\cdot\Delta(\sG_1)+\frac{1}{r_2}\cdot\Delta(\sG_2)-\frac{r_2\cdot r_1}{r}\cdot
\big(\frac{\ch_1(\sG_1)}{r_1}\big)^2.
\end{multline*}
The generalized Hodge index formula
(see for example \cite[Variant 1.6.2]{L}) implies that
$$
(\frac{\ch_1(\sG_1)}{r_1}.\ch_1(\sN)^{n-1})^2 \geq 
((\frac{\ch_1(\sG_1)}{r_1})^2.\ch_1(\sN)^{n-2})\cdot \ch_1(\sN)^{n}.
$$
Since $\mu(\sF)=\mu(\sG_1)=\mu(\sG_2)=0$, the left hand side is zero, hence
$$
(\frac{\ch_1(\sG_1)}{r_1})^2.\ch_1(\sN)^{n-2}=0,
$$
and
$$
\frac{1}{r}\cdot\delta(\sF)=\frac{1}{r_1}\cdot\delta(\sG_1)+\frac{1}{r_2}\cdot\delta(\sG_2).
$$
By Bogomolov's inequality (see for example \cite[7.3.1]{HL}) $\delta(\sF)$ and $\delta(\sG_i)$
are non-negative, hence $\delta(\sF)=0$ if and only if both, $\delta(\sG_1)$ and $\delta(\sG_2)$
are zero. So the equivalence of the first three conditions in d) follows by induction on the number of stable direct factors in $\mathfrak{gr}_{{\rm JH}}\sF$.

The equivalence of 1) and 4) follows in the same way by induction on the length of the socle filtration. Finally the equivalence of 4), 5) and 6) is just a special case of the 
equivalence of 1), 2) and 3), applied to $\sS_{i-1}(\sF)/\sS_i(\sF)$.
\end{proof}

\begin{proposition}\label{st.4}
Under the assumptions made in \ref{sl.5} let $\V$ be an irreducible $\C$ 
variation of Hodge structures whose logarithmic Higgs field factors through 
$$
\theta:E^{1,0}\>>>  E^{0,1}\otimes \Omega \>>> E^{0,1}\otimes \Omega^1_Y(\log S),
$$
for a poly-stable subsheaf $\Omega$ of $\Omega^1_Y(\log S)$ of 
slope $\mu(\Omega^1_Y(\log S))$. Let 
$$
\Omega=\Omega_1\oplus \cdots \oplus \Omega_s,
$$
be the decomposition in a direct sum of stable sheaves,
and assume that $S^m(\Omega_i)$ remains stable for $i=1,\ldots,s$ and for all $m>0$.
Assume that
$$
\mu(\V)=\mu(E^{1,0})-\mu(E^{0,1}) = \mu(\Omega^1_Y(\log S))
$$
and that either one of the following conditions holds true:
\begin{enumerate}
\item[i.] $\delta(E^{1,0})=0$.
\item[ii.] $\Omega_i$ is invertible for $i=1,\ldots,s$.
\end{enumerate} 
Then for some $\iota \in \{1, \ldots , s\}$ the Higgs field factors through
$$
\theta:E^{1,0}\>>>  E^{0,1}\otimes \Omega_\iota \>>> E^{0,1}\otimes \Omega
\>>> E^{0,1}\otimes \Omega^1_Y(\log S),
$$
and the sheaves $E^{1,0}$ and $E^{0,1}$ are stable.
For $T_\iota=\Omega_\iota^\vee$ the Higgs field $\theta$ induces an injection
$\theta'_\iota:E^{1,0}\otimes T_\iota \to E^{0,1}$ of stable sheaves of the same rank and slope. 
Writing $\ell=\rk(E^{1,0})$, the sheaves
$$
S^\ell(E^{1,0})\otimes \det(E^{1,0})^{-1} \mbox{ \ \ and \ \ }
E^{1,0}\otimes {E^{1,0}}^\vee
$$ 
are unitary.
\end{proposition}
If the assumption $\delta(\V)={\rm Min}\{ \delta(E^{1,0}), \delta(E^{0,1})\}=0$
in Proposition \ref{in.4}, i), holds true, either $\V$ or $\V^\vee$ satisfies the assumption
i) in \ref{st.4}. So Proposition \ref{in.4} follows immediately from Proposition \ref{st.4}.

In case there exists an invertible sheaf $\sL$ on $Y$ with
$\sL^\ell=\det(E^{1,0})$ one also finds that $E^{1,0}\otimes \sL$ is unitary.
At the moment we are unable to verify the existence of such an invertible sheaf. 
Nevertheless, \ref{st.4} is sufficient to show that in Theorem \ref{in.5} 
the factor $M_i$ of $\tilde{U}$ corresponding to $\Omega_i$ is a ball, and this 
will allow in Sections \ref{nc} and \ref{co} to get hold of $\sL$ on some \'etale 
covering of $U$.

\begin{proof}[Proof of Proposition \ref{st.4}, i)] Consider the slope filtrations
$\sS_\bullet(E^{1,0})$ and $\sS_\bullet(E^{0,1})$ and a stable saturated subsheaf 
$\sF$ of $\mathfrak{gr}_\sS (E^{1,0})$ of slope $\mu_0=\mu(E^{1,0})$. By \ref{st.2}
one can refine $\sS_\bullet(E^{1,0})$ to obtain a Jordan-H\"older filtration
${\rm JH}_\bullet (E^{1,0})$ such that $\sF$ is a direct factor of the corresponding
graded sheaf. 
 
By \ref{st.3}, c) and d), for all $\nu>0$ 
$$
\delta(S^{\nu\cdot\ell}(\mathfrak{gr}_{{\rm JH}} (E^{1,0}))\otimes \det(E^{1,0})^{-\nu})
=\delta(S^{\nu\cdot\ell}(\mathfrak{gr}_{{\rm JH}} (E^{1,0})))=0.
$$
Moreover by \ref{ya.2} one has
$\mu(S^{\nu\cdot\ell}(\mathfrak{gr}_{{\rm JH}} (E^{1,0}))\otimes \det(E^{1,0})^{-\nu})=0$,
and the condition ($**$) in \ref{in.3}, or Lemma \ref{ya.3}, d), imply that the sheaf
$$
S^{\nu\cdot\ell}(\mathfrak{gr}_{{\rm JH}} (E^{1,0}))\otimes \det(E^{1,0})^{-\nu}
$$ 
is unitary. So all direct factors of $S^{\ell}(\mathfrak{gr}_{{\rm JH}} (E^{1,0}))\otimes \det(E^{1,0})^{-1}$
have a trivial first Chern class, hence
$$
\Upsilon(E^{1,0})=\frac{\ch_1(E^{1,0})}{\ell}=\frac{\ch_1(S^\ell(\sF))}{\ell\cdot\rk(S^\ell(\sF))}=
\frac{\ch_1(\sF)}{\rk(\sF)}=\Upsilon(\sF).
$$ 
Repeating this for all the stable saturated subsheaves one gets: 
\begin{claim}\label{st.5}
For saturated stable subsheaves $\sF_1$ and $\sF_2$ of $\mathfrak{gr}_\sS (E^{1,0})$
of slope $\mu(E^{1,0})$ one has $\Upsilon(\sF_1)=\Upsilon(\sF_2)$.
\end{claim}
\begin{claim}\label{st.6} 
Let $\sF$ be a saturated subsheaf of $\mathfrak{gr}_\sS (E^{1,0})$, stable of slope
$\mu_0$. Then for all stable direct factors $T_i=\Omega_i^\vee$ of
$T=\Omega^\vee$ the sheaf $\sF\otimes T_i$ is stable.
\end{claim}
\begin{proof}
Let us choose again a Jordan-H\"older filtration refining $\sS_\bullet(E^{1,0})$
such that $\sF$ occurs as a direct factor of 
$\mathfrak{gr}_{{\rm JH}} (E^{1,0})$.

The condition ($*$), or Lemma \ref{ya.3}, c), imply that $\sF\otimes T_i$ is poly-stable.
Let $\sB$ be one of its direct factors.
We have inclusions $\BP(\sB) \to \BP(\sF \otimes T_i)$ and
$\BP(\sF)\times \BP(T_i)\to \BP(\sF \otimes T_i)$.
Let us write $\pi:\BP(\sF)\times \BP(T_i)\to \BP(T_i)$ for the projection and
$$
Z=\BP(\sB) \cap \BP(\sF)\times \BP(T_i).
$$ 
Let $\sF'$ be any stable direct factor of 
$S^{\nu\cdot\ell}(\sF)$, of slope $\nu\cdot\ell\cdot\mu(E^{1,0})$.

Since $S^{\nu\cdot\ell}(\mathfrak{gr}_{{\rm JH}} (E^{1,0}))\otimes \det(E^{1,0})^{-\nu}$
is unitary, $\sF'\otimes \det(E^{1,0})^{-\nu}$ is unitary and by assumption irreducible.

The Addendum \ref{ya.5} forces the sheaf $\sF'\otimes \det(E^{1,0})^{-\nu}\otimes S^{\nu\cdot\ell}(T_i)$
to be stable. Hence the composite
$$ 
\sF'\otimes S^{\nu\cdot\ell}(T_i) \>>> S^{\nu\cdot\ell}(\sF)\otimes S^{\nu\cdot\ell}(T_i) \>>>   S^{\nu\cdot\ell}(\sB)
$$
is either injective, or zero. In different terms,
the bihomogeneous ideal of $Z$, regarded as a subscheme of $\BP(\sF)\times \BP(T_i)$,
is generated by 
$$
\sI(\nu\cdot\ell) \otimes \pi^*(\sO_{\BP(T_i)}(\nu\cdot\ell))
$$ 
for some $\nu>0$, and some sheaf of ideals $\sI$ on $\BP(\sF)$. 
Then $Z$ must be of the form $Z'\times \BP(T_i)$ for some $Z'\subset \BP(\sF)$.

Let us restrict everything to a general point $\eta={\rm Spec}(\overline{\C(Y)})$ of $Y$.
The embedding $\BP(\sB)_\eta \to \BP(\sF \otimes T_i)_\eta$ is linear,
and if $p$ is a point in $\BP(T_i)_\eta$ the same holds true for
$$
\BP(\sF)_\eta\times \{p\} \>>> \BP(\sF \otimes T_i)_\eta.
$$
Then $Z'_\eta=\BP(\sF)_\eta\times \{p\}\cap \BP(\sB)_\eta$ is a linear subspace,
and the projection 
$$
\sF_\eta \otimes {T_i}_\eta \> \tau >> \sB_\eta
$$
must have a kernel of the form $\sK_\eta\otimes {T_i}_\eta$. If $\sK_\eta\neq 0$ the non-zero morphism
$\sF \to \sB \otimes \Omega_i$ has a non-trivial kernel $\sK$. Since $\sF$ is stable and $\sB\otimes \Omega_i$ poly-stable, both of the same slope, one finds $\sK=\sF$, obviously a contradiction.
\end{proof}
Let $\mathfrak{gr}_\sS\theta: \mathfrak{gr}_\sS  E^{1,0} \to \mathfrak{gr}_\sS  E^{0,1}\otimes
\Omega$ denote the Higgs field, and let 
$$
\mathfrak{gr}_\sS\theta': \mathfrak{gr}_\sS  E^{1,0}\otimes T \>>> 
\mathfrak{gr}_\sS  E^{0,1}\otimes \Omega \otimes T\>>> 
\mathfrak{gr}_\sS  E^{0,1}
$$
be induced by the trace $\Omega \otimes T \to \sO_Y$. 
For a stable subsheaf $\sF$ of $\mathfrak{gr}_\sS  E^{1,0}$ of slope $\mu_0$
we write $\theta_{\sF,i}$ for the restriction of $\mathfrak{gr}_\sS\theta'$
to $\sF\otimes T_i$. The image of $\mathfrak{gr}_\sS\theta'$ is not necessarily saturated.
As in \ref{ya.7} we define ${\rm Im}'\theta_{\sF,i}$ to be the saturated hull of the image.
\begin{claim}\label{st.7}
Let $\sF$ and $\sF'$ be stable subsheaves of $\mathfrak{gr}_\sS  E^{1,0}$ of slope $\mu_0$. Then
for $i\neq j$ one has ${\rm Im}'\theta_{\sF,i} \cap {\rm Im}'\theta_{\sF',j}=0$.
\end{claim}
\begin{proof}
Otherwise Claim \ref{st.6} implies that ${\rm Im}'\theta_{\sF,i}={\rm Im}'\theta_{\sF',j}$. Then
$$
\sF\otimes T_i \oplus \sF'\otimes T_j
$$ 
contains a direct factor $\sK$ whose image under $\mathfrak{gr}_\sS\theta'$ is zero. By \ref{st.6} both,
$\sF\otimes T_i$ and $\sF'\otimes T_j$ are stable, hence the saturated image $\sK_\iota$ of $\sK$ under 
the projections is $\sF\otimes T_i$ for $\iota=i$ or $\sF\otimes T_j$ for $\iota=j$. 

If $\mu=\mu_\sN$ for an ample invertible sheaf, $\equiv$ will stand for the equality in ${\rm NS}(Y)$. If $\mu=\mu_{\omega_Y(S)}$, then $\equiv$ stands for the equality in ${\rm NS}(Y)/{\rm NS}_0$, for the subspace ${\rm NS}_0$, introduced in Lemma \ref{ya.6}, v).

Since $\Upsilon(\sF)=\Upsilon(\sF')$ (using Lemma \ref{ya.9} if $\mu=\mu_{\omega_Y(S)}$) one finds that 
$$
\ch_1(\sK)\equiv \ch_1(\sK_\iota)\mbox{ \ \ and \ \ }
\Upsilon(\sK)\equiv \Upsilon(\sK_\iota)=\Upsilon(\sF)-\Upsilon(\Omega_\iota).
$$
So $\ch_1(\Omega_i)$ is a rational multiple of $\ch_1(\Omega_j)$, up to classes in ${\rm NS}_0$, contradicting Lemma \ref{sl.5}, v).
\end{proof}
\begin{claim}\label{st.8}
For $i\neq j$ the saturated images $\sB_i$ of $E^{1,0}\otimes T_i$ and $\sB_j$ of $E^{1,0}\otimes T_j$ are disjoint
in $E^{0,1}$.
\end{claim}
\begin{proof}
By Lemma \ref{sl.7} the Higgs field $\theta$ respects the slope filtration. Hence if for some $i\neq j$
Claim \ref{st.8} is wrong, the intersection of the saturated images of $\mathfrak{gr}_\sS (E^{1,0})\otimes T_i$ and of $\mathfrak{gr}_\sS (E^{1,0})\otimes T_j$ in $\mathfrak{gr}_\sS (E^{0,1})$ contains a stable subsheaf
$\sC$ of slope $\mu_0$. So one finds two stable subsheaves $\sF$ and $\sF'$ of $\mathfrak{gr}_\sS (E^{1,0})$
which violate the Claim \ref{st.7}.
\end{proof}
\begin{claim}\label{st.9}
There exists some $\iota$ such that the Higgs field $\theta$ factors through
$E^{1,0} \to E^{0,1}\otimes \Omega_\iota$.
\end{claim}
\begin{proof}
By Claim \ref{st.8} the higgs field $\theta':E^{1,0}\otimes T \to E^{0,1}$ decomposes as a direct sum 
of morphisms
$$ 
\theta':\bigoplus_{i=1}^s (E^{1,0}\otimes T_i) \>>> \bigoplus_{i=1}^s \sB_i\>\subset >> E^{0,1}.
$$
Since $\theta':E^{1,0}\otimes T \to E^{0,1}$ is a morphism between semi-stable sheaves of 
the same slope, the cokernel $\sC$ of ${\rm Im}'(\theta')\hookrightarrow E^{0,1}$ has to be zero.
Otherwise it would be a semi-stable sheaf of slope $\mu(E^{0,1})<0$, contradicting Proposition \ref{sl.4}.

As at the end of the proof of Claim \ref{st.7} we choose for ${\rm NS}_0$ the subspace of ${\rm NS}(Y)$ introduced in Lemma \ref{sl.5}, v), if $\mu=\mu_{\omega_Y(S)}$, and ${\rm NS}_0=0$ if $\mu=\mu_\sN$ for $\sN$ ample. In both cases $\equiv$ stand for the equality of Chern classes in ${\rm NS}(Y)/{\rm NS}_0$. 

Since $\sB_1\oplus\cdots\oplus\sB_s$ is a subsheaf of $E^{0,1}$ with a torsion cokernel and since both have the same slope their first Chern classes are equal in ${\rm NS}(Y)/{\rm NS}_0$.

The morphism $\theta$ factors through
$$
\theta:E^{1,0} \>\subset >> \bigoplus_{i=1}^s \sB_i\otimes \Omega_i \>
\subset >> E^{0,1} \otimes \Omega^1_Y(\log S).
$$
The sheaf $\mathfrak{gr}_\sS  E^{0,1}$ contains the direct sum of the sheaves
$\mathfrak{gr}_\sS \sB_i$, and again both have the same rank and slope, hence the same first Chern class in ${\rm NS}(Y)/{\rm NS}_0$.

We choose an index set $I$, consisting of pairs 
$(\sF,i)$ with $\sF$ a stable subsheaf of $\mathfrak{gr}_\sS E^{1,0}$ 
of slope $\mu_0$ and with $i\in\{1,\ldots,s\}$ such that ${\rm Im}'\theta_{\sF,i}=\mathfrak{gr}_\sS \theta'(\sF\otimes T_i)\neq 0$.

At present we do not know whether ${\rm Im}'\theta_{\sF,i}={\rm Im}'\theta_{\sF',i}$ implies that $\sF=\sF'$. So if this is not the case, we include one $(\sF,i)$ in $I$, but not the other. So writing
$I_j$ for the set of tuples of the form $(\sF,j)\in I$ one has an inclusion 
$$
\bigoplus_{I_j} \sF\otimes T_j \>>> \mathfrak{gr}_\sS \sB_j.
$$
Again both sides have the same rank and slope, hence 
\begin{multline}\label{eqst.2}
 \ \ \ \ \ \ \ \ch_1(E^{1,0})=-\ch_1(E^{0,1})\equiv\\ -\sum_{i=1}^s \ch_1(\sB_i)\equiv
-\sum_I \big( \rk(T_i)\cdot\ch_1(\sF) - \rk(\sF)\cdot \ch_1(\Omega_i)\big).
\end{multline}
Recall that $\rk(\sF)\cdot \ch_1(E^{1,0})= \rk(E^{1,0})\cdot \ch_1(\sF)$, hence (\ref{eqst.2})
implies that
$$
\big(1+\sum_I  \frac{\rk(T_i)\cdot\rk(\sF)}{\rk(E^{1,0})}\big)\cdot \ch_1(E^{1,0}) \equiv
\sum_I \rk(\sF)\cdot \ch_1(\Omega_i).
$$
Let us assume that for some stable saturated subsheaf $\sF'$ of $E^{1,0}$ the sheaf 
$\sB'=\theta'(\sF'\otimes T_1)$ is non zero. By definition, the index set $I$ contains a tuple $(\sF'',1)$ with ${\rm Im}'\theta'_{\sF',1} = {\rm Im}'\theta'_{\sF'',1} $, hence with
$\ch_1(\sF')\equiv\ch_1(\sF'')$. One obtains
\begin{multline}\label{eqst.3}
\ch_1(\sB')\equiv - \rk(\sF')\cdot \ch_1(\Omega_1)+\rk(T_1)\cdot\ch_1(\sF') \equiv \\
- \rk(\sF')\cdot \ch_1(\Omega_1)+\frac{\rk(T_1)\cdot\rk(\sF')}{\rk(E^{1,0})}
\cdot\ch_1(E^{1,0}) \equiv \\
- \rk(\sF')\cdot \ch_1(\Omega_1)+\beta\cdot \sum_I \rk(\sF)\cdot \ch_1(\Omega_i),
\end{multline}
where $\beta$ is a positive rational number. 

If Claim \ref{st.9} is wrong, for some $i>1$ there exists a direct factor
$\sF$ of $\mathfrak{gr}_\sS E^{1,0}$ with $\mathfrak{gr}_\sS \theta'(\sF\otimes T_i)\neq 0$,
say for $i=2$. Then $\rk(\sF)\cdot \ch_1(\Omega_2)$ occurs on the right hand side of 
(\ref{eqst.3}). So we may write
$$
\ch_1(\sB') \equiv \sum_{i=1}^s \beta_i \cdot \ch_1(\Omega_i),
$$
with $\beta_2 >0$ and with $\beta_i \geq 0$ for $i=3, \ldots, s$. Recall that 
\begin{gather*}
\omega_Y(S)=\bigotimes_{i=1}^s \det(\Omega_i) \mbox{ \ \ and that}\\ 
\ch_1(\Omega_1)^{\gamma_1}. \cdots . \ch_1(\Omega_s)^{\gamma_s}=
\left\{\begin{array}{ll} \alpha\cdot \ch_1(\omega_Y(S))^n > 0 & \mbox{if for all } i \ \gamma_i = n_i\\
0 & \mbox{otherwise. }\\
\end{array}\right. 
\end{gather*}
for some $\alpha>0$. Hence (using the definition of ${\rm NS}_0$ if $\mu=\mu_{\omega_Y(S)}$)
the degree of the intersection
\begin{multline*}
\ch_1(\sB').\ch_1(\omega_Y(S) \otimes \det(\Omega_1)^{\nu})^{n-1}=\\
\sum_{i=1}^s \beta_i \cdot \ch_1(\Omega_i).\ch_1(\omega_Y(S) \otimes \det(\Omega_1)^{\nu})^{n-1}=\\
\sum_{i=1}^s \beta_i \cdot \ch_1(\Omega_i).(\ch_1(\omega_Y(S)) + \nu \ch_1(\Omega_1))^{n-1},
\end{multline*}
as a polynomial in $\nu$ is $n_1$ with highest coefficient 
$$
\sum_{i=2}^s \beta_i \cdot \ch_1(\omega_Y(S))^n >0.
$$
On the other hand, the sheaf $\omega_Y(S) \otimes \det(\Omega_1)^{\nu}$ is nef for all $\nu >0$, 
and since $\sB'$ is a subbundle of $E^{0,1}$, the intersection number
$$
\ch_1(\sB').\ch_1(\omega_Y(S) \otimes \det(\Omega_1)^{\nu})^{n-1}
$$
can not be positive. 
\end{proof}
To end the proof of \ref{st.4}, i), let us assume that in \ref{st.9} one has $\iota=1$. 

Let $\sR$ be a stable unitary subsheaf of $T_1\otimes \Omega_1$. By
Addendum \ref{ya.5} $\sR\otimes T_1$ is stable, hence the composite
$$
\sR\otimes T_1 \>>> T_1\otimes \Omega_1\otimes T_1 \> {\rm tr}\otimes {\rm Id} >>
T_1
$$
must be an isomorphism and $\sR$ has to be invertible with $\sR^{n_1}=\sO_Y$.
Since $T_1$ is stable and since by Theorem \ref{ya.4} it remains stable on all
\'etale coverings of $Y$, the subsheaf $\sR$ has to be the subsheaf $\sO_Y\subset T_1\otimes \Omega_1$, given by the homotheties. 

Consider two stable subsheaves $\sF$ and $\sF'$ 
of slope $\mu_0$ in $\mathfrak{gr}_\sS E^{1,0}$. The assumption made in \ref{st.4}, i), and Lemma \ref{st.3}, d), imply that $\delta(\sF)=\delta(\sF')=0$. Since $\delta(\sF)=\delta(\sF^\vee)$ one obtains by Lemma \ref{st.3}, a), that $\delta(\sF^\vee \otimes \sF')=0$. 

If ${\rm Im}'\theta_{\sF,1}={\rm Im}'\theta_{\sF',1}$, then $\rk(\sF)=\rk(\sF')$ and
$\mu(\sF^\vee \otimes \sF')=0$. Hence the conditions ($*$) and ($**$) in the Set-up \ref{in.3} or Lemma \ref{ya.3}, c) and d), imply that $\sF^\vee \otimes \sF'$ is unitary. Then the image of
$$
\sF^\vee \otimes \sF' \>>> T_1\otimes \Omega_1
$$
has to be the factor $\sO_Y$, hence $\sF\cong \sF'$. Altogether, if ${\rm Im}'\theta_{\sF,1}={\rm Im}'\theta_{\sF',1}$ and if $\sF\neq \sF'$, one can change the decomposition of $\mathfrak{gr}_\sS E^{1,0}$ in such a way, that $\sF$ lies in the kernel of 
$$
\mathfrak{gr}_\sS \theta: \mathfrak{gr}_\sS E^{1,0} \>>> \mathfrak{gr}_\sS E^{0,1}\otimes \Omega_1
$$
contradicting the injectivity of $\mathfrak{gr}_\sS \theta$.
 
Therefore the injectivity of $\theta_{\sF,1}$, for all $\sF$, implies the injectivity of
$$
\mathfrak{gr}_\sS \theta'_1:\mathfrak{gr}_\sS E^{1,0}\otimes T_1 \>>> \mathfrak{gr}_\sS E^{0,1}.
$$
hence of $\theta'_1$. One obtains $\rk(E^{0,1})=\ell'=\ell\cdot n_1$.

Since $\ell\cdot \mu(E^{1,0})=-\ell'\cdot\mu(E^{0,1})$ 
the Arakelov equality says that
$$
\mu(\Omega_i)=\mu(E^{1,0})-\mu(E^{0,1})= \\
\frac{\ell+\ell'}{\ell'}\mu(E^{1,0})=
\frac{n_1+1}{n_1}\mu(E^{1,0}).
$$
Let $\sF$ be a stable subsheaf of $E^{1,0}$ of slope $\mu(E^{1,0})$. 
One obtains for $\sF\otimes T_1$
$$
\mu(\sF\otimes T_1)=\mu(\sF)-\mu(\Omega)=\frac{-1}{n_1}
\mu(E^{1,0})=\frac{-1}{n_1}\mu(\sF).
$$
The sheaf $\sF\otimes T_1$ is stable, and we choose $\sB$ to be its saturated hull in $E^{0,1}$.
Then $\mu(\sB)=\mu(\sF\otimes T_1)$ and $\sF\oplus \sB$ defines a saturated sub-Higgs sheaf of $E$ with
\begin{multline*}
\big(\ch_1(\sF) + \ch_1(\sB)\big).\ch_1(\omega_Y(S))^{n_1-1}=
\rk(\sF)\cdot\mu(\sF)+n_1\cdot\rk(\sF)\cdot\mu(\sB)=\\
\rk(\sF)\cdot\mu(\sF)-\rk(\sF)\cdot\mu(\sF)=0.
\end{multline*}
By Proposition \ref{sl.4} it gives rise to a local sub-system of $\V$. Since we assumed the latter 
to be irreducible, $\sF=E^{1,0}$ and $E^{0,1}=\sB$. By Claim \ref{st.6} $\sB$ contains $E^{1,0}\otimes T_1$ as a subsheaf with torsion cokernel.
\end{proof}
\begin{proof}[Proof of Proposition \ref{st.4}, ii)]
Note that Claim \ref{st.6} is obviously true in this case, but we have to replace \ref{st.5} by a different argument. Let us assume that
$\ell=\rk(E^{1,0}) \leq \ell'=\rk(E^{0,1})$.

By assumption, $\Omega$ is the direct sum of invertible sheaves
$\Omega_1\oplus\cdots\oplus\Omega_s$. Let us write 
$$
\ch_1(\omega_Y(S))=L_1 + \cdots + L_s +R
$$
where $L_i=\ch_1(\Omega_i)$ and where $R$ is the first Chern class of the direct 
factor of $\Omega^1_Y(\log S)$ complementary to $\Omega$. Consider 
$$
\theta'_i: E^{1,0}\otimes T_i \>>> E^{0,1}
$$
with saturated image $\sB_i$. The kernel of $\theta'_i$ is a semi-stable subsheaf, and since $T_i$
is invertible, we can write it as $\sK_i\otimes T_i$. So we find a quotient
$\sF_i$ of $E^{1,0}$ such that $\theta'_i$ factors through 
$$
E^{1,0}\otimes T_i\>>> \sF_i\otimes T_i \>\subset >> \sB_i.
$$
Let us assume that there is a stable subsheaf $\sF$ of $E^{1,0}$ of slope $\mu(E^{1,0})$ which is not 
contained in $\sK_1$. Then the saturated hull $\sB$ of $\sF\otimes T_1$ is a stable subsheaf of $\sB_1$ of slope $\mu(E^{0,1})$ and
\begin{multline*}
\ch_1(\sF).L_1.\ch_1(\omega_Y(S))^{n-2}=\\
(\ch_1(\sF)-L_1).L_1.\ch_1(\omega_Y(S))^{n-2}=\ch_1(\sB).L_1.\ch_1(\omega_Y(S))^{n-2}.
\end{multline*}
Since the right hand side is non-positive, $\ch_1(\sF).L_1.\ch_1(\omega_Y(S))^{n-2}\leq 0$.
For some $\alpha,\beta >0$ one can write
$$
\ch_1(\omega_Y(S))^{n-1}= \alpha\cdot L_1.\ch_1(\omega_Y(S))^{n-2} + \beta\cdot C_1
$$
where $C_1=L_2 . \cdots . L_s . R^{n-s}$.
If $\sF\not\subset \sK_2$ the same argument shows that for the saturated hull $\sB'$ of $\sF\otimes T_2$
$$
\ch_1(\sF).C_1= 
(\ch_1(\sF)-L_2).C_1=\ch_1(\sB').C_1 \leq 0,
$$
hence $\ch_1(\sF).\ch_1(\omega_Y(S))^{n-1}\leq 0$ and $\mu(\sF)\leq 0$. On the other hand,
the $\mu$-semi-stability implies that $\mu(\sF)=\mu(E^{1,0})$ and the Arakelov equality implies
that $\mu(E^{1,0})$ is a positive multiple of $\ch_1(\omega_Y(S))^n$, a contradiction. 

Therefore $\sF$ has to lie in the kernels $\sK_i$ for $i\neq 1$ and
one has a factorization of the restriction of $\theta'$ to $\sF$ like
$$
\sF \otimes T \>>> \sF \otimes T_1 \>>> \sB \subset E^{0,1}.
$$
Then $\sF \oplus \sB$ is a sub-Higgs sheaf. Since both sheaves have the same rank
$$
\mu(E^{1,0})+ \mu(E^{0,1})=\mu(\sF)+
 \mu(\sB)=\frac{1}{\rk(\sF)}(\ch_1(\sF)+\ch_1(\sB)) \leq 0.
$$
The equality
$$
0=\ell\cdot \mu(E^{1,0})+ \ell' \cdot \mu(E^{0,1})
= \ell'(\mu(E^{1,0})+ \mu(E^{0,1}))-
(\ell'-\ell)\mu(E^{1,0}),
$$
together with the positivity of $\mu(E^{1,0})$ imply
that $\ell=\ell'$ and that 
$$
\mu(\sF)+  \mu(\sB) = 0.
$$
By Proposition \ref{sl.4} $\sF \oplus \sB$ is induced by a local sub-system of $\V$. 
Since we assumed the latter to be irreducible, one finds
$\sF=E^{1,0}$ and $\sB=E^{0,1}$. So both are stable, and the Higgs field
is of the form asked for in \ref{st.4}. 

It remains to verify that $S^\ell(E^{1,0})\otimes \det(E^{1,0})^{-1}$
and $E^{1,0}\otimes{E^{1,0}}^\vee$ are unitary, or equivalently that 
$$
\delta(S^\ell(E^{1,0}))=\delta(E^{1,0}\otimes{E^{1,0}}^\vee)=0.
$$ 
By \ref{st.3}, a) and c), it is sufficient to show that $\delta(E^{1,0})=0$.
Since $E$ is the Higgs bundle of a local system, $\ch_\iota(E^{1,0}\oplus E^{0,1})=0$ 
for $\iota=1,\ 2$. This implies
$$
\ch_1(E^{1,0})=-\ch_1(E^{0,1})=-\ell\cdot\ch_1(T_1)+\ch_1(E^{1,0}),
$$
hence $\ch_1(E^{1,0})$ and $\ch_1(E^{0,1})$ are both rational multiples of $\ch_1(T_1)$.
By \ref{ya.6} $\ch_1(T_1)^2$ is numerically trivial, and therefore  
$$
\ch_1(E^{1,0})^2=\ch_1(E^{0,1})^2=\ch_1(E^{1,0}).\ch_1(E^{0,1})\equiv 0.
$$
The last equality implies that 
$$
0\equiv \ch_2(E^{1,0}\oplus E^{0,1})= \ch_2(E^{1,0}) + \ch_2(E^{0,1})=
\ch_2(E^{1,0}) + \ch_2(E^{1,0}\otimes T_1).
$$
Since $T_1$ is invertible, and $\ch_1(T_1)^2\equiv 0$ one finds $\ch_2(E^{1,0}\otimes T_1)\equiv \ch_2(E^{1,0})$,
hence $\ch_2(E^{1,0})\equiv 0$ and $\delta(E^{1,0})=0$.
\end{proof}

\section{Families of Abelian varieties}\label{av}

For a smooth family $f:V\to U$ of Abelian varieties consider the variation of Hodge structures
$R^1f_*\C_V$, with Higgs bundle $(F^{1,0}\oplus F^{0,1},\tau)$. If the local monodromies are
uni-potent, the Kodaira-Spencer map induces a morphism
\begin{equation}\label{eqav.1}
T^1_Y(-\log S) \>>> F^{0,1}\otimes {F^{1,0}}^\vee={F^{0,1}}^{\otimes 2}.
\end{equation}
As well known, it factors through $S^2(F^{0,1})$. If $\varphi: U \to \sA_g$ denotes the morphism to the
moduli stack, induced by $f:V\to U$, the sheaf $S^2(F^{0,1})|_U$ is nothing but the
pullback $\varphi^*(T^1_{\sA_g})$.

\begin{lemma}\label{av.1}
Let $f:X \to Y$ be a smooth family of Abelian varieties, and assume that for some $N\geq 3$ and for some \'etale
covering $Y'\to Y$ the morphism $\varphi: Y \to \sA_g$ lifts to an embedding $\varphi':Y'\to \sA_g^{(N)}$,
where $\sA_g^{(N)}$ denotes the moduli scheme of Abelian varieties with a level $N$-structure. Then $\Omega^1_Y$ is nef and $\omega_Y$ is ample.
\end{lemma}
\begin{proof}
One may assume that $Y=Y'$. Recall that the bundle $F^{1,0}$ is nef. Then
the sheaf $\varphi^*\Omega^1_{\sA_g^{(N)}}=S^2(F^{1,0})$
is nef, hence the same holds true for the image of $\rho:\varphi^*\Omega^1_{\sA_g}\to\Omega^1_Y$.
Since we assumed that $\varphi$ is an embedding, the restriction map $\rho$ is surjective, and $\Omega_Y^1$ is nef.

A similar result, for submanifolds of arbitrary period domains can be found in \cite{Zu}.
There, in the proof of Lemma 2.2, one also finds the necessary calculations for the ampleness of $\omega_Y$.
One knows already that the Chern form of $T^1_Y$ is negative semi-definite. Hence in order to
see that $\det(\Omega_Y^1)$ is ample, one just has to show that in each point $y\in Y$ this Chern form is
strictly negative. This is shown in \cite[2.2]{Zu}, provided the differential of the period map is injective
in $y$, a condition which is satisfied for all $y\in Y\subset \sA_g^{(N)}$.
\end{proof}

Let us return to the Set-up \ref{in.3}, hence $U\subset Y$ is the complement of a normal crossing divisor.
Consider the decomposition of the variation of Hodge structures $R^1f_*\C_V$ in irreducible
$\C$-sub-variations of Hodge structures.
If one of the factors, say $\V_i$, is not defined over $\R$, we write $\bar\V_i$ for its conjugate.
Hence numbering the factors such that the $\V_i$ for $r' < i \leq r$ are exactly the ones which
are defined over $\R$ and irreducible over $\C$ we may write
$$
R^1f_*\C_V=\V_1 \oplus \bar\V_1 \oplus \cdots \oplus \V_{r'} \oplus \bar\V_{r'} \oplus \V_{r'+1}
\oplus \V_{r'+2} \oplus \cdots \oplus \V_r.
$$
Let $(E_i^{1,0}\oplus E^{0,1}_i,\theta_i)$ denote the logarithmic Higgs bundle of $\V_i$.
For $i\leq r'$ the sheaf $E_i^{{1,0}^\vee}\subset \varphi^*(E_{\sA_g}^{0,1})$ 
intersects $E_i^{0,1}$ in $0$. Renumbering, we may assume that 
$$
\ell_i=\rk(E_i^{1,0}) \leq \ell'_i=\rk(E^{0,1}_i)\mbox{ \ \ for \ \ } i=1,\ldots,r'
$$
and, of course, $\ell_i=\rk(E_i^{1,0}) = \rk(E^{0,1}_i)$ for $i=r'+1,\ldots,r$.
The sheaf $H=\bigoplus_{i=1}^rH_i$ with
$$
\begin{array}{ll}
H_i=E^{0,1}_i\otimes {E^{1,0}_i}^\vee & \mbox{if \ \  }i\in \{1,\ldots,r'\} \\
H_i=S^2(E_{i}^{0,1}) & \mbox{if \ \  }i\in \{r'+1,\ldots,r\}
\end{array}
$$
is a direct factor of $S^2(F^{0,1})$, and the image of the Kodaira-Spencer map (\ref{eqav.1})
lies in $H$.

Assume that $\omega_Y(S)$ is nef and ample with respect to $U$,
and let $T_1,\ldots,T_s$ be the stable direct factors of $T^1_Y(-\log S)$.
If the map $\varphi:U \to \sA_g$, induced by the smooth family
$f:V\to U$, is generically finite, the Higgs field induces an injection
$T^1_Y(-\log S) \to H$. Hence for each $T_j$ there is some $i$ such that the composite
$$
T_j \> \subset >> T^1_Y(-\log S) \>>> H=\bigoplus_{i=1}^r H_i \>{\rm pr}_i >> H_i
$$
is non-zero. On the other hand, if $\V_i$ is not unitary, there exists some $j$ such that
$T_j\to H_i$ is non-zero. So we can restate Proposition \ref{in.4} in the following form:
\begin{proposition}\label{av.2} Assume that $\varphi:U\to \sA_g$ is generically finite
and that $Y$, $U$ and $f$ satisfy the assumptions made in Set-up \ref{in.3}. Assume moreover that
$S^m(T_j)$ is stable for all $m>0$ and all $j$. If $\V_i$ is not unitary there exist a unique $j$ with
$T_j\to H_i$ non-zero. Moreover $T_j$ is a direct factor of $H_i$.
\end{proposition}
\begin{proof}
We can apply Proposition \ref{st.4} to $\V_i$, hence we know that the sheaves $H_i$ are poly-stable and that there exists exactly one $j$ with $T_j\to H_i$ non-zero. By the Arakelov equality for $\V_i$ the slope of the sheaves $T_j$ and $H_i$ coincide, hence $T_j\to H_i$ has to be injective with a splitting image.
\end{proof}
\begin{remark}\label{av.3}
The dual of the Higgs field of $\V_{\sA_g}$ is given by a tautological map
$$
\theta_{\sA_g}^\vee: E_{\sA_g}^{1,0}\otimes S^2(E_{\sA_g}^{0,1}) \>>> E_{\sA_g}^{0,1}.
$$
For its description recall that $E_{\sA_g}^{0,1}$ is dual to $E_{\sA_g}^{1,0}$, so locally we may
choose a basis $e_1,\ldots , e_g$ of $E_{\sA_g}^{1,0}$ and for $E_{\sA_g}^{0,1}$ a dual basis
$e_1^\vee,\ldots,e_g^\vee$. Then
\begin{equation}\label{eqav.2}
\theta^\vee_{\sA_g}(e_i\otimes e^\vee_j\cdot e^\vee_k)= \frac{1}{2}(e^\vee_j(e_i)e^\vee_k+e^
\vee_k(e_i)e^\vee_j)=\delta_{i,j}e^\vee_k+\delta_{i,k}e^\vee_j.
\end{equation}
Returning to the decomposition of $\varphi^*\V_{\sA_g}$, choose some $i>r'$. Then
$\theta^\vee_{\sA_g}$ induces 
$$
\theta_i^\vee: E_i^{1,0}\otimes S^2(E_i^{0,1}) \>>> E_i^{0,1}
$$
again of the form in (\ref{eqav.2}), with $g$ replaced by $\ell_i=\rk(E_i^{1,0})$.

For $i\leq r'$ remark that a basis of $E_0^{1,0}\otimes (E_0^{{1,0}^\vee}\otimes E_0^{0,1})$ 
is given by
$e_i\otimes e^\vee_j\cdot e^\vee_k=e_i\otimes e^\vee_j\otimes e^\vee_k$ with $1\leq i,j \leq \ell_i$ 
and $\ell_i < k \leq \ell_i+\ell'_i$. The image of such an element is $\delta_{i,j}e^\vee_k$.
\end{remark}
\begin{remark}\label{av.4}
Under the assumption made in Theorem \ref{in.6} assume that $\varphi:U\to \sA_g$ is generically finite. By Lemma \ref{av.2} the tangent sheaf $T_Y^1(-\log S)$ is a direct factor of $H$, hence
$$
T^1_Y(-\log S)\>>> \varphi^*(T^1_{\overline{\sA}_g}(-\log(\overline {\sA}_g \setminus \sA_g)))|_Y
$$
splits and $U\to \sA_g$ is \'etale. Moreover $U$ is totally geodesic in $\sA_g$,
and we can apply B. Moonen's characterization of bounded symmetric domains in $\sA_g$
\cite{Mo}. So we see already that the Arakelov equality (\ref{eqin.3}) forces $U$ to be a bounded symmetric domain. However in order to see that $U$ is a Shimura subvariety
one needs in \cite{Mo} the existence of a CM-point. This condition enters the scene, since it forces by \cite{An} the monodromy group to coincide with the derived Mumford-Tate group
(see Section \ref{sh}).

So we will argue in a different way in the next sections. We first verify the explicite description of the variation of Hodge structures, stated in Theorem \ref{in.5}. It will allow in
Section \ref{rg} to determine the possible Hodge cycles, and in Section \ref{sh}
to describe the Mumford-Tate group. As it will turn out, this description implies that
the monodromy group is equal to the derived Mumford-Tate group, and $U$ must be a Shimura variety.
So finally we will obtain the existence of CM points as a Corollary.
\end{remark}
\begin{remark}\label{av.5}
Assume that a variant of Proposition \ref{in.4} holds true, which does not require the conditions i) and ii) and which allows direct factors $\Omega_i$ with $S^m(\Omega_i)$ non-stable
for some $m>1$. Then for each irreducible $\C$-sub-variation $\V$ of Hodge structures
in $R^1f_*\C_V$ one would obtain exactly one direct factor $\Omega_j$ of $\Omega_Y^1(\log S)$,
as in Lemma \ref{av.2}.

If $\Omega_j$ is invertible or if $S^m(\Omega_j)$ remains stable and if $\delta(\V)=0$
one can apply the methods of the next Sections to describe $\V$. 
It remains to study the case, where for some stable direct factor
$\Omega_j$ of $\Omega^1_Y(\log S)$ there is some $m>1$ with $S^{m}(\Omega_j)$ non-stable. Then the corresponding factor $M_j$ of $\tilde{U}$ is locally Hermitian symmetric domain of $\rk > 1$ and the Superrigidity Theorem of Margulis applies. So in this case, one should be able to understand the corresponding variation of Hodge structure by different methods.
\end{remark}
\section{The structure of $\tilde{U}$ in Theorem \ref{in.5}}\label{pr}
In this section we will show, that the assumptions of Theorem \ref{in.5}
imply that $M_i$ is a complex ball.Although this is obvious
if $\Omega_i$ is invertible, we will not exclude this case in the beginning, and  we will allow that 
the slope and the discriminant is chosen according to the Assumptions \ref{sl.5}. 

So we will use the following Set-up. The sheaf $\Omega=\Omega_i$ is a stable direct factor of $\Omega^1_Y(\log S)$
of rank $n'$, and $T=\Omega^\vee$. Let $\V$ be a sub-variation of Hodge structures in $R^1f_*\C_V$, with Higgs bundle 
$$
(E^{1,0}\oplus E^{0,1},\theta: E^{1,0}\to E^{1,0}\otimes \Omega \to E^{1,0}\otimes \Omega^1_Y(\log S)).
$$ 
Assume that $E^{1,0}$ and $E^{0,1}$ are both stable. Writing  
$\ell=\rk(E^{1,0})$ and $\ell'=\rk(E^{0,1})$ one may assume by the condition c) in \ref{in.4} that $\ell'=n'\cdot \ell$. 

Assume moreover, that $S^m(\Omega)$ remains stable for all $m>0$, hence that in \ref{ya.4} one has
$i\leq s'$. We will use the condition $\delta(\V)=0$, to show that if $\Omega$ is not invertible the equation (\ref{eqya.1}) in Theorem \ref{ya.4}, c), holds true.

Note that 
$$
\ell\cdot \Upsilon(E^{1,0})+\ell\cdot n' \cdot \Upsilon(E^{0,1})=\ch_1(E)=0.
$$
Hence $\displaystyle\mu(E^{1,0})-\mu(E^{0,1})= \frac{n'+1}{n'}\cdot \mu(E^{1,0})$ and
the Arakelov equality says that
$$
\mu(E^{1,0})=\frac{n'}{n'+1}\cdot\mu(\Omega) \mbox{ \ \ and \ \ }
\mu(E^{0,1})=\frac{-1}{n'+1}\cdot\mu(\Omega). 
$$
The local system $\V$ induces local systems $\bigotimes^\ell \V$ and $\bigwedge^\ell \V$. In \cite[p. 40-43]{Sim2} one finds the construction of the corresponding Higgs bundle and Higgs fields.
In particular for some $\theta$ obtained as the direct sum of
$$ 
\theta_m: \bigwedge^{\ell-m}(E^{1,0}) \otimes \bigwedge^m (E^{0,1}) \>>>
\bigwedge^{\ell-m-1}(E^{1,0}) \otimes \bigwedge^{m+1} (E^{0,1}) \otimes \Omega,
$$ 
the Higgs bundle of $\bigwedge^\ell \V$ is
$$
(E^{(\ell)},\theta)=\big(\bigoplus_{m=0}^{\ell}E^{\ell-m,m},\oplus_{m=0}^{\ell-1}\theta_m\big)
$$
with $E^{\ell-m,m}=\bigwedge^{\ell-m}(E^{1,0}) \otimes \bigwedge^m E^{0,1}$.
Define $G=\bigoplus_{m=0}^\ell G^{\ell-m,m}$ as the saturated image
$$
G^{\ell-m,m}={\rm Im}'\big(\det(E^{1,0})\otimes S^m(T) \>>> \bigwedge^{\ell-m}
(E^{1,0}) \otimes \bigwedge^{m} (E^{0,1})\big). 
$$
Since 
$$
\det(E^{1,0})\otimes S^m(T) \>>> \bigwedge^{\ell-m}
(E^{1,0}) \otimes \bigwedge^{m} (E^{0,1})
$$
is a morphism between semi-stable sheaves of the same slope,
and since $S^m(T)$ is stable, the sheaves $G^{\ell-m,m}$ are either $0$ or 
isomorphic to $$\det(E^{1,0})\otimes S^m(T).$$
Obviously  
$$
\theta(G^{\ell-m,m}) \subset \bigwedge^{\ell-m-1}(E^{1,0}) \otimes \bigwedge^{m+1} (E^{0,1}) \otimes 
\Omega,
$$
must be contained in $G^{\ell-m-1,m+1}\otimes \Omega$ and we obtain:
\begin{claim}\label{pr.1}
$G$ is a saturated sub-Higgs sheaf of $E^{(\ell)}$.
\end{claim}

\begin{claim}\label{pr.2}
The sheaf $G^{\ell-m,m}$ is a direct factor of $$\bigwedge^{\ell-m}(E^{1,0}) \otimes \bigwedge^{m} 
(E^{0,1}).$$
\end{claim}
\begin{proof}
$\bigwedge^{\ell-m}(E^{1,0}) \otimes \bigwedge^{m} (E^{0,1})$ is poly-stable, $G^{\ell-m,m}$ is stable,
and both have the same slope.
\end{proof}
\begin{claim}\label{pr.3}
$G$ is a direct factor of the Higgs bundle $E^{(\ell)}$. In particular
$\ch_1(G)$ and $\ch_2(G)$ are both zero.
\end{claim}
\begin{proof}
As we have just seen, $G^{\ell-m,m}$ is a direct factor $E^{(\ell)}$. It remains to show, that the projections 
$$
\Phi_m:\bigwedge^{\ell-m}(E^{1,0})\otimes \bigwedge^m(E^{0,1}) \>>> G^{\ell-m,m},
$$
to $G^{\ell-m,m}$ can be chosen such that 
$G= G^{\ell,0} \oplus G^{\ell-1,1} \oplus \cdots G^{\ell-r,r}$ 
is a quotient Higgs bundle. We will construct the splittings by descending induction in such a way
that the diagram
$$
\begin{CD}
\displaystyle \bigwedge^{\ell-m}(E^{1,0})\otimes \bigwedge^m(E^{0,1}) \>\theta^{\ell-m,m} >> 
\displaystyle \bigwedge^{\ell-m-1}(E^{1,0})\otimes 
\bigwedge^{m+1}(E^{0,1})\otimes \Omega\\
\V \Phi_m VV \V V \Phi_{m+1} V\\
G^{\ell-m,m} \>>> G^{\ell-m-1,m+1} \otimes \Omega
\end{CD}
$$
commutes. As long as $G^{\ell-m-1,m+1}=0$ there is nothing to construct, and we can choose
$\Phi_m$ to be any splitting, existing by \ref{pr.2}. 

If $r$ is the largest integer with $G^{\ell-r,r}\neq 0$, assume by induction, that we found 
the $\Phi_{m'}$ for all $m'>m$ and that $m<r$.

So $\theta^{\ell-m,m}|_{G^{\ell-m,m}}$ is non-zero. Since $G^{\ell-m,m}$ is stable and since
$\theta^{\ell-m,m}$ a morphism between poly-stable sheaves of the same slope, one finds
$$
G^{\ell-m,m}\> \theta^{\ell-m,m} > \subset >
G^{\ell-m-1,m+1}\otimes \Omega \subset \bigwedge^{\ell-m-1}(E^{1,0})\otimes 
\bigwedge^{m+1}(E^{0,1})\otimes \Omega.
$$
So the saturated image of $G^{\ell-m,m}$ under $\Phi_{m+1}\circ \theta^{\ell-m,m} $ is isomorphic to
$G^{\ell-m,m}$, and 
$$
G^{\ell-m,m} \> \subset >> \bigwedge^{\ell-m}(E^{1,0})\otimes \bigwedge^m(E^{0,1}) 
\>\Phi_{m+1} \circ\theta^{\ell-m,m}>>  G^{\ell-m,m}.
$$
defines a splitting $\Phi_m$ with the desired properties. 

So $G$ splits as a sub-Higgs bundle of $E^{(\ell)}$, hence it is itself a Higgs bundle arising from
a local system. Then all the Chern classes $\ch_i(G)$ are zero. 
\end{proof}
\begin{proposition}\label{pr.4} One has $G^{0,\ell}\neq 0$.
\end{proposition}
\begin{proof}
Let us write again $r$ for  the largest integer with $G^{\ell-r,r}\neq 0$. For $0\leq m\leq r$
the sheaf $G^{\ell-m,m}$ is a stable sheaf of slope 
\begin{multline*}
(\ell-m)\cdot \mu(E^{1,0}) + m\cdot \mu(E^{0,1})=
\ell \cdot \mu(E^{1,0}) - m\cdot\mu(\Omega_Y^1(\log S))=\\
\big(\frac{\ell \cdot\ell'}{\ell+\ell'}-m\big)\cdot \mu(\Omega_Y^1(\log S)),
\end{multline*}
and of rank $g^{\ell-m,m}=\binom{n'+m-1}{m}$. 
By \ref{pr.3} the first Chern class of $G$ is zero, hence
\begin{multline*}
0=\frac{\mu(\det(G))}{\mu(\Omega^1_Y(\log S))}=\sum_{m=0}^r 
\binom{n'+m-1}{m}\cdot \big(\frac{\ell \cdot\ell'}{\ell+\ell'}-m\big)= \\
\left(\frac{\ell \cdot\ell'}{n'\cdot(\ell+\ell')} 
- \frac{r}{n'+1}\right)\cdot (r+1)\cdot\binom{r+n'}{r+1}.
\end{multline*}
Since $\ell'=n'\cdot \ell$ one finds that
$$
0=\frac{\ell \cdot\ell'}{n'\cdot(\ell+\ell')} 
- \frac{r}{n'+1}= \frac{\ell}{n'+1} 
- \frac{r}{n'+1},
$$
and $r=\ell$.
\end{proof}
Finally we will show that the factor $M$ in the universal covering $\tilde{U}$ corresponding to $\Omega$ is a complex ball. This is obvious if $n'=1$. For $n'>1$ we just have to verify the condition (\ref{eqya.1}) in \ref{ya.4}, c). This is done in the next Proposition.
\begin{claim}\label{pr.5} Assume that $\rk(\Omega)>1$. Then
the condition $\ch_2(G)=0$ implies that $2 \cdot (n'+1) \cdot \ch_2(\Omega)
- n' \cdot \ch_1(\Omega)^2=0$.
\end{claim}
\begin{proof}
The claim follows by a formal calculation of Chern numbers. 
Hence we may replace $Y$ by any finite covering, and assume 
that there exists an invertible sheaf $\sL$ with
$\det(E^{1,0})=\sL^\ell$.  Or we may calculate with $\Q$-Chern classes. 
Consider the sheaf
$$
F=F^{1,0}\oplus F^{0,1} \mbox{ \ \ with \ \ }
F^{1,0}=\sL, \ \ F^{1,0}=\sL \otimes T.
$$
Then $S^\ell(F)$ is a Higgs bundle with $\sL^\ell\otimes S^m(T)$ in bidegree $(\ell-m,m)$, hence
isomorphic to $G$. By \ref{pr.3} the first Chern class of $G$ is zero, hence $\ch_1(F)$ as well.
On the other hand,
$$
\ch_1(F)= \ch_1(\sL) + n'\cdot \ch_1(\sL) - \ch_1(\Omega)=
\frac{n'+1}{\ell}\ch_1(E^{1,0}) - \ch_1(\Omega),
$$
and $\ch_1(\sL)=\frac{1}{n'+1}\ch_1(\Omega)$. For the second Chern class
it is easier to calculate the discriminant
$$
\Delta(\sF)=2 \cdot \rk(\sF) \cdot \ch_2(\sF) - (\rk(\sF)-1)\cdot \ch_1(\sF)^2.
$$
By \ref{st.3}, a), the discriminant is invariant under tensor products with invertible sheaves, 
hence $\Delta(\sL\oplus \sL\otimes T) = \Delta(\sO_Y \oplus T)$.

Claim \ref{pr.3} implies that $\ch_1(G)^2=\ch_2(G)=0$, hence $\Delta(G)=0$, and 
from \ref{st.3}, b), one obtains $\Delta(F)=0$. Then
$$
0=\Delta(\sO_Y \oplus T)=2 \cdot (n'+1) \cdot \ch_2(T)
- n' \cdot \ch_1(T)^2,
$$
as claimed.
\end{proof}

\section{Higgs bundles in the non-compact case}\label{nc}
Consider a family of Abelian varieties 
$f:V\to U$ and an irreducible $\C$-sub-variation of Hodge structures $\V$
of $R^1f_*\Q_V$. We will assume in that the assumptions made in Theorem \ref{in.5} hold true
for $\V$, in particular the conditions a)--c) in Proposition \ref{in.4}.
Lemma \cite[3.2]{VZ1} allows to find a finite Galois extension $K$ of $\Q$ with Galois group $\Gamma$,
and a $K$-sub-variation of Hodge structures $\V_K$ in $R^1f_*K_V$ such that
$\V=\V_K\otimes_K\C$. 

We will assume that $\ell=\rk(E^{1,0})\leq \rk(E^{0,1})$, so by assumption there exists a unique $\iota=\iota(\V) \in \{1,\ldots,s\}$ with $E^{0,1}=E^{1,0}\otimes T_\iota$, and the Higgs field $\theta$ is given by the natural embedding
$$
E^{1,0} \>>> E^{1,0}\otimes \sE nd (\Omega_\iota) =E^{1,0}\otimes T_\iota\otimes \Omega_\iota.
$$
\begin{assumption}\label{nc.1}
Let $S=S_1\cup\ldots\cup S_\eta$ be the decomposition of $S$ in irreducible components and let
$\gamma_j\in \pi_1(U,*)$ be the image of a generator of the local fundamental group
of a small neighborhood of a general point of $S_j$. Assume that 
the image of $\langle \gamma_1,\ldots,\gamma_\eta\rangle$ 
under the representation corresponding to $\V$ is non-trivial.   
\end{assumption}
Let us remark, that the local monodromies of $\V$ are uni-potent.
Hence if the image of $\gamma_j$ under the representation is non-trivial,
it has to have infinite order.

Note that the assumption \ref{nc.1} and the description of the Higgs field of $\V$ given above remain true if one replaces $Y$ by an \'etale covering and $\V$ by an irreducible direct factor of its pullback. So by abuse of notations we will assume that the pullback of $\V$ remains irreducible on all \'etale coverings.

Let us consider the $\C$-variation of Hodge structures $\END(\V)$
with Higgs bundle
$$
(\End=\bigoplus_{p=-1}^{1}\End^{-p,p} , \ \rho=\bigoplus_{p=-1}^{1}\rho_{-p,p}).
$$ 
Then $\End^{1,-1}=E^{1,0}\otimes {E^{1,0}}^\vee\otimes \Omega_\iota$,
$$
\End^{0,0}=E^{1,0}\otimes {E^{1,0}}^\vee \oplus
E^{1,0}\otimes {E^{1,0}}\otimes T_\iota \otimes \Omega_\iota,
$$ 
and $\End^{-1,1}=E^{1,0}\otimes {E^{1,0}}^\vee\otimes T_\iota$.
Again the Higgs field is induced by the tautological maps
on $T_\iota$ and $\Omega_\iota$. In particular, the kernel $\sK$ of
$$
\rho_{0,0}:\End^{0,0}\>>> \End^{-1,1}\otimes \Omega_\iota
$$
is isomorphic to $E^{1,0}\otimes {E^{1,0}}^\vee$, diagonally embedded into
$\End^{0,0}$, and we obtain the first part of:
\begin{lemma}\label{nc.2}
$\END(\V)$ contains a unitary local sub-system $\U$ of rank $\ell^2$. 
It is the largest unitary local sub-system. Moreover $\U$ is defined over some number field,
as well as the decomposition $\END(\V)=\U\oplus \M$. 
\end{lemma}
\begin{proof}
The explicite description of the Higgs field given above, shows that
$\END(\V)=\U \oplus \M$ where the Higgs field of the Higgs bundle
for $\M$ is given by
\begin{equation}\label{eqnc.1}
E^{1,0}\otimes {E^{1,0}}^\vee\otimes \Omega_\iota \>>>\\
(E^{1,0}\otimes {E^{1,0}}^\vee \oplus
E^{1,0}\otimes {E^{1,0}}\otimes T_\iota \otimes \Omega_\iota)/\sK \otimes \Omega_\iota
\end{equation}
and by 
\begin{equation}\label{eqnc.2}
(E^{1,0}\otimes {E^{1,0}}^\vee \oplus
E^{1,0}\otimes {E^{1,0}}\otimes T_\iota \otimes \Omega_\iota)/\sK
\>>> E^{1,0}\otimes {E^{1,0}}^\vee\otimes T_\iota\otimes \Omega_\iota.
\end{equation}
Obviously both morphisms are injective, hence $\U$ is maximal.
 
To find the field of definition, we argue as in the proof of \cite[3.3]{VZ1}. 
Consider a family
$\{\M_t\}_{t\in \Delta}$ of local sub-systems of $\END(\V)$ defined over a small disk
$\Delta$, with $\M_0=\M$. For $t\in \Delta$ let $(F_t,\theta_t)$ denote the Higgs bundle
of an irreducible direct factor of $\M_t$. Then 
$$
F_t^{1,-1} \>>> F_t^{0,0}\otimes \Omega_\iota \mbox{ \ \ and \ \ }
F_t^{0,0} \>>> F_t^{-1,1}\otimes \Omega_\iota 
$$
are both injective for $t$ sufficiently small. If the composite
$$
\rho:\M_t \> \subset >> \END(\V) \>>> \U 
$$ 
is non-zero, the complete reducibility of local systems coming from
variations of Hodge structures implies that $\M_t$ and $\U$
contain a common direct factor. Since the Higgs field of $\U$ is trivial,
one obtains a contradiction.

So $\M$ is rigid as a local sub-system, hence it can be defined over some 
number field $K$. As in the proof of \cite[3.3]{VZ1} this implies that
$\U$ is also defined over a number field.
\end{proof}
\begin{lemma}\label{nc.3}
Let us keep the notations and assumptions from \ref{nc.1}. Let $\BS$ be an irreducible
direct factor of $\END(\V)$. Then the image of $\langle \gamma_1,\ldots, \gamma_\eta\rangle$ under the representation corresponding to $\BS$ is trivial, if and only if $\BS\subset \U$.
\end{lemma}
\begin{proof}
If $\BS\subset \U$ then ``unitary and uni-potent'' implies that the image of each $\gamma_j$ is trivial.

On the other hand, the assumption \ref{nc.1} implies that the local system $\V$ can not be extended as a local system to $Y$. This implies in particular, that $\Omega_\iota\subset \Omega^1_Y(\log S)$ does not factor through an inclusion
$\Omega_\iota\subset \Omega^1_Y$. 

The components of the Higgs field of $\M$ in (\ref{eqnc.1}) and (\ref{eqnc.2}) are injective morphisms between semi-stable sheaves of the same slope, and the second one is an isomorphism. Those two properties carry over to
all direct factors $F$ of the Higgs bundle. In fact, if $F$ is non-zero, the injectivity implies that
$F^{-p,p}\neq 0$ for $p=-1,0,1$. The surjectivity of $F^{0,0}\to F^{-1,1}\otimes \Omega_\iota$
implies that the Higgs field has non-trivial poles along at least one component of $S$. 
Otherwise one would obtain an inclusion 
$$
F^{-1,1}\otimes \Omega_\iota\> \subset >> F^{-1,1}\otimes \Omega^1_Y,
$$
contradicting the Assumption \ref{nc.1}. 
\end{proof}
Recall that $\V=\V_K\otimes_K\C$ for $K$ a number field. Let $\W_\Q$ be the Weil restriction of
$\V$, and $\W=\W\otimes_\Q \C$. Writing $\V=\V_1,\V_2,\ldots,\V_r$ for the different conjugates of $\V$
under the Galois group $\Gamma$ the local system $\END(\W)$ contains 
$$
\E=\END(\V)\oplus \END(\V_2) \oplus \cdots \oplus \END(\V_r)
$$
as a local sub-system, obviously invariant under conjugation by $\Gamma$, and hence defined over $\Q$.
Since the characterization of the maximal unitary direct factor in Lemma \ref{nc.3} is invariant under conjugation,
one finds that each $\V_i$ contains a maximal unitary direct factor $\U_i$ of rank $\ell^2$, and that those
are conjugate under conjugation by $\Gamma$. The local system
$\T=\U\oplus \U_2 \oplus \cdots \oplus \U_r$ is invariant under $\Gamma$, hence defined over $\Q$.
\begin{corollary}\label{nc.4}
Over some \'etale covering $\phi:Y'\to Y$, one has
$$
\phi^* \T \cong \bigoplus \sO_{Y'}.
$$
\end{corollary}
\begin{proof}
The proof of \cite[4.1]{VZ1} does not use the fact that the base is a curve.
One just needs, that for each irreducible local sub-system of the non-unitary
part of a local system, there is some $S_j\subset S$ such that the image of $\gamma_j$
in the corresponding representation of $\pi_1(U,*)$ is non-trivial.
This holds true for the complement of $\T$ in $\E$.
One finds that $\T$ is defined over $\Q$, and that it extends to a unitary local system on $Y$. 

Since $\T$ is a local sub-system of a $\Q$-local system, as in \cite[4.3]{VZ1} one can define a $\Z$-structure on $\T$. A unitary local system on $Y$ with a $\Z$ structure will be trivial over some \'etale covering. 
\end{proof}
Remark that the maximal unitary part $\U$  of $\END(\V)$
is a direct factor of $\T$, hence its pullback to $Y'$ is trivial. Then
$\END(\phi^*\V)$ has $\ell^2$ linear independent global sections, one of them given by homotheties.
Since we required $\V$ to remain irreducible under \'etale coverings of $Y$, the local system
$\phi^*\V$ is irreducible. If $\ell >1$, one finds one section with a non-trivial kernel. 
So $\ell=1$ and $E^{1,0}$ is invertible. By \ref{st.4} one has $E^{0,1}=E^{1,0}\otimes T_\iota$.
Then
$$
\det(E^{1,0})^{n_1+1}\otimes \det(\sO_Y\oplus T_\iota)=
\det(E^{1,0})^{n_1+1}\otimes \det(T_\iota)
$$
is equal to $\det(E^{1,0}\oplus E^{0,1})$. This is the determinant of a local system with uni-potent local monodromy
around the components of $S$, hence an element of ${\rm Pic}^0(Y)$. The later is divisible, and we found the sheaf
$\sL_\iota$, asked for in Theorem \ref{in.5}, a): 
\begin{proposition}\label{nc.5}
Under the assumptions made in Theorem \ref{in.5}, assume that \ref{nc.1} holds true for $\V$.
Then, replacing $Y$ by an \'etale covering, there exists an invertible sheaf $\sL_\iota$ with
$\sL_\iota^{n_\iota+1}=\det{\Omega_\iota}$. Moreover the Higgs bundle
$\sL_\iota^{n_\iota+1}\oplus \sL_\iota^{n_\iota+1}\otimes T_\iota$ with Higgs field
$$
\sL_\iota^{n_\iota+1} \>>> \sL_\iota^{n_\iota+1}\otimes T_\iota \otimes \Omega_\iota
$$
induced by the homotheties, is the Higgs field of a variation of Hodge structures $\L_\iota$.
For some rank one local system $\U_\iota$ on $Y$ one has $\V=\U_\iota\otimes \L_\iota$.
\end{proposition}
It remains to consider irreducible $\Q$-sub-variations of Hodge structures
which violate \ref{nc.1}.
\begin{proposition}\label{nc.6}
Let $\V$ be an irreducible $\C$-sub-variation of Hodge structure in $R^1f_*\C_V$, not satisfying the 
Assumption \ref{nc.1}. Then, replacing $Y$ by an \'etale covering
there exists a morphism $\psi:Y\to Y_1$, such that $Y_1$ is a projective manifold, and  $\V=\psi^*\V_1$
for a $\C$-variation of Hodge structures $\V_1$ satisfying again the assumptions made in 
Theorem \ref{in.5}.
\end{proposition}
\begin{proof}
Again we may assume that $\V$ is defined over a number field, and we consider the Weil restriction
$\W_\Q$ of $\V$. If $\V$ violates the Assumption \ref{nc.1}, the all the conjugates of $\V$ violate
\ref{nc.1} and the local system $\W_\Q$ extends to a local system on
$Y$. Moreover, since $\W_\Q$ is a local sub-system of the variation of Hodge structures of a family of Abelian varieties, $\W_\Q$ has a $\Z$ structure. Then $\W_\Q$ is induced by a family
$g:Z\to Y$ of Abelian varieties.

One obtains a morphism $\psi:Y\to Y_1$ for $Y_1$ a closed subscheme of a suitable moduli space 
of polarized Abelian varieties.

Replacing $Y$ by an \'etale covering we may assume that $\psi$ lifts to a fine moduli space, 
hence that $g:Z\to Y$ is the pullback of a family $f_1:X_1\to Y_1$. Let us write
$(G,\varrho)$ for the Higgs bundle of $g:Z\to Y$.

By construction of $Y_1$ the Kodaira-Spencer map for $f_1$ is injective and an isomorphism
in a general point of $Y_1$. Proposition \ref{st.4} 
applied to the original family, implies that the Kodaira-Spencer map of $\W_\Q$ factors through 
$$
G^{1,0} \>>> G^{0,1}\otimes \bigoplus \Omega_j \> \subset >> G^{0,1}\otimes
\psi^*\Omega^1_{Y_1}\> \subset >> \Omega^1_Y(\log S),
$$
where the direct sum in the middle is over some subset of $\{1,\ldots,s\}$, say
over $\{1,\ldots,m\}$.

The sheaves $\bigoplus_{j=1}^m \Omega_j$ and  $\psi^*\Omega^1_{Y_1}$ have
the same rank. On the other hand, the first one is a direct factor 
of $\Omega^1_Y(\log S)$, so both sheaves are isomorphic. 

This implies that $Y_1$ is non-singular. Since $\Omega^1_Y(\log S)$ is nef, the direct
factor $\psi^*\Omega^1_{Y_1}$ is nef, hence $\Omega^1_{Y_1}$, as well.
The universal covering of $Y_1$ has to be the product $M_1\times \cdots \times M_m$, where the $M_i$ are 
factors of the universal covering of $U$. By Proposition \ref{av.1} $\omega_{Y_1}$ must be ample. 

Let $\V_1$ be the direct factor of $R^1f_{1*}\C_{X_1}$ whose pullback is $\V$.
Lemma \ref{ya.6} implies that $\V_1$ satisfies again the numerical 
properties asked for in Theorem \ref{in.5}. In fact,
$$
\ch_1(\Omega_{m+1})^{n_{m+1}}. \cdots . \ch_1(\Omega_s)^{n_s}
$$
restricted to a fibre of $\psi$ is an effective zero cycle, and any (in)equality for 
$\mu_{\omega_{Y_1}}$ or $\delta_{\omega_{Y_1}}$ applied to bundles on $Y_1$ is equivalent to the corresponding
(in)equality for $\mu_{\omega_{Y}(S)}$ or $\delta_{\omega_{Y}(S)}$ applied to the pullback of the bundles.
Obviously the conditions a)--c) in Proposition \ref{in.3} remain true for $\V_1$, if they hold for $\V$. 
\end{proof}
The local system $\V_1$ in \ref{nc.6} will be studied in the next Section.
\begin{corollary}\label{nc.7}
Assume in Theorem \ref{in.5} that there is no \'etale covering of $U$ with a surjective
morphism to a projective manifold of dimension larger than zero.
Then $\V'=\phi^*\V$ or its dual has the Higgs bundle
$$
\phi^*\det(\Omega_i)^{\frac{1}{n_i+1}} \ \oplus \
\big(\phi^*\det(\Omega_i)^{\frac{1}{n_i+1}}\otimes \phi^*\Omega^\vee_i, \ \tau\big).
$$
\end{corollary}
\begin{proof} Since we excluded the existence
of a non-trivial morphism $U\to Y_1$ each irreducible $\C$-sub-variation of Hodge structures $\V$ in $R^1f_*\C_V$ satisfies the Assumption \ref{nc.1}, and \ref{nc.7} follows from \ref{nc.5}.
\end{proof}
\section{Higgs bundles in the compact case} \label{co}
If in Theorem \ref{in.5} $U=Y$ is a compact ball quotient, by \cite[9.1]{Sim} there is a uniformizing projective 
variation of Hodge structures. As in \cite[4.1]{Lo} one can replace $Y$ by an \'etale covering,
such that the natural invertible bundle $\L$ on the ball descends to $Y$, giving an invertible 
sheaf $\sL$ with $\sL^{n+1}=\omega_Y$. Then the uniformizing
variation of Hodge structures has a Higgs bundle of the form
$$
(\sL \ \oplus \ T^1_Y\otimes \sL)
\mbox{ \ or \ } (\Omega_Y^1\otimes \sL \ \oplus \
\sL).
$$
We will extend this result below to the case $s>1$, keeping the assumption $S=\emptyset$.
The next Lemma and its proof are due to F. Bogomolov.
\begin{lemma}\label{co.1}
Let $\phi:G' \to G$ be a finite surjective morphism between $\Q$-algebraic groups,
where $G'$ is an algebraic subgroup of $\Gl(\nu,\Q)$, for some $\nu$. Let
$\Gamma\subset G$ be a finitely generated subgroup. Then there exists a
subgroup $\Gamma'\subset \Gamma$ of finite index, such that the inclusion
$\Gamma' \subset G$ lifts to $\Gamma'\subset G'$.
\end{lemma}
\begin{proof}
Since $\Gamma\subset G$ is finitely generated and since $\phi: G'\to G$ is 
surjective with finite kernel $K$, the pre-image $\phi^{-1}(\Gamma)$ is 
finitely generated. For example, it is generated by $K$ and by the pre-image of any system 
of generators of $\Gamma$.  Since $\phi^{-1}(\Gamma)$ is a finitely generated subgroup in a 
matrix group $G'$, it is well known that $\phi^{-1}(\Gamma)$ is residue finite, i.e. that there exists a sequence 
$\phi^{-1}(\Gamma)=\Gamma_1\supset\Gamma_2\supset\ldots \supset\Gamma_i\supset \ldots$
of subgroups of finite index with
$$
\bigcap_{i=0}^{\infty}
\Gamma_i=\{e\}
$$  
(see for example \cite[Proof of Lemma 6.4]{Sim3}).
Since $K$ is finite, there must exist some $i$ with
$\Gamma_i\cap K=\{e\}$. Then one may choose $\Gamma'=\phi(\Gamma_i)$, and
$\Gamma_i$ lifts $\Gamma'$ to $G'$.
\end{proof}
\begin{lemma}\label{co.2}
Let $U$ be a smooth manifold such that the universal covering $\tilde U$ is a product
$M_1\times \cdots \times M_s$. Assume that for some $i$ the factor $M_i$ is an $n_i$-dimensional complex ball. Let $\Omega_1\oplus \cdots \oplus \Omega_s$ be the corresponding decomposition of the cotangent bundle of $U$. Then, replacing $U$ by an \'etale covering, there exists an invertible sheaf $\sL_i$ on $U$ with
$\sL_i^{n_i+1}=\det(\Omega_i)$, for $n_i=\dim(M_i)=\rk(\Omega_i)$ .
\end{lemma}
\begin{proof} Let us assume that $i=1$. 
Recall the description of the complex ball $M_1$ in \cite[1.8]{Lo}.
Consider a $\C$ vector space $W_1$ of dimension $n_1+1$, equipped with a Hermitian form 
$\psi_1$ of signature $(1,n_1)$. Then $M_1\subset \BP(W_1)$ is the open subset
defined by $\psi_1(w,w)>1$. The action of $\pi_1(U,*)$ on $M_1\times \cdots \times M_s$
is given by 
$$
\rho:\pi_1(U,*) \>>> {\rm Aut}(M_1)\times \cdots \times {\rm Aut}(M_s)
$$
and the first factor of the right hand side is $\BP {\rm U} (\psi_1) \hookrightarrow \BP {\Sl}(n_1+1)$.

Replacing $\pi_1(U,*)$ by a subgroup of finite index, hence replacing $U$ by some \'etale 
covering, Lemma \ref{co.1} allows to lift $\rho$ to
$$
\rho':\pi_1(U,*) \>>> {\Sl}(n_1+1)\times {\rm Aut}(M_2)\times \cdots \times {\rm Aut}(M_s).
$$ 
For $\gamma\in \pi_1(U,*)$ let $\gamma_1$ denote the first component of
$\rho'(\gamma)$. Then up to the multiplication with a constant $\gamma_1$ 
lies in the unitary group for $\psi_1$. Since it lies in ${\Sl}(n_1+1)$ the constant has to have
absolute value one. Hence $\gamma_1$ lies in ${\rm SU}(n_1+1)$. As in \cite[4.1]{Lo} we may 
replace $U$ again by an \'etale covering, and assume that the subgroup of $\C^*$, generated by the 
eigen values of $\rho'(\pi_1(U,*))$ is torsion free. The group $\rho'(\pi_1(U,*))$ acts on the line bundle $\sN_1={\rm pr}_1^*\sO_{\BP(W_1)}(-1)_{M_1}$
and it descends to a line bundle $\sL^{-1}_1$ on the quotient. By \cite[4.1]{Lo}
the canonical sheaf ${\rm pr}_1^*\omega_{\BP(W_1)}|_{M_1}$ is ${\rm SU}(\psi_1)$-equivariantly
isomorphic to $\sN_1^{n_1+1}$. The latter descends to $\det(\Omega_1)$ on $U$. 
\end{proof}
\begin{corollary}\label{co.3}
In \ref{co.2} assume that $S=\emptyset$, hence $U=Y$ projective. Then the Higgs bundle 
$$
F=(\sL_i^{-1}\otimes \Omega_i \oplus \sL_i^{-1}, {\rm id}_{\sL_i^{-1}\otimes \Omega_i})
$$
is the Higgs bundle of a complex variation of Hodge structures $\L_i$.
\end{corollary}
\begin{proof}
Obviously the first Chern class of $F$ is zero, and as in the proof of \ref{pr.5}
one sees that 
$$
\Delta(F)=\Delta(\Omega_i \oplus \sO_U)=2\cdot (n_i+1)\ch_2(\Omega_i)-n_i\ch_1(\Omega_i)^2.
$$ 
If $n_i=1$ one has $\Delta(F)=0$, hence $\delta(\sF)=0$, as well. For $n_i>0$ one finds by Theorem
\ref{ya.4} and Lemma \ref{ya.6}, vii), that $\delta(F)=0$. Then \ref{co.3} follows from \cite{Sim2}.
\end{proof}
\begin{corollary}\label{co.4}
Under the assumptions made in Theorem \ref{in.5} assume that $S=\emptyset$, hence that
$f:X\to Y$ is smooth. Replacing $Y$ by an \'etale covering, there exists a unitary local system $\U_i$, regarded as a variation of Hodge structures of bidegree $(0,0)$ such that $\V$ or its dual is isomorphic to $\U_i\otimes \L_i$.
\end{corollary}
\begin{proof} 
We may assume that $\ell=\rk(E^{1,0}) \leq \rk(E^{0,1})$. Then 
$$
\ch_1(E^{1,0})=-\ch_1(E^{0,1})= - n_i\cdot\ch_1(E^{1,0}) + \ell \cdot \ch_1(\Omega_i),
$$
hence $\displaystyle\ch_1(E^{1,0})= \frac{\ell}{n_i+1}\cdot\ch_1(\Omega_i)=\ell\cdot \ch_1(\sL_i)$.
One finds $\ch_1(E^{1,0}\otimes \sL_i^{-1})=0$ and $\delta(E^{1,0}\otimes \sL_i^{-1})=0$.

Then $\sU=E^{1,0}\otimes \sL_i^{-1}$ together with the trivial Higgs field must be the Higgs bundle of
a unitary bundle $\U_i$. The explicite descriptions of the Higgs fields
of $\L_i\otimes \U_i$ and $\V$ in \ref{co.3} and \ref{in.3} show, that
the Higgs fields of
$$
(\sU,0)\otimes (\sL_i \ \oplus \ \sL_i\otimes T_i, \rho)
\mbox{ \ \ and  \ \ }
(E^{1,0}\oplus E^{0,1}=E^{1,0}\otimes T_i, \theta)
$$
coincide, so $\L_i\otimes \U_i \cong \V$.
\end{proof}

\begin{proof}[Proof of Theorem \ref{in.5}] We have shown already in Section \ref{pr}
that the assumptions made in Theorem \ref{in.5} imply that $M_i$ is a complex ball.
If the Assumption \ref{nc.1} holds true, we verified the conditions a) and b) in \ref{in.5}
in Corollary \ref{nc.5}.

Otherwise, we know by Proposition \ref{nc.6} that $\V$ is the pullback of a variation of Hodge structures
$\V_1$ under a surjection $\psi:U\to Y_1$ with $Y_1$ a projective manifold. Moreover $Y_1$ satisfies again the 
assumptions made in Theorem \ref{in.5}. So Corollary \ref{co.3} applies to $\V_1$ and $Y_1$, and 
the conditions a) and b) hold on $Y_1$. Obviously they are compatible with pullback.
\end{proof}
In fact, we did not use up to now that $\mu=\mu_{\omega_Y(S)}$, we used just the assumptions 
stated in \ref{sl.5}. So we obtained:
\begin{variant}\label{co.5}
The Theorem \ref{in.5} remains true for $\mu=\mu_\sN$ and $\delta=\delta_\sN$, provided $\sN$ is an ample invertible sheaf and $\Omega_Y^1(\log S)$ is $\mu_\sN$-poly-stable.
\end{variant}
\section{Generalized Hilbert modular varieties and surfaces}\label{su}
We will call $U$ in the Set-up \ref{in.3} a generalized Hilbert modular variety,
if the universal covering $\tilde{U}$ is the product of complex one dimensional balls.
We allow $U$ to be a product $U_1\times U_2$ of two generalized Hilbert modular surfaces, for example
$U$ could be the product of curves of genus $g>1$. 
\begin{corollary}\label{su.1} Assume that $\mu$ is chosen according to Assumption \ref{sl.5}.
Assume that
$$
\mu(R^1f_*\C_V) = \mu(\Omega^1_Y(\log S))
\mbox{ \ \ and \ \ } \delta(R^1f_*\C_V)=0.
$$
Then $U$ is a generalized Hilbert modular variety, i.e. its universal covering is isomorphic to the
product of one dimensional complex balls.
\end{corollary}
\begin{proof}
Let $\V_1 \oplus \cdots \oplus \V_\nu$ be the decomposition of $R^1f_*\C_V$ in irreducible $\C$-sub-variations of Hodge structures, and let $E_i$ be the Higgs bundle of $\V_i$. Hence $F=E_1\oplus\cdots\oplus E_\nu$ is the Higgs bundle of $R^1f_*\C_V$. Lemma \ref{sl.1}, a), 
and Proposition \ref{sl.6} imply that 
$$
\mu(\Omega^1_Y(\log S))=\mu(F^{\doubledot}) \leq {\rm Max}\{\mu(E_i^{\doubledot}); i=1,\ldots,\nu\} \leq \mu(\Omega^1_Y(\log S)).
$$
So these are equalities, and applying \ref{sl.1}, a), again, one finds that each
of the $E_i$ satisfies the Arakelov equality, hence they are all semi-stable of the same slope.
Moreover, for all $i$
\begin{equation}\label{eqsu.1}
\frac{\rk(E_i^{1,0})}{\rk(E_i^{0,1})}=\frac{\rk(E_1^{1,0})}{\rk(E_1^{0,1})}.  
\end{equation}
By assumption $\delta(F^\doubledot)=0$ and Lemma \ref{st.3}, d), implies that
$\delta(E_i^\doubledot)=0$. So Theorem \ref{in.5} or its variant \ref{co.5} applies.
By the explicite description of the irreducible direct factors $\V_i$ given there, (\ref{eqsu.1}) can only hold for $\V_i$ and for its complex conjugate, if both are isomorphic. So
$\rk(E^{1,0}_i)=\rk(E^{0,1}_i)$ and $U$ is a generalized Hilbert modular variety.
\end{proof}

Note that this result, as well as Theorems \ref{in.5} and \ref{in.6} rely on the conditions ($*$) and ($**$) hidden in the Set-up \ref{in.3}. So strictly speaking, as long as the announced
article by Sun and Yau does not exist, the results only apply if $\omega_Y(S)$ is ample.
This condition excludes in particular all generalized Hilbert modular varieties with
$U\neq Y$. For surfaces one can replace the polarization $\mu$ by 
$\mu_N$, for some small twist $N$ of $\omega_Y(S)$.

\begin{setup}\label{su.2}
$Y$ is a surface, $U\subset Y$ the complement of a normal crossing divisor. $\Omega_Y^1(\log S)$ is nef, and $\omega_Y(S)$ is ample with respect to $U$. Let $\sH$ be an ample invertible sheaf, 
$$
N=\ch_1(\omega_Y(S))+\epsilon\cdot\ch_1(\sH)
$$ 
and $\mu_\epsilon=\mu_N$. Assume that there exists some $\epsilon_0>0$ such that for all $\epsilon_0 \geq \epsilon \geq 0$ and for all $m>0$ the sheaf $S^m(\Omega^1_Y(\log S))$ is $\mu_\epsilon$-poly-stable. 

Let $f:V\to U$ be a smooth family of polarized $g$-dimensional Abelian varieties with uni-potent local monodromy around the components of $S$, such that the induced morphism $\varphi:U\to \sA_g$ is generically finite. 
\end{setup}
Let us first verify, that for each $Y$ we can find an ample invertible sheaf $\sH$ as in the Set-up \ref{su.2}.

By Theorem \ref{ya.4} we know that $S^m(\Omega_Y^1(\log S))$ is $\mu_0$-poly-stable. Note that the case c), iii), can not occur. In fact, by \cite{Y} if $S^m(\Omega^1_Y(\log S))$ is not stable, for some $m>0$, then $\Omega^1_Y(\log S)= \Omega_1\oplus \Omega_2$ with $\Omega_i$ invertible.

Assume first that $S^m(\Omega^1_Y(\log S))$ is stable with respect to $\mu_0$. As in the proof of
b) in Lemma \ref{ya.3} one finds that $\Omega^1_Y(\log S)$ is $\mu_\epsilon$-semi-stable
for $\epsilon$ sufficiently small. If $\Omega^1_Y(\log S)$ is not $\mu_\epsilon$-stable, 
there exists a subsheaf $\sG$ of $S^m(\Omega_Y^1(\log S)$, with
$$
\mu_\epsilon(\sG)=\mu_0(\sG) + \epsilon\cdot\upsilon(\sG).\ch_1(\sH) = \mu_0(\Omega^1_Y(\log S))
+ \epsilon\cdot\upsilon(\Omega^1_Y(\log S)).\ch_1(\sH).
$$
By assumption $\mu_0(\sG) < \mu_0(\Omega^1_Y(\log S))$, hence for all $\epsilon$ sufficiently small
$$
0<\mu_0(\Omega^1_Y(\log S))-\mu_0(\sG)= \epsilon\cdot (\upsilon(\sG)-\upsilon(\Omega^1_Y(\log S))).\ch_1(\sH),
$$
a contradiction. So $S^m(\Omega^1_Y(\log S))$ remains $\mu_\epsilon$-stable.

Assume next that $\Omega^1_Y(\log S)= \Omega_1\oplus \Omega_2$. Let $\sA$ be any ample invertible sheaf on $Y$.
We know by Lemma \ref{ya.6} that $\ch_1(\Omega_1)^2=\ch_1(\Omega_2)^2=0$ and
$\beta=\ch_1(\Omega_1).\ch_1(\Omega_2) >0$. If $\ch_1(\sA).\ch(\Omega_1) > \ch_1(\sA).\ch(\Omega_2)$
choose 
$$
\alpha=\ch_1(\sA).\ch(\Omega_1) - \ch_1(\sA).\ch(\Omega_2)+ \beta.
$$
Then $\sH=\sA^\beta\otimes \Omega_1^\alpha \otimes \Omega_2^\beta$ is ample and
$$
\ch_1(\sH).\ch_1(\Omega_1)=\ch_1(\sH).\ch_1(\Omega_2)=\beta\cdot \ch_1(\sA).\ch_1(\Omega_1)+\beta^2.
$$
So $S^m(\Omega_Y^1(\log S))$ as the direct sum invertible sheaves of the same slope
is $\mu_\epsilon$-poly-stable. So we obtained:
\begin{lemma}\label{su.3}
Let $Y$ be a non-singular projective surface and $U$ the complement of a normal crossing divisor 
$S$. If $\omega_Y(S)$ is nef and ample with respect to $U$, than one can find an ample invertible sheaf $\sH$ and some $\epsilon_0>0$ such that for all $m>0$ and for $\epsilon_0\geq \epsilon \geq 0$ the sheaf $S^m(\Omega_Y^1(\log S))$ is $\mu_\epsilon$-poly-stable. 
\end{lemma}

On obtains the following variant of Theorems \ref{in.5} and \ref{in.6}:
\begin{variant}\label{su.4}
In the Set-up \ref{su.2} one has for some $\epsilon_0>0$ and all $\epsilon_0\geq \epsilon\geq 0$
and all non-unitary irreducible $\C$-sub variations of Hodge structures $\V$ of $R^1f_*\C_V$ with
Higgs bundle $(E^{1,0}\oplus E^{0,1}, \theta)$ the inequality
$\mu_\epsilon(\V)\leq \mu_\epsilon(\Omega_Y^1(\log S))$. If equality holds, the sheaves
$E^{1,0}$ and $E^{0,1}$ are both semi-stable and $\delta(\V)\geq 0$. In addition one has:
\begin{enumerate}
\item[I.] If for all $\V$ one has the equalities $\mu_\epsilon(\V)= \mu_\epsilon(\Omega_Y^1(\log S))$ and 
if $\delta(\V)=0$ then $U$ is either a ball quotient or a generalized Hilbert modular surface.
\item[II.] Assume that $\Omega^1_Y(\log S)$ is the direct sum of two line bundles $\Omega_1$ and $\Omega_2$ with $\mu_\epsilon(\V)=\mu_\epsilon(\Omega_1)=\mu_\epsilon(\Omega_2)$ for all $\C$-sub variations of Hodge structures
in $R^1f_*\C_V$. Then $U$ is a generalized Hilbert modular surface. 
\item[III.] Assume that $\mu_\epsilon(R^1f_*\C_V)= \mu_\epsilon(\Omega_Y^1(\log S))$ and 
$\delta(R^1f_*\C_V)=0$. Then $U$ is a generalized Hilbert modular surface and
$\Omega^1_Y(\log S)$ is the direct sum of two line bundles $\Omega_1$ and $\Omega_2$
of the same slope.
\item[IV.] In II) or III), replacing $U$ by an \'etale covering, there exist invertible sheaves $\Omega_i^{\frac{1}{2}}$, and
\begin{equation}\label{eqsu.2}
R^1f_*\C_V =\W_\Q\otimes \C= \L_1\otimes \U_1\oplus\L_2\otimes\U_2,
\end{equation}
where the $\L_i$ are the uniformizing variations of Hodge structures with Higgs bundle
$$
\big( \Omega_i^{\frac{1}{2}} \oplus \Omega_i^{\frac{-1}{2}}, \ \theta:\Omega_i^{\frac{1}{2}}
\>>> \Omega_i^{\frac{-1}{2}} \otimes \Omega_i \subset \Omega_i^{\frac{-1}{2}} \otimes \Omega^1_Y(\log S) \big),
$$
and where the $\U_i$ are unitary local systems.
\item[V.] If in I) $\Omega_Y^1(\log S)$ is $\mu_0$-stable, $U$ is a ball quotient. Replacing $U$ by an \'etale covering, there exists an invertible sheaf $\omega_Y(S)^{\frac{1}{3}}$, and
$$
R^1f_*\C_V= \L\otimes \U \oplus \bar\L\otimes \bar\U
$$ 
for a unitary local system $\U$ concentrated in bidegree $(0,0)$ and for $\L$ with Higgs bundle
$$
(\omega_Y(S)^{-\frac{1}{3}}\otimes \Omega^1_Y(\log S) \ \oplus \ \omega_Y(S)^{-\frac{1}{3}}, \ {\rm id}). 
$$
\end{enumerate}
\end{variant}
\begin{proof}
I) has been shown in Proposition \ref{sl.6} and II) is a special case of \ref{su.1}. Using I), Part III) is obvious. The explicite form of the variation of Hodge structures in IV) and V)
follows from Proposition \ref{st.4} and from Variant \ref{co.5}. 
\end{proof}
Let us consider the sheaves in part IV) of \ref{su.4} a bit closer.
\begin{lemma}\label{su.5}
The decomposition (\ref{eqsu.2}) is defined over a finite Galois 
extension $K$ of $\Q$ with Galois group $G$; i.e. $\L_i$ and $\U_i$ are defined over $K$
and the decomposition exists for $\W_\Q\otimes K$. 
\end{lemma}
\begin{proof}
In order to see that $\L_i\otimes \U_i$ is defined over $\bar \Q$, one just has to repeat the argument
used to prove \ref{nc.2} or \cite[3.3]{VZ1}. For the tensor product decompositions one argues as in 
\cite[3.7, iii)]{VZ1}).
\end{proof}
\begin{lemma}\label{su.6}
Assume that for all $\tau\in G\setminus
\{{\rm id}\}$ the local system $\L_1^\tau$ is unitary.
Then the representation $\rho_1$ of $\L_1$ is discrete, and 
some \'etale covering of $U$ is a product of two curves. 
\end{lemma}
\begin{proof}
Consider the adjoint representation $\END (\W_\Q)$. Obviously
it has a $\Z$-structure. Moreover $\END _0(\L_1)$ is a direct factor
of $\END (\W_\Q)\otimes \bar\Q$. Hence for the ring $\sO$ of integers in some algebraic number 
field $K$ the system $\END _0(\L_1)$ inherits an $\sO$-structure.

By assumption, the Weil restriction $\sW( \END _0(\L_1))$ contains  
only one non compact factor, $\END _0(\L_1)$.  
Since $\sW( \END _0(\L_1))$ has a $\Z$-structure $\END _0(\L_1)$ must be discrete.

Consider the adjoint representation $\Sl_2 \to {\rm Aut}(sl_2)$.  
Since its kernel is finite, and since $ \END _0(\L_1)$ is discrete, one finds $\L_1$ to be discrete.
\end{proof}
\begin{corollary}\label{su.7}
If in the decomposition (\ref{eqsu.2}) $\rk(\U_1)\neq \rk(\U_2)$, then
some \'etale covering of $U$ is the product of two curves.
\end{corollary}

\begin{proof} Assume that $\rk\U_1=\nu < \rk\U_2=g-\nu$.
Consider the  $g$th wedge product
$$
\bigwedge^g(\W)=\bigwedge^g(\L_1\otimes \U_1 \oplus \L_2\otimes \U_2).
$$
It has one direct factor $S^{\nu}(\L_1)\otimes S^{n-\nu}(\L_2)$ 
and all other direct factors are tensor products of
$S^s(\L_1)$, $S^t(\L_2)$ and of unitary local systems, where $s\leq \nu$ and
$t \leq g-\nu$. 

For $\tau\in {\rm Gal}(\bar\Q/\Q)$, 
one has $\L_2^\tau\not\cong \L_1$. Otherwise, since $\bigwedge^g(R_1f_*\bar\Q_X)$
is defined over $\Q$, it would have a direct factor
of the form $S^{g-\nu}(\L_1)\otimes \U$, contradicting $\nu<g-\nu$.

Hence the local system $\L_2^\tau$ is either isomorphic to $\L_2$ or it is unitary.
The Weil scalar restriction $\sW(\L_2)$ has an $\Q-$structure, and
it is the direct sum over all local systems, conjugate to $\L_2$. Hence
except of $\L_2$ all the direct factors of $\sW(\L_2)$ are unitary.
By Lemma \ref{su.6} some finite \'etale covering of $U$ is a product of two curves,
contradicting the assumption made.
\end{proof}

Recall that a generalized Hilbert modular surface $U$ is a Hilbert modular surface in the usual 
sense, if and only if the local system $\L_1\oplus \L_2$ in Lemma \ref{su.5} is defined
over $\Q$, whereas each of the $\L_\iota$ is defined over a real quadratic extension
$K$ of $\Q$. Moreover $\L_1\oplus \L_2$ has a $\Z$ structure, hence it is the variation of Hodge 
structures of a family $h:Z\to U$ of Abelian surfaces. 
As well known, for such a Hilbert modular surface one has $S=Y\setminus U\neq \emptyset$.

\begin{corollary}\label{su.8}
Assume that $U$ is a generalized Hilbert modular surface with $S \neq \emptyset$,
and that no \'etale covering of $U$ is the product of two curves.
Then, replacing $U$ by an \'etale covering, the unitary local systems $\U_\iota$ in
(\ref{eqsu.2}) are trivial, and $\L_\iota$ is defined over a real quadratic extension
$K$ of $\Q$. In particular $U$ is a Hilbert modular surface, and $f:V\to U$ is isogenous
to $Z\times_U \cdots \times_U Z$. The fibres of $g:Z\to U$ have real multiplication.
\end{corollary}
\begin{proof}
Since no finite \'etale covering of $Y\setminus S$ is a product of two curves, the Galois group 
${\rm Gal}(\bar\Q/\Q)$ permutes $\L_1$ and $\L_2$. By Proposition \ref{nc.5},
we may assume that $\U_1$ and $\U_2$ are both trivial. Then $\L_1\oplus \L_2$
has a $\Z$ structure, and $V\to U$ is isogenous to
$Z\times_U \cdots \times_U Z$ for a family of Abelian surfaces $g:Z \to U$. 

The general fibre $F_\eta$ of $g$ must be a simple Abelian surface, since otherwise it would be 
isogenous to the product of elliptic curves and the $\L_i$ would be defined over $\Q$.
Moreover, since
$$
\END(\L_1\oplus \L_2)= \C^2 \oplus \L_1\otimes \L_2^\vee \oplus
\L_1^\vee\otimes \L_2,
$$ 
the dimension of $\End(F_\eta)\otimes \Q$ is two. By the well known classification
of endomorphisms of Abelian surfaces (see for example \cite[5.5.7]{BL}) this implies that $F_\eta$ 
either has real multiplication, or that $\End(F_\eta)\otimes \Q$ is an imaginary
quadratic extension of $\Q$. However, by \cite[Example 6.6 in Chapter 9]{BL}
there are only finitely many of surfaces of the second type.
\end{proof}
Of course, there are generalized Hilbert modular surfaces with
$S=\emptyset$ (see \cite{Sh} or \cite{Gr}), and for some of them 
the variation of Hodge structures has rank bigger than two, hence the unitary systems
$\U_i$ will be non-trivial.\vspace{.2cm}

Variations of Hodge structures, uniformizing certain ball quotients $U=Y\setminus S$,
have been constructed in \cite{DM} with $S\neq\emptyset$ and with $S=\emptyset$. 
For example, the moduli scheme of $5$ points in $\BP^1$ is an example of the second kind
(in \cite[p. 86]{DM} there seems to be a misprint in example 5). In \cite{VZ2} it is 
shown that this example, a compact two dimensional ball quotient in the moduli scheme of
$4$-dimensional Jacobian varieties, is a Shimura variety. We will give a generalization in the
next two sections.
 
\section{The decomposition of certain wedge products}\label{rg}
In the next two sections we will use the assumptions made in Theorem \ref{in.6}. To show that $U$ is a Shimura variety we will determine the possible Hodge cycles for self products of $f:V\to U$, hence the possible trivial (or unitary) local sub-systems
in wedge products of the local systems described in Theorem \ref{in.5}, b).
In this section we will just state one application, the rigidity of the family of Abelian varieties
in Theorem \ref{in.6}. The Lemma \ref{rg.4} will be needed again in Section \ref{sh}.

Let $\V$ be a variation of Hodge structures of weight $k$ with Higgs bundle
$$
\big( E=\bigoplus_{p+q=k} E^{p,q},\theta=\bigoplus_{p+q=k} \theta_{p,q}\big).
$$
Let $q_0$ be the smallest integer with $E^{k-q_0,q_0}\neq 0$
The $i$-th iterated cup product with the Kodaira-Spencer map defines a morphism,
the Griffiths-Yukawa coupling,
\begin{multline*}
\theta^i:E^{k-q_0,q_0} \>\theta_{k-q_0,q_0}>> E^{k-q_0-1,q_0+1}\otimes \Omega^1_S
\>\theta_{k-q_0-1,q_0+1}>>\\ E^{k-q_0-2,q_0+2} \otimes S^2(\Omega^1_S) \>>> \cdots
\>\theta_{k-q_0-i+1,q_0+i-1}>> E^{k-q_0-i,q_0+i}\otimes S^i(\Omega^1_S).
\end{multline*}
We define its length or the length of $\V$ to be
$$
\varsigma(\V)={\rm Min}\{i\geq 1 ; \ \theta^i=0\}-1.
$$
If $\V$ and $\W$ are two variations of Hodge structures, one has
$$\varsigma(\V\otimes \W)=\varsigma(\V) + \varsigma(\W).$$

Let $\L_j$ be one of the uniformizing variations of Hodge structures in Addendum
\ref{in.5}, say with Higgs bundle $\sL_j \oplus \sL_j \otimes T_j$,
where $\Omega_j$ is one of the stable direct factors of $\Omega^1_Y(\log S)$, of rank $n_j$,
where $T_j=\Omega_j^\vee$ and where $\sL_j^{n_j+1}=\det(\Omega_j)$. The Higgs field is given by the
homotheties $\sO_Y\to \sE nd(\Omega_j)$, tensorized with ${\rm id}_{\sL_j}$.
\begin{lemma}\label{rg.1}
For $1\leq k \leq n_j$, the variations of Hodge structures $\bigwedge^k\L_j$
and $\bigwedge^k\L_j^\vee$ are concentrated in two degrees, and the Higgs fields of $\bigwedge^k\L_j$
and $\bigwedge^k\L_j^\vee$ are given by injections
\begin{gather*}
\sL_j^k \otimes \bigwedge^{k-1}T_j \>>> \sL_j^k \otimes \bigwedge^k T_j \otimes \Omega_j
\mbox{ \ \ and \ \ }\\
\sL_j^{-k}\otimes\bigwedge^{k} \Omega_j \>>> \sL_j^{-k} \otimes \bigwedge^{k-1}\Omega_j\otimes \Omega_j
\end{gather*}
respectively. For $k=n_j+1$ the local systems $\bigwedge^k\L_j$ and $\bigwedge^k\L_j^\vee$ are
both of rank one and of bidegree $(1,n_j)$ and $(n_j,1)$.
\end{lemma}
\begin{proof}
The Higgs fields of $\bigwedge^k\L_j$ and $\bigwedge^{k'}\L_j^\vee$ are induced by natural
direct factors
$$
\bigwedge^{k-1} T_j \subset \bigwedge^{k}T_j\otimes \Omega_j \mbox{ \ \ and \ \ }
\bigwedge^{k'} \Omega_j \subset \bigwedge^{k'-1}\Omega_j\otimes \Omega_j
$$
(see for example \cite[(6.9), p. 79]{FH})
tensorized with $\sL_j^k$ and $\sL_j^{-k'}$, respectively.
\end{proof}
\begin{corollary}\label{rg.2}
Let $\U$ be a unitary local system of rank $\ell$. Then
$$
\varsigma(\bigwedge^{k}(\L_j\otimes \U))=
\varsigma (\bigwedge^{k}(\L_j\otimes \U)^\vee)=\ell,
$$
for $\ell \leq k \leq \ell\cdot n_j$ whereas
$$
1\leq \varsigma(\bigwedge^k(\L_j\otimes \U))=\varsigma (\bigwedge^k(\L_j\otimes \U)^\vee) < \ell
$$
for $0 < k < \ell$ and for $\ell\cdot n_j < k < \ell\cdot(n_j+1)$.
For $k=\ell\cdot(n_j+1)$ both, $\bigwedge^k(\L_j\otimes \U)$ and $\bigwedge^k(\L_j\otimes \U)^\vee$
are unitary local systems, concentrated in bidegree $(\ell,\ell\cdot n_j)$
and $(\ell\cdot n_j,\ell)$, respectively.
\end{corollary}
\begin{proof} It is sufficient to consider $\bigwedge^k(\L_j\otimes \U)$.
The length of a Higgs field can be calculated in a general point, so by abuse of notations
we may assume that $\U=\C^\ell$, and
$$
\bigwedge^k(\L_j\otimes \U)=\sum_{k_1+\cdots+k_\ell=k}\bigotimes_{i=1}^\ell \bigwedge^{k_i}\L_j.
$$
By \ref{rg.1}
$$
\varsigma\big(\bigotimes_{i=1}^\ell \bigwedge^{k_i}\L_j\big)=\ell
$$
if and only if none of the $k_i$ is zero or equal to $n_j+1$. This will hold true for at least one
of the direct factors, whenever $\ell\leq k \leq \ell\cdot n_j$.
\end{proof}

Recall that by Proposition \ref{in.4} each irreducible $\C$-sub-variation of Hodge structures
in $R^1f_*\C_V$ has a Higgs field involving only one of the stable direct factors of
$\Omega^1_Y(\log S)$. Theorem \ref{in.5} allows to write for a $\Q$-sub-variation $\W_\Q$ of Hodge structures in $R^1f_*\Q_{V}$
\begin{gather}\label{eqrg.1}
\W=\W_\Q\otimes \C=\bigoplus_{i=1}^{s}\hat\V_i\mbox{ \ \ with \ \ }\\ \notag
\hat\V_i =\left\{
\begin{array}{lll}
\V_i &\mbox{for}& i=1,\ldots,s''\\
\V_i\oplus\V_i^\vee, &\mbox{for}& i=s''+1,\ldots,s,
\end{array}\right.
\end{gather}
for the local systems $\V_i=\L_i\otimes \U_i$ with $\L_i$ as in Theorem \ref{in.5} and with $\U_i$ unitary. We denote the natural antisymmetric form on $R^1f_*\Q_V$ by $Q$ and
we write $\ell_i=\rk(\U_i)$.
\begin{lemma}\label{rg.3} \
\begin{enumerate}
\item[a.] The local systems $\hat\V_i$ and the decomposition (\ref{eqrg.1}) are defined over a
a totally real number field $K$ and are orthogonal with respect to the form $Q$.
Moreover we may choose $K$ such that $\W$ decomposes as a direct sum of irreducible 
$K$-sub-variation of Hodge structures $\T_j$ which remain irreducible over $\R$.
\item[b.] In a) assume that $\T_j$ decomposes over $\C$ as the direct sum
of two non-trivial sub-variations of Hodge structures.
Then this decomposition is defined over $K(\sqrt{a})$ for some $a\in K$, totally negative.
\item[c.] For $i=s''+1,\ldots,s$ there exist totally negative
elements $a_i\in K$ such that $\V_i$ is defined over $K(\sqrt{a_i})$. The involution
$\iota_i$ of $K(\sqrt{a_i})$ over $K$ interchanges $\V_i$ and $\V_i^\vee$.
\end{enumerate}
\end{lemma}
\begin{proof}
As in the proof of Lemma \ref{nc.2} we start by copying the argument from \cite[3.3]{VZ1}.

Let $\W$ be any variation of Hodge structures defined over a totally real number field
$K_0$, and allowing a decomposition as the one in (\ref{eqrg.1}) over $\R$.
Consider a family $\{\M_t\}_{t\in \Delta}$ of local sub-systems of $\W$
defined over a small disk $\Delta$, with $\M_0=\hat\V_s$. Let us write
$$
\X=\bigoplus_{i=1}^{s-1}\hat\V_i
$$
For $t\in \Delta$ let $(F_t,\theta_t)$ denote the Higgs bundle
of an irreducible direct factor of $\M_t$. Then
$F_t^{1,0} \to F_t^{0,1}\otimes \Omega_s$
is injective for $t$ sufficiently small. If the composite
$$
\rho:\M_t \> \subset >> \W \>>> \X
$$
is non-zero, the complete reducibility of local systems coming from
variations of Hodge structures implies that $\M_t$ and $\X$
contain a common direct factor. Since the Higgs field
$G^{1,0}\to G^{0,1}\otimes \Omega^1_Y(\log S)$
of $\X$ factors through
$$
G^{1,0}\>>> G^{0,1}\otimes \bigoplus_{i=1}^{s-1}\Omega_i\>>> G^{0,1}\otimes \Omega_Y^1(\log S),
$$
one obtains a contradiction.

So $\M_0$ is rigid as a local sub-system, hence it can be defined over $\bar\Q$.
The complex conjugation maps the local system $\V_s$ to its dual, hence to
$\V_s$ if $s''=s$ or to $\V_s^\vee$, otherwise. Then $\hat\V_s$ is invariant under complex conjugation,
hence it can be defined over $L\subset \R\cap\bar\Q$.

As in the proof of \cite[3.3]{VZ1} this implies that
$\X$ is also defined over $L$, and that the splitting $\W=\V_s\oplus \X$ can
be chosen to be orthogonal. By induction we may assume that the decomposition
(\ref{eqrg.1}) is defined over $L$ and orthogonal.

If $\hat\V_s$ decomposes as a direct sum of irreducible $\C$-sub-variations of Hodge structures,
\cite[3.2]{VZ1} allows us to choose the decomposition to be defined over $\bar\Q$.

Taking the sum over complex conjugates, we obtain a decomposition over $\R\cap \bar\Q$ in factors, which 
remain irreducible over $\R$. Enlarging $L$ we will assume that this decomposition is defined over $L$.

Let $\T$ be an irreducible $L$-sub-variation of Hodge structures in $\V_s$. So
$\bar\T=\T$. Since $Q(v,\bar v)\neq 0$ for all local sections of $\W$, one
has $\bar\T=\T^\vee$, and the restriction $Q$ to $\T$ is non degenerate.

Let $\sigma:L\to \C$ be any embedding, and $\bar\sigma$ its conjugate.
So $\T^{\bar\sigma}$ is equal to $\overline{\T^\sigma}$. Since
$Q$ is defined over $\Q$, one finds
$$
\T^{\bar\sigma}=(\T^{\sigma})^\vee=(\T^\vee)^\sigma=\T^\sigma
$$
and $\T^\sigma$ is defined over $\sigma(L)\cap \R$. In different terms,
if $K\subset L$ is a minimal field of definition of $\T$, it must be totally real.

Again, \cite[3.3]{VZ1} allows to find a splitting of $\T$ in $\W$ which is orthogonal
and defined over $K$. By induction on the rank of $\W$ one obtains \ref{rg.3}, a).

Assume now, that $\T=\T'\oplus\T''$ is a non-trivial decomposition over $\C$.
Then $\T'$ can be defined over some quadratic extension
$K(\sqrt{a})$. It remains to verify that $a$ is totally negative.
If not, there is an embedding $\gamma:K(\sqrt{a})\to \R$, and ${\T'}^\gamma$
is defined over $K(\sqrt{a})$. The above argument for this variation of
Hodge structures, tells us, that it is defined over a totally real subfield, hence over
$K$. Then $\T'$ is defined over $\R$, a contradiction.

Finally part c) follows from a) and b).
\end{proof}
The local system $\bigwedge^k\W$ decomposes over $\bar\Q$ as a direct sum of local systems of the form
$$
\W_{\underline{k}}=\big(\bigotimes_{i=1}^{s''}\bigwedge^{k_i}\V_i\big)\otimes
\big(\bigotimes_{i=s''+1}^{s}\bigwedge^{k_i}\V_i\otimes \bigwedge^{k'_i} \V_i^\vee\big)
$$
for some tuple $\underline{k}=(k_1,\ldots , k_{s''}, k_{s''+1}, k'_{s''+1}, \ldots , k_s, k'_s)$,
with
\begin{equation}\label{eqrg.2}
\sum_{i=1}^{s''}k_i+\sum_{s''+1}^s(k_i+k'_i)=k.
\end{equation}
\begin{lemma}\label{rg.4}
Assume that $\varrho \in H^0(Y, \W_{\underline{k}})$ is non-zero. Then $\varrho$ is concentrated
in one bidegree $(p=p(\varrho),q=q(\varrho))$. Moreover,
\begin{enumerate}
\item[a.] if $k_i=k'_i$ for $s'' < i \leq s$, then $p=q$.
\item[b.] otherwise, $p(\varrho)=q(\bar{\varrho})$ and $q(\varrho)=p(\bar{\varrho})$, where
$\bar\varrho$ is the complex conjugate of $\varrho$.
\end{enumerate}
\end{lemma}
\begin{proof}
Let $(F_i,\tau_i)$ denote the Higgs bundle of $\bigwedge^{k_i}(\L_i\otimes \U_i)$, hence
$$
F_i^{m_i,k_i-m_i}=\bigwedge^{m_i}(\sL_i\otimes \U_i) \otimes \bigwedge^{k_i-m_i}(\sL_i\otimes T_i \otimes \U_i).
$$
For $s'' < i \leq s$ we write $(F'_i,{\tau'}_i)$ for the Higgs bundle of
$\bigwedge^{k'_i}(\L_i\otimes \U_i)^\vee$, hence
$$
F_i^{m'_i,k'_i-m'_i}=\bigwedge^{m'_i}(\sL_i^\vee\otimes \U_i^\vee \otimes \Omega_i) \otimes 
\bigwedge^{k'_i-m'_i}(\sL_i^\vee\otimes \U_i^\vee).
$$
The section $\varrho$ defines a local sub-system $\C\subset \W_{\underline{k}}$, hence a
direct factor $\sO_Y$ of the Higgs bundle $(F_{\underline{k}},\tau_{\underline{k}})$.
Remark that $F^{m,k-m}_{\underline{k}}$ decomposes in a direct sum of factors
\begin{gather}\label{eqrg.3}
\big(\bigotimes_{i=1}^{s}F_i^{m_i,k_i-m_i}\big)\otimes
\big(\bigotimes_{i=s''+1}^{s}{F'}_i^{m'_i,k'_i-m'_i}\big)\\
\mbox{with \ \ } \notag
m=\sum_{i=1}^s m_i + \sum_{i=s''+1}^s m'_i \mbox{ \ \ and \ \ }
k=\sum_{i=1}^s k_i + \sum_{i=s''+1}^s k'_i .
\end{gather}
The sheaves in (\ref{eqrg.3}) are tensor products of poly-stable sheaves.
One finds for each direct factor $\sF$ of the sheaf in (\ref{eqrg.3}) 
\begin{multline*}
\mu(\sF)=\sum_{i=1}^s \big( m_i\cdot\ch_1(\sL_i)+(k_i-m_i)\cdot(\ch_1(\sL_i)-
\Upsilon(\Omega_i))\big).\ch_1(\omega_Y(S))\\
- \sum_{i=s''+1}^s \big((k'_i-m'_i)\cdot\ch_1(\sL_i)+m'_i\cdot(\ch_1(\sL_i)-
\Upsilon(\Omega_i))\big).\ch_1(\omega_Y(S))=\\
= \sum_{i=1}^s \big( (k_i-(k_i-m_i)\frac{n_i+1}{n_i})\cdot\mu(\sL_i)\big)-
\sum_{i=s''+1}^s \big( (k'_i-m'_i\frac{n_i+1}{n_i})\cdot\mu(\sL_i)\big).
\end{multline*}
Assume that $\sF=\sO_Y$. Since $\mu(\sL_i)>0$  this implies for $i=1,\ldots,s''$ that
$0=k_i-2(k_i-m_i)$, hence $k_i=2m_i$. For $s'' < i \leq s$ one finds
$$
0=n_i k_i-(k_i-m_i)(n_i+1)-n_i k'_i+m'_i(n_i+1)= (m_i+m'_i)(n_i+1)-k_i-n_ik'_i.
$$
If in either one of those cases one has $k_i=k'_i$, then
$m_i+m'_i=k_i$ and one finds $k=2m$, as claimed in a). In general,
$$
m_i+m'_i=\frac{1}{n_i+1}(k_i+n_ik'_i)
$$
is uniquely determined by $\underline{k}$, hence $m$ as well, and one obtains b).
\end{proof}

If for a family of $g$-dimensional Abelian varieties  $f:V\to U$ the length of
$R^g f_*\C_V=\bigwedge^g R^1f_*\C_V$ is $g$, then the family is rigid (see \cite[Section 3]{VZ2}
and the references given there). This criterion will only apply, if in Theorem \ref{in.5} all
the stable direct factors of $\Omega_Y^1(\log S)$ are invertible, hence if $n_i=1$ for all $i$.

Nevertheless the rigidity holds true, even if in Theorem \ref{in.5} one has
direct factors $\Omega_i$ with $n_1>0$.

\begin{lemma}\label{rg.5}
All global endomorphism $\End(R^1f_*\Q_V)$ are pure of bidegree $(0,0)$. In particular
the family $f:V\to U$ is rigid.
\end{lemma}
\begin{proof}
By \cite{Fa} the second part follows from the first one. To keep the notations consistent with those
of the proof of \ref{rg.4} we will show that each global section $\varrho$ of
$\bigotimes^2\W$ is of pure bidegree $(1,1)$, for $\W=R^1f_*\C_V$. Assume that
there is a section of a different bidegree, let us say of bidegree $(2,0)$. Then
$\varrho$ gives rise to some trivial direct factor of
\begin{gather*}
\sF_1=(\sL_i\otimes \U_i)\otimes (\sL_j\otimes\U_j),\\
\mbox{of \ \ }
\sF_2=(\sL_i^\vee\otimes \U_i^\vee\otimes \Omega_i)\otimes (\sL_j^\vee\otimes\U_j^\vee\otimes\Omega_j),\\
\mbox{or of \ \ }
\sF_3=(\sL_i\otimes \U_i)\otimes (\sL_j^\vee\otimes\U_j^\vee\otimes\Omega_j).
\end{gather*}
One has $\Upsilon(\sF_1)=\ch_1(\sL_i)+\ch_1(\sL_j)$, 
\begin{gather*}
\Upsilon(\sF_2)=\Upsilon(\Omega_i)-\ch_1(\sL_i)+\Upsilon(\Omega_j)-\ch_1(\sL_j) 
= \frac{1}{n_i}\ch_1(\sL_i)+\frac{1}{n_j}\ch_1(\sL_j)\\
\mbox{ and \ \ }
\Upsilon(\sF_3)=\ch_1(\sL_i)+\Upsilon(\Omega_j)-\ch_1(\sL_j) 
= \ch_1(\sL_i)+\frac{1}{n_j}\ch_1(\sL_j).
\end{gather*}
Obviously, none of those poly-stable sheaves can have a trivial direct factor.
\end{proof}

\section{Shimura varieties}\label{sh}
Let $F$ be an Abelian variety and let $Q$ be the polarization, i.e. an
antisymmetric non-degenerate form on $H^1(F,\Q)$. The
Hodge group $\Hg(F)$ is defined in \cite{Mum1} as the smallest $\Q$-algebraic
subgroup of $\Sp(H^1(F,\Q),Q)$, whose extension to $\R$ contains the complex structure
$$
u:S^1 \>>> \Sp(H^1(F,\R),Q).
$$
(see also \cite{Mum2}), where
$z$ acts on $(p,q)$ cycles by multiplication with $z^p\cdot \bar z^q$.

In a similar way, one defines the Mumford-Tate group $\MT(F)$.
The complex structure $u$ extends to a morphism of real algebraic groups
$$
h:{\rm Res}_{\C/\R}\G_m \>>>\Gl(H^1(F,\R)),
$$
and $\MT(F)$ is the smallest $\Q$-algebraic subgroup of $\Gl(H^1(F,\Q))$, whose extension to $\R$ contains
the image (see \cite{De}, \cite{De2}, \cite{Mo} and \cite{Scho}). Let us recall some of its properties, stated in \cite{Mo} and \cite{De} with the necessary references. The group $\MT(F)$ is reductive, and it preserves the intersection form $Q$ up to scalar multiplication.

Equivalently $\MT(F)$ is the largest $\Q$-algebraic subgroup of the linear group $\Gl(H^1(F,\Q))$,
which leaves all Hodge cycles of $F\times \cdots \times F$
invariant, hence all elements
$$
\eta\in H^{2p}(F\times \cdots \times F,\Q)^{p,p}
=\big[\bigwedge^{2p}(H^1(F,\Q)\oplus \cdots \oplus
H^1(F,\Q))\big]^{p,p}.
$$
For a smooth family of Abelian varieties $f:V \to U$ there exist a union $\Sigma$ of countably many
proper closed subvarieties of $Y$ such that $\MT(f^{-1}(y))$ is independent of $y$ for
$y\in U\setminus \Sigma$ (see \cite{De2}, \cite{Mo} or \cite{Scho}). Let us fix such a very general
point $y\in U\setminus \Sigma$ in the sequel and $F=f^{-1}(y)$. Then the Mumford-Tate group $\MT$ of $R^1f_*\Q_{V}$ is the Mumford-Tate group of $F$.

Consider Hodge cycles $\eta$ on $F$ which remain Hodge cycles under
parallel transform. Then $\MT$ is the largest $\Q$-subgroup of $\Gl(H^1(F,\Q))$
which leaves all those Hodge cycles invariant (\cite[\S 7]{De2} or \cite[2.2]{Scho}).

Let $\Mon^0$ be the algebraic monodromy group, i.e. the connected component of the Zariski
closure of the image of the monodromy representation. Let us recall two results from
\cite{De2} and \cite{An} (see \cite[1.4]{Mo}).

\begin{proposition}\label{sh.1} \
\begin{enumerate}
\item[a.] $\Mon^0$ is a normal subgroup of the derived subgroup $\MT^{\rm der}$
of $\MT$.
\item[b.] If for some $y'\in Y$ the fibre $f^{-1}(y')$ has complex multiplication,
then $\Mon^0=\MT^{\rm der}$.
\end{enumerate}
\end{proposition}
\begin{lemma}\label{sh.2} Let $K$ be a totally real Galois extension of $\Q$ with Galois group $\Gamma$, and let $R^1f_*K_V=\W_1\oplus\cdots \oplus \W_\nu$ 
be the decomposition in irreducible $K$-sub-variations of Hodge structures. Then
$\MT\otimes K$ is conjugate to a subgroup of $\Gl(W_1)\times\cdots\times\Gl(W_\nu)$
where $W_i$ denotes the fibre of $\W_i$ over $y$.
\end{lemma}
\begin{proof}
Since the decomposition $\W_1 \oplus \cdots \oplus \W_\nu$ is defined over a subfield of $\R$, one can decompose the complex structure correspondingly
as a sum of $h_i:{\rm Res}_{\C/\R}\G_m \to \Gl(W_i\otimes_K\R)$. Then
$\Gl(W_1)\times\cdots\times\Gl(W_\nu)$, extended to $\R$, contains the image of $h$.
\end{proof}
\begin{proposition}\label{sh.3}
If $f:V\to U$ satisfies the assumptions made in Theorem \ref{in.5} one has
$\Mon^0=\MT^{\rm der}$.
\end{proposition}
\begin{proof}
Recall that $\MT=\MT(F)$ for a very general fibre $F$ of $f$.
By \cite[4.4]{Sim2} $\Mon^0$ is reductive. By \cite[3.1 (c)]{De} it is sufficient to show that each tensor
$$
\eta\in \bigwedge^{k}\big(H^1(F,\Q)\oplus \cdots \oplus
H^1(F,\Q)\big)= H^{k}(F\times \cdots \times F, \Q)
$$
which is invariant under $\Mon^0$ is also invariant under $\MT^{\rm der}$.
By abuse of notations, let us replace $F\times \cdots \times F$ by $F$.

So we will consider sections $H^{k}(F, \Q)$. Since each section which is invariant under
$\Mon^0$ is a sum of global sections
$$
\eta\in H^0\big(U, \bigwedge^k \W_\Q\big),
$$
for $\Q$-irreducible local sub-system $\W_\Q\subset R^1f_*\Q_V$, it is sufficient to
show that such $\eta$ are invariant under $\MT^{\rm der}$.

Let $L$ be a Galois extension of $\Q$ with Galois group $\Gamma$, containing all the fields
$K(\sqrt{a_i})$ constructed in \ref{rg.3}.
Over $L$ the section $\eta$ decomposes as
$$
\eta=\sum_{I} \eta_{\underline k}
$$
where $I$ is a set of tuples $\underline{k}=(k_1,\ldots , k_{s''}, k_{s''+1}, k'_{s''+1}, \ldots , k_s, k'_s)$, satisfying the equation (\ref{eqrg.2}), and where
$$
\eta_{\underline k}\in\W_{\underline{k}}=\big(\bigotimes_{i=1}^{s''}\bigwedge^{k_i}\V_i\big)\otimes
\big(\bigotimes_{i=s''+1}^{s}\bigwedge^{k_i}\V_i\otimes \bigwedge^{k'_i} \V_i^\vee\big).
$$
If in the decomposition (\ref{eqrg.1}) one has $s=s''$,
or more generally if $k_i=k_i'$ for $i=s''+1,\cdots,s$ and for all $\underline{k}$ with
$\eta_{\underline k}\neq 0$, then by \ref{rg.4} the section $\eta$ has bidegree $(p(\eta),p(\eta))$, hence it is a Hodge cycle and invariant under $\MT^{\rm der}$.

Otherwise choose some ${\underline k}^{(0)}$, say with $k^{(0)}_s\neq {k'}^{(0)}_s$, and
$\eta_{\underline{k}^{(0)}}\neq 0$. Consider the fix group $\Gamma'$ of $\W_{\underline{k}^{(0)}}$.
Replacing $\eta_{\underline{k}^{(0)}}$ by the sum over its conjugates under $\Gamma$, we may
assume that $\Gamma'$ is also the fix group of $\eta_{\underline{k}^{(0)}}$.

Let $\eta'$ be the sum over all different conjugates of $\eta_{\underline{k}^{(0)}}$
under the action of $\Gamma$, then $\eta=\eta'+\eta''$, with $\eta'$ and $\eta''$ defined over $\Q$.

Hence it is sufficient to consider $\eta=\eta'$, hence to assume that
$\eta$ is equal to the sum over all different conjugates of $\eta_{\underline{k}^{(0)}}$. Choosing the index set $I$ to be minimal, one has
$$
I\simeq \{ \eta_{\underline{k}^{(0)}}, \ldots , \eta_{\underline{k}^{(\nu)}}\},
$$
hence a transitive action of $\Gamma$ on $I$.  We write this action as
$\underline{k}\mapsto \gamma(\underline{k})$. So for each $\iota$ there is some $\gamma_\iota \in \Gamma$ with $\eta_{\underline{k}^{(\iota)}}=\gamma_\iota(\eta_{\underline{k}^{(0)}})$.

The section $\eta$ gives rise to
$$
\sigma'=\bigwedge_{\iota=1}^{\nu}\eta_{\underline{k}^{(\iota)}}.
$$
The Galois group $\Gamma$ permutes the different components of $\eta$ and
${\sigma'}^\gamma=\pm \sigma'$. This defines homomorphism $\chi:\Gamma \to \{\pm 1 \}$. Choose a generator $\beta$ of the Galois extension of $\Q$, defined by this homomorphism, such that
$\Gamma$ acts on $\beta$ by multiplication with $\chi$. 
Then $\sigma=\beta\cdot\sigma'$ is invariant under $\Gamma$.

By \ref{rg.4} each $\eta_{\underline{k}}$ is concentrated in a unique bidegree $(p(\eta_{\underline{k}}),
q(\eta_{\underline{k}}))$. Posing the conditions
$$
p(\eta_{\underline{k}})- q(\eta_{\underline{k}}) <0 \mbox{ \ \ or \ \ }
p(\eta_{\underline{k}})- q(\eta_{\underline{k}}) >0
$$
defines two disjoint subsets $I^+$ and $I^-$ of $I$ of the same cardinality.
If $\eta_{\underline{k}}\in I^+$, then its complex conjugate lies in $I^-$, and vice versa.
So \ref{rg.4} implies that the sum over all $p(\eta_{\underline{k}})$ with
$\eta_{\underline{k}}\in I^+$ coincides with the sum over all
$q(\eta_{\underline{k}})$ with $\eta_{\underline{k}}\in I^-$.  

Then $\sigma$ is pure of bidegree $(p,p)$ for some $p$. 
Finally remark that $\sigma$ is again a section of some tensor bundle,
hence a Hodge cycle and invariant under $\MT$.

Let $\T_L$ denote the Weil restriction of the one dimensional subspace $\langle \eta_{\underline{k}^{(0)}}\rangle_L$ of $\bigwedge^k\W_L$. So
$\T_L$ is generated by the $\eta_{\gamma(\underline{k})}$ for
$\gamma\in \Gamma$, and a basis of $\T_L$ is given by the sections 
$\eta_{\underline{k}^{(0)}}, \ldots , \eta_{\underline{k}^{(\nu)}}$.
The group $\MT$ leaves $\sigma$ invariant, hence the subspace $\T_L$ as well.

By Lemma \ref{rg.3} and Lemma \ref{sh.2} the group $\MT$ respects the decomposition
of $R^1f_*K_V$ in irreducible $K$-sub-variation of Hodge structures. 
Hence if one considers the decomposition in $L$ irreducible direct factors,
it can for each $i$ only permute $\V_i$ and $\V_i^\vee$. Since $\MT^{\rm der}$ lies in the kernel
of the corresponding morphism $\MT\to \{\pm 1\}^{\oplus s}$ the group $\MT^{\rm der}$ respects
all the $\V_i$ and $\V_i^\vee$, hence $\W_{\underline{k}^{(0)}}$. 

On the other hand, since the fix groups of
$\W_{\underline{k}^{(0)}}$ and of $\eta_{\underline{k}^{(0)}}$ in $\Gamma$ coincide, the intersection
$W_{\underline{k}^{(0)}}\cap \T_L$
is generated by $\eta_{\underline{k}^{(0)}}$. In particular one finds for $h\in \MT^{\rm der}$ some
$\alpha^{(0)}(h) \in L$ with $h(\eta_{\underline{k}^{(0)}})= \alpha^{(0)}(h) \cdot \eta_{\underline{k}^{(0)}}$. Then
$$
h(\eta)=\sum_{\iota=0}^{\nu} \alpha^{(\iota)}(h) \cdot \eta_{\underline{k}^{(\iota)}},
$$
where $\gamma(\alpha^{(0)}(h))=\alpha^{(\iota)}(h)$ if and only if
$\gamma(\underline{k}^{(0)})=\underline{k}^{(\iota)}$. Obviously, for $g,h\in \MT^{\rm der}$,
$$
\alpha^{(0)}(g\circ h)=\alpha^{(0)}(g)\circ \alpha^{(0)}(h),
$$
and one obtains a homomorphism
$$
\Psi:\MT^{\rm der} \>>> L^*\times \cdots \times L^*.
$$
By definition of $\MT^{\rm der}$ such a morphism must be trivial, and
$$
\eta=\eta_{\underline{k}^{(0)}}+ \cdots + \eta_{\underline{k}^{(\nu)}}
$$
is invariant under $\MT^{\rm der}$, as claimed.
\end{proof}
As mentioned in Remark \ref{av.4}, the poly-stability of $E^{1,0}$ and $E^{0,1}$ for all
direct factors of the variation of Hodge structures allows us to apply B. Moonen's characterization of
bounded symmetric domains in $\sA_g$ \cite{Mo}. There one uses the existence of at least
one CM point in $U$. Then \ref{sh.1}, b), would imply the equality $\Mon^0= \MT^{\rm der}$. 

After we established such an equality by different arguments, one can use \cite[3.8]{Mo}
to deduce that $U$ is a Shimura subvariety of the moduli stack $\sA_g$. 
Let us sketch the argument, using a slightly different language.

The Hodge group $\Hg=\Hg(F)$ is contained in 
$$\MT\cap \Sp(H^1(F,\Q_F),Q).$$ 
By \ref{sh.2} the induced real group $\Hg_\R$ is conjugate to a subgroup of
$$
\Sp(\hat V_1, Q) \times \cdots \times \Sp(\hat V_{s''},Q) \times \Sp(\hat V_{s''+1},Q) \times \cdots \times
\Sp(\hat V_s,Q),
$$
where again $\hat V_i$ is the fibre of $\hat\V_i$ at $y\in U$.

Mumford constructs in \cite{Mum1} a Shimura variety $\sX(\Hg,u)$, as the image of
$$
\Phi: \Hg_\R \>>> \Sp(H^1(F,\Q),Q)/\big( \begin{array}{c} \mbox{centralizer of the} \\ 
\mbox{complex structure } u \end{array}\big) =\tilde\sA_g \>>> \sA_g.
$$
The morphism $\Phi$ factors through the quotient of $\Hg_\R$ by a maximal compact subgroup.

The monodromy group $\Mon^0$ is contained in $\Hg$, hence equal to
$\Hg^{\rm der}=\MT^{\rm der}$. We may replace in Mumford's construction $\Hg$ by 
the isogenous group $\Hg^{\rm der}$, and the dimension of $\sX(\Hg,u)$ is the dimension of $\Phi(\Hg^{\rm der})$.
\begin{corollary}\label{sh.4}
$\displaystyle \dim(\sX(\Hg,u)) \leq \sum_{i=1}^s n_i$.
\end{corollary}
\begin{proof}
The variation of Hodge structures comes from a representation with values in the real group
\begin{multline*}
G=({\rm U}(1,1)\times {\rm U}(\ell_{1}))\times \cdots \times ( {\rm U}(1,1)\times {\rm U}(\ell_{s''}))\\
\times ({\rm U}(n_{s''+1},1)\times {\rm U}(\ell_{s''+1})) \times \cdots \times ({\rm U}(n_s,1)\times 
{\rm U}(\ell_s)).
\end{multline*}
Since $\Mon^0=\Hg^{\rm der}$ is contained in $G$, the image of $\Phi$ lies in the quotient of
$G$ by a maximal compact subgroup. Since ${\rm U}(\ell_i)$ is compact, the latter is isogenous to
$$
{\rm U}(1,1)/K_1 \times \cdots \times {\rm U}(1,1)/K_{s''}\times {\rm U}(n_{s''+1},1)/K_{s''+1}
\times \cdots \times {\rm U}(n_s,1)/K_s.
$$
for maximal compact $K_i\subset {\rm U}(n_{i},1)$. However, ${\rm U}(n_{i},1)/K_{i}$ is a ball
quotient of dimension $n_i$.
\end{proof}
\begin{proof}[Proof of Theorem \ref{in.6}]
By Proposition \ref{in.4} we know that for each irreducible $\C$-sub-variation of Hodge structures
the conditions a)--c), stated there, hold true. So we can apply Theorem \ref{in.5} and deduce that
after replacing $Y$ by an \'etale covering, all such $\V$ are of the form $\U_i\otimes \L_i$
with $\U_i$ unitary and $\L_i$ as in Theorem \ref{in.5}, b). 

The structure of the irreducible components of the variation of Hodge structures was used in Section \ref{rg} and in this Section to show the Corollary \ref{sh.4}. The rigidity has been verified in Lemma \ref{rg.5}. So it remains to show, that Corollary \ref{sh.4} together with Lemma \ref{av.2}
implies that $U$ is a Shimura variety.
 
Since $\sX(\Hg,u)$ is a moduli variety for Abelian varieties with Hodge group contained in
$\Hg$, the morphism $\varphi:U\to \sA_g$ factors through
$$
U \>\varphi' >> \sX(\Hg,u) \>\subset >> \sA_g.
$$
By assumption, $\varphi$ is generically finite, hence \ref{sh.4} implies that $\varphi'$ is dominant.
On the other hand, writing $F^{1,0}\oplus F^{0,1}$ for the Higgs bundle of $R^1f_*\C_V$,
we had seen in Lemma \ref{av.2} that the inclusion 
$$T^1_U \>>> \varphi^*T^1_{\sA_g}=S^2(F^{0,1})|_U$$
splits. So $T^1_U \to {\varphi'}^*T^1_{\sX(\Hg,u)}$ splits as well, and $\varphi$ is \'etale.
\end{proof}


\begin{thebibliography}{X-X00} 
\bibitem[A92]{An} Andr\'e, Y.: Mumford-Tate groups of mixed Hodge structures and the theorem of 
the fixed part.  Compositio Math. {\bf 82} (1992) 1--24,
MR1154159, Zbl 0770.14003.
\bibitem[B-L92]{BL} Birkenhage, Ch., Lange, H.: Complex Abelian Varieties.
Grundlehren d. math. Wiss. {\bf 302} (1992) Springer-Verlag,
Berlin-Heidelberg,
MR1217487, Zbl 0779.14012.
\bibitem[D72]{De2} Deligne, P.: La conjecture de Weil pour les surfaces $K3$.
Invent. math. {\bf 15} (1972) 206--226,
MR0296076, Zbl 0219.14022.
\bibitem[D82]{De} Deligne, P.: Hodge cycles on abelian varieties.
(Notes by J. S. Milne). Springer Lecture Notes in Math. {\bf 900}
(1982) 9--100,
MR0654325,  Zbl 0537.14006.
\bibitem[D-M86]{DM} Deligne, P., Mostow, G.:Monodromy of hypergeometric functions and non-lattice
integral monodromy. IHES {\bf 63} (1986) 5--90,
MR0849651, Zbl 0615.22008.
\bibitem[Do87]{D} Donaldson, S.K.: Infinite determinants, stable bundles and curvature.
Duke Math. J. {\bf 54} (1987) 231--247,
MR0885784, Zbl 0627.53052.
\bibitem[E-V92]{EV} Esnault, H., Viehweg, E.: Lectures on Vanishing
Theorems. DMV-Seminar {\bf 20} (1992), Birkh\"auser,
Basel-Boston-Berlin, MR1193913, Zbl 0779.14003.
\bibitem[F83]{Fa} Faltings, G. Arakelov's theorem for abelian varieties.
Invent. Math. {\bf 73} (1983) 337--347,
MR0718934, Zbl 0588.14025.
\bibitem[F-H91]{FH} Fulton, W., Harris, J.:
Representation Theory. A first course. Graduate Texts in Math.
{\bf 129} (1991) Springer-Verlag, New-York,
MR1153249, Zbl 0744.22001.
\bibitem[G02]{Gr} Granath, H.: On Compact Shimura Surfaces. preprint (2002), 15 pages, math.AG/0211290
\bibitem[H-L97]{HL} Huybrechts, D., Lehn, M.: The Geometry of Moduli Spaces of Sheaves.
Aspects of Math. {\bf E31} Friedr. Vieweg \& Sohn, Braunschweig, 1997,
MR1450870, Zbl 0872.14002.
\bibitem[L04]{La} Langer, A.:
Semistable sheaves in positive characteristic. Ann. of Math. {\bf 159} (2004) 251--276,
MR2051393, Zbl 1080.14014.
\bibitem[Lz04]{L} Lazarsfeld, R.: Positivity in Algebraic Geometry.\\
Part I. Ergebnisse der Math. 3. Folge. {\bf 48} (2004),
Springer, Berlin-Heidelberg-New York,
MR2095471, Zbl 1066.14021.\\
Part II. Ergebnisse der Math. 3. Folge. {\bf 49} (2004),
Springer, Berlin-Heidelberg-New York,
MR2095472, Zbl pre02134815.
\bibitem[Lo03]{Lo} Looijenga, E.: Compactifications defined by arrangements. I. The ball quotient case. Duke Math. J.  {\bf 118} (2003) 151--187,
MR1978885, Zbl 1052.14036. 
\bibitem[Mo98]{Mo} Moonen, B.: Linearity properties of Shimura varieties.
Part I. J. Algebraic Geom. {\bf 7} (1998) 539--567,
MR1618140,  Zbl 0956.14016.
\bibitem[M66]{Mum1} Mumford, D.: Families of Abelian varieties.
Proc. Sympos. Pure Math. {\bf 9} (1966) 347--351,
MR0206003, Zbl 0199.24601.
\bibitem[M69]{Mum2} Mumford, D.: A note of Shimura's paper: Discontinuous
groups and Abelian varietes.  Math. Ann. {\bf 181}  (1969) 345--351,
MR0248146, Zbl 0169.23301.
\bibitem[Sc96]{Scho} Schoen, C.: Varieties dominated by product varieties.
Int. J. Math. {\bf 7} (1996) 541--571, MR1408839, Zbl 0907.14002.
\bibitem[Sh78]{Sh} Shavel, I.H.: A class of algebraic surfaces of general type constructed from quaternion algebras. Pacific J. Math. {\bf 76} (1978) 221--245,
MR0572981, Zbl 0422.14022.
\bibitem[S88]{Sim} Simpson, C.: Constructing variations of Hodge structure
using Yang-Mills theory and applications to uniformization.
Journal of the AMS {\bf 1} (1988) 867--918, MR0944577, Zbl 0669.58008.
\bibitem[S90]{Sim1} Simpson, C.: Harmonic bundles on noncompact curves.
Journal of the AMS {\bf 3} (1990) 713--770, MR1040197, Zbl 0713.58012.
\bibitem[S92]{Sim2} Simpson, C.: Higgs bundles and local systems. Publ. Math. I.H.E.S
{\bf 75} (1992) 5--95,
MR1179076, Zbl 0814.32003.
\bibitem[S93]{Sim3} Simpson, C.: Some families of local systems over smooth projective varieties.  Ann. of Math. {\bf 138} (1993) 337--425,
MR1240576, Zbl 0813.14014.
\bibitem[U-Y86]{UY} Uhlenbeck, K., Yau, S.-T.: On the existence of Hermitian-Yang-Mills connections
in stable vector bundles. Comm. Pure Appl. Math. {\bf 39} (1986) 257--293,
MR0861491, Zbl 0615.58045.
\bibitem[V-Z04]{VZ1} Viehweg, E., Zuo, K.: A characterization of certain Shimura curves
in the moduli stack of abelian varieties. J. Diff. Geom. {\bf 66} (2004) 233--287,
MR2106125, Zbl 1078.11043.
\bibitem[V-Z05]{VZ2} Viehweg, E., Zuo, K.: Complex multiplication, Griffiths-Yukawa couplings, and rigidy for families of hypersurfaces. J. Alg. Geom. {\bf 14} (2005) 481--528,
MR2129008, Zbl pre02204424.
\bibitem[Y88]{Y2} Yau, S.T.: Uniformization of geometric structures. Proc. of Symp.
in Pure Math. {\bf 48} (1988) 265--274,
MR0974340, Zbl 0656.53060.
\bibitem[Y93]{Y} Yau, S.T.: A splitting theorem and an algebraic geometric characterization of
locally Hermitian symmetric spaces. Comm. in Analysis and Geometry. {\bf 1} (1993) 473--486,
MR1266476, Zbl 0842.53035.
\bibitem[Z00]{Zu}  Zuo, K.: On the negativity of kernels of Kodaira-Spencer maps on Hodge bundles and applications. Asian J. Math. {\bf 4} (2000) 279--301,
MR1803724, Zbl 0983.32020.
\end{thebibliography}
\end{document}